\documentclass{article}

\usepackage[margin=0.8in]{geometry}
\usepackage[utf8]{inputenc}
\usepackage{graphicx}
\usepackage{amsmath,amsthm}
\usepackage{amsfonts}
\usepackage{amssymb}
\newtheorem{proposition}{Proposition}[section]

\usepackage[title]{appendix}
\usepackage{bm}
\usepackage{tikz}
\usepackage{eurosym}
\usepackage{xcolor}
\usepackage{hyperref}
\usepackage{alphalph}
\usepackage{etoolbox}
\usepackage{cleveref}
\crefname{appsec}{Appendix}{Appendices}
\usepackage{url}
\usepackage{booktabs}
\usepackage{subcaption,caption}
\usepackage{multirow}
\usepackage{longtable}
\usepackage{algorithmic}
\usepackage{algorithm}
\usepackage{mdframed}
\newmdtheoremenv{prop}{Proposition}

\newcommand{\mtc}{\mathcal}

\AtBeginDocument{%
  \AtBeginEnvironment{subequations}{%
  }
}
\allowdisplaybreaks

\date{}

\begin{document}
\title{Exact solutions to a carsharing pricing and relocation problem under uncertainty}

\author{Giovanni Pantuso \\
  Department of Mathematical Sciences, University of Copenhagen,\\
  Universitetsparken 5, 2100, Copenhagen, Denmark\\
  gp@math.ku.dk.}
\maketitle              

\begin{abstract}
  In this article we study the problem of jointly deciding carsharing prices and vehicle relocations. We consider carsharing services operating in the context of multi-modal urban transportation systems.
  Pricing decisions take into account the availability of alternative transport modes, and customer preferences with respect to these. In order to account for the inherent uncertainty in customer preferences, the problem is formulated as a mixed-integer two-stage stochastic program with integer decision variables at both stages. We propose an exact solution method for the problem based on the integer L-Shaped method which exploits an efficient exact algorithm for the solution of the subproblems. Tests on artificial instances based on the city of Milan illustrate that the method can solve, or find good solutions to, moderately sized instances for which a commercial solver fails. Furthermore, our results suggest that, by adjusting prices between different zones of the city, the operator can attract significantly more demand than with a fixed pricing scheme and that such a pricing scheme, coupled with a sufficiently large fleet, significantly reduces the relevance of staff-based relocations. A number of issues, that remain to be addressed in future research, are pointed out in our conclusions. 
\end{abstract}%


\section{Introduction}\label{sec:intro}
High dependency on private vehicles and low occupancy rates
increase car usage and congestion in many cities of the world, contributing to pollution and poor urban air quality \cite{ChoYAB18}.
Improvements of public transport \cite{PruKLF12} and road pricing measures
\cite{CarG11,Dud13} have, to a large extent, failed to provide sustainable solutions \cite{ChoYAB18,BarFDS12,MayN96}.
In this context, shared mobility, and particularly carsharing, has emerged as a viable alternative, linked to, e.g.,
a decrease in congestion \cite{CraEHN12}, pollution \cite{MarS16}, land used \cite{ShaC13} and transport costs \cite{Dun11,MouPC19}.

In reaction to the high flexibility demanded by users, modern carsharing services are commonly designed for \textit{on-demand}, \textit{short-term}, \textit{one-way} usage \cite{IllH19}.
That is, users are allowed to rent a car without reservation and return it as soon as, and wherever (within the operating area), their journey is completed.
The drop-off location/station is thus typically different from the pick-up location/station.
Such configuration poses new planning challenges to \textit{carsharing operators} (CSOs).
On-demand rentals make the CSO unaware of when, where and for how long new rentals will occur.
One-way rentals create frequent imbalances in the distribution of vehicles, that is an accumulation of vehicles in low-demand zones, and vehicle shortage in high-demand zones \cite{BarTX04,BoyZG15} with levels of service dropping accordingly. A central task for a CSO is to provide a distribution of vehicles in the business area compatible with demand tides and oscillations \cite{WasJB13,WeiB15}.

As a prime form of response to these challenges, CSOs initiate staff-based vehicle relocations between stations/zones of the city before shortages occur and customer satisfaction levels drop \cite{JorC13,IllH19}. That is, CSO's staff reach designated cars and drives them to different places.
This gives rise to the so called \textit{Vehicle Relocation Problem} (VReP), which consists of determining the relocations to perform in order to prepare for future demand. The research literature covers several problem settings, different levels of detail and granularity of decisions, as well as different mathematical approaches, see e.g., \cite{BoyZG17,WeiB15,BoyZG15,BruCL14,JorCB14,KekCMF09,KekCC06,MukW05,BarTX04,BarT99,FolHFAP20,HelJMAFP21}. These studies are thoroughly reviewed and classified in a recent survey, see \cite{IllH19}.

The survey also reports that staff-based rebalancing could be complemented by manipulating demand through dynamic pricing.
In fact, users of urban mobility services typically choose among different transport modes (e.g., metro, carsharing, bikesharing) that vary in a number of key attributes including price,
see e.g., \cite{ZoeK16,HanP18, Pan20}. However, the authors of the survey comment that, at present times, this issue is solely ``an important part of future research'' \cite{IllH19}.

Compared with relocation decisions, carsharing pricing strategies have received limited attention in the research literature though the number of available studies is growing. 
We can distinguish two main categories of pricing strategies, which we refer to as \textit{individual} and \textit{collective}, according to their end recipient.
Individual pricing strategies are targeted to individual users. They require an interaction between the CSO and the individual user by means of which the final trip price, pick-up or drop-off location are agreed upon. As an example, the operator sends the individual user offers in the form of discounts or bonuses in exchange to a trip configuration which the operator deems beneficial for the entire system. Collective pricing strategies are, instead, targeted to the entire user base.
They have the scope of influencing the cumulative rental demand by, e.g., decreasing the price of rentals to/from selected zones,
but do not require an interaction with the individual user, nor their reply as to whether the price is accepted or not. The approach proposed in this article belongs to the latter category.

Several pricing strategies can be classified as \textit{individual}.
In \cite{WagWBN15} a method is developed that identifies vehicles placed in low demand zones using idle time as a proxy.
The method then offers the user a drop-off location with low expected idle time in exchange for a discount (e.g., free minutes of usage).
The method is evaluated in a simulation framework based on the city of Vancouver. The authors report that the average vehicle idle time is decreased by up to $16$ percent.
In \cite{WasJ16} it is assumed that each customer is interested in a trip between a specific origin and destination and is sensitive to price.
The operator then offers a price for the given trip in order to, ideally, incentivize/prevent favorable/unfavorable car movements.
They model the carsharing system as a continuous-time Markov chain where a pricing policy is input to the model.
In \cite{DifSS18} a user-based relocation method is presented in which the users are offered to leave the car
in a location different from the one planned in exchange for a fare discount.
The authors formulate the decision problem as mixed-integer nonlinear programming problem,
and model customers preferences with respect to the offer of alternative drop-off locations expressing the corresponding utility as a functions of the distance between the desired and offered drop-off locations.
In \cite{StoG21} a predictive, user-based, relocation strategy is introduced for station-based carsharing services.
They assume that the CSO can offer each returning customers an incentive to relocate.
That is, upon the arrival of a customer, the CSO determines whether to offer an incentive, what the incentive should be and where the vehicle should be relocated.
Therefore, an optimization problem is solved upon the arrival of every customer. Estimated customer preferences are used to model their reaction to incentive offers.
Station-based services are considered also in \cite{LiuXC21} where users send real-time trip requests to the CSO, specifying their origin and destination.
Upon receiving requests, the CSO assigns vehicles to users, plans staff-based relocations and determines incentives for customers, to whom a car has already been assigned, in exchange
for a change of destination.
The focus on \cite{WuLSP21} is instead on free-floating services.
The system they consider is organized as follows. Each user sends a request for a vehicle either on-demand or as a reservation for a future rental.
With the request they are required to specify intended pick-up and drop-off locations and departure time.
The system elaborates the available requests and responds to users with proposed service options.
Each service option includes a pick-up and drop-off location, pick-up time, and price.
Finally, in \cite{WanJL21} a pricing scheme to induce user-based relocations is introduced, with focus on station-based one-way services.
The pricing scheme consists of adding or deducting a fixed expense to the
original expense to adjust the users' preferential pick-up and drop-off location.
To study the user behavior a questionnaire is carried out. The pricing problem is finally
formulated as an optimization problem.

A number of articles have also focused on \textit{collective} pricing strategies. 
A mixed-integer nonlinear programming model is provided in \cite{XuML18} for the joint problem of deciding fleet size,
trip pricing and staff-based relocations for an electric station-based carsharing service.
The authors consider demand elasticity with respect to prices using a logit-based function. The authors test the method on a case study based in Singapore. 

In \cite{WanM19} a pricing scheme is developed with the scope of influencing user demand and keeping the distribution of vehicle at a given balance level.
After modeling the relationship between carsharing price and demand, the authors develop a nonlinear optimization model to define pricing schemes which minimize the deviation from inventory upper and lower bounds at each charging station.

A station-based electric one-way service is considered in \cite{XieWWDM19}.
The CSO is to decide charging schedules and service prices.
The authors assume that rental demand is influenced by prices and adopt a linear elasticity function to express demand as a function of price.
Electric carsharing services are considered also in \cite{RenLLHL19} where a dynamic pricing scheme is proposed with the scope of solving imbalances in the
distribution of vehicles, as well as facilitating vehicle-grid-integration. For each origin and destination station the operator can
influence demand using two price adjustment levels. Rental demand is connected to prices via a price elasticity to a reference demand for
a default price. Pricing decisions are made using a mixed-integer nonlinear program.
Finally, in \cite{LuCZLL21} a bilevel nonlinear mathematical programming model is proposed to determine carsharing prices and staff-based relocations.
In the upper level, the carsharing operator determines vehicles relocations and prices.
In the lower level, travelers choose travel modes from a cost-minimization perspective and demand is computed using a logit model.
The authors assume a one-way station-based carsharing system in competition with private cars.

The method presented in the present paper can be classified as a collective pricing strategy and extends the available literature in a number of ways. 
First, the available methods do not take into account the impact of alternative transport modes on transport choices.
Carsharing services live in multi-modal transport systems and failure to model this heterogeneity may result into myopic models of customers behavior.
As an example, compared to \cite{LuCZLL21}, in the present paper customers choices are not limited to private or shared cars.
Rather, we assume that customers may choose among any number of available transport services (e.g., bicycle and bus).
This, in turn, can help the CSO set prices from/to a given zone also as a function of the alternatives available in the zone. 
Second, the articles that include demand elasticity, typically limit their attention to elasticity with respect to prices, see e.g., \cite{XuML18,LuCZLL21} and \cite{StoG21}.
In this article we allow the operator to model customers preferences with respect to any number of both exogenous and endogenous characteristics of the service such as,
but not limited to, travel time and waiting time. Such elasticity is modeled using utility functions which yield a linear optimization problem as long as
the function is linear in the endogenous characteristics of the service e.g., price.
Finally, extending all available methods, we explicitly account for uncertainty with respect to customer preferences.
That is, we consider that a portion of the preferences of each customer is unknown to the operator and is, as such, handled by means of a stochastic program.

The contributions of this article can be stated as follows.
\begin{enumerate}
\item We propose a two-stage integer stochastic programming model for the joint pricing and relocation problem.
  The central idea is to influence demand by acting on prices, and performing preventive relocations accordingly,
  in order to maximize expected profits.
  To model the interplay between pricing decisions and customers choices we follow the recipe first provided by \cite{BieS16}
  for integrating demand models within mixed integer linear optimization models.
  This framework consists of modeling user preferences, and the uncertainty therein, by means of utility functions.
  A discretization of the unknown portion of the utility, and the adoption of utility functions which are linear in the decision variables of the model,
  ensure that the resulting optimization model is linear.
  This framework has also been used by e.g., \cite{PanABG17} in the context of parking services and \cite{HanP18} in the context of carsharing,
  and a more general description is provided in \cite{PanBGA21}.
\item To solve the resulting stochastic program with integer variables at both stages,
  we propose an exact L-Shaped method that exploits a compact reformulation and efficient exact algorithm for the integer subproblems. 
\item We provide empirical evidences on the performance of the algorithm and on the solutions obtainable,
  based on artificial instances built on data from the city of Milan. An instance generator is made available online.
\end{enumerate}

The remainder of this article is organized as follows. 
In \Cref{sec:prob} we define the problem and clarify modeling assumptions.
In \Cref{sec:model} we provide an extensive formulation of the problem which has the scope of explicitly defining the relationship between pricing decisions and customer choices.
We introduce the model, particularly the second-stage, in an extensive and discursive manner in order to make the interplay between pricing decisions and customers demand explicit.
In \Cref{sec:ls} we describe an integer L-Shaped method to find exact solutions to the problem.
This is enabled by a compact and more tractable reformulation of the second-stage problem where customers choices are pre-processed,
and by an exact greedy algorithm for solving the second-stage problem, both described in \Cref{sec:ls}. 
In \Cref{sec:instances} we present a set of artificial instances based on the carsharing services offered in the Italian city of Milan. The same instances as well as an instance generator are made available online.
In \Cref{sec:results} we present the results of a computational study. We shed lights on the efficiency of the algorithm and comment on the solutions obtainable by means of the model introduced.
In \Cref{sec:conclusions} we draw final conclusions, point out existing limitations of this work and discuss possible avenues of future research.

\section{Problem definition and assumptions} \label{sec:prob}
A CSO offers one-way, reservation-free, carsharing services and is faced with the problem of jointly deciding the prices to charge and relocations to perform in order to comply with demand.
The characteristics of the service, and the perimeter of the corresponding decision problem, are clarified by following assumptions.
\begin{figure}
 \centering
 \includegraphics[width=\textwidth]{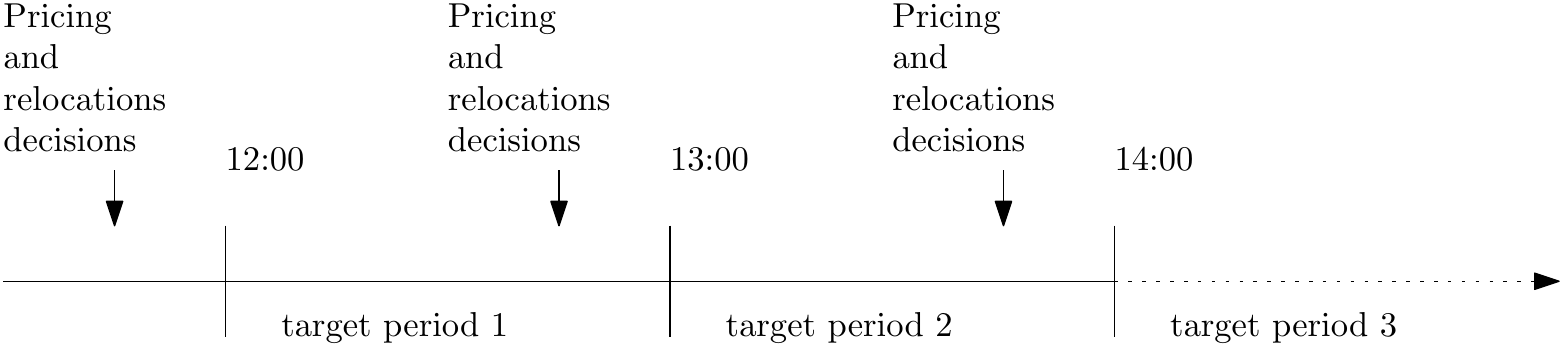}
  \caption{Target periods and decision timing. In this example target periods have a length of one hour and pricing and relocation decisions are made and implemented before the beginning of each the target period.}
  \label{fig:target_period}
\end{figure}

\begin{figure}
 \centering
 \includegraphics[width=\textwidth]{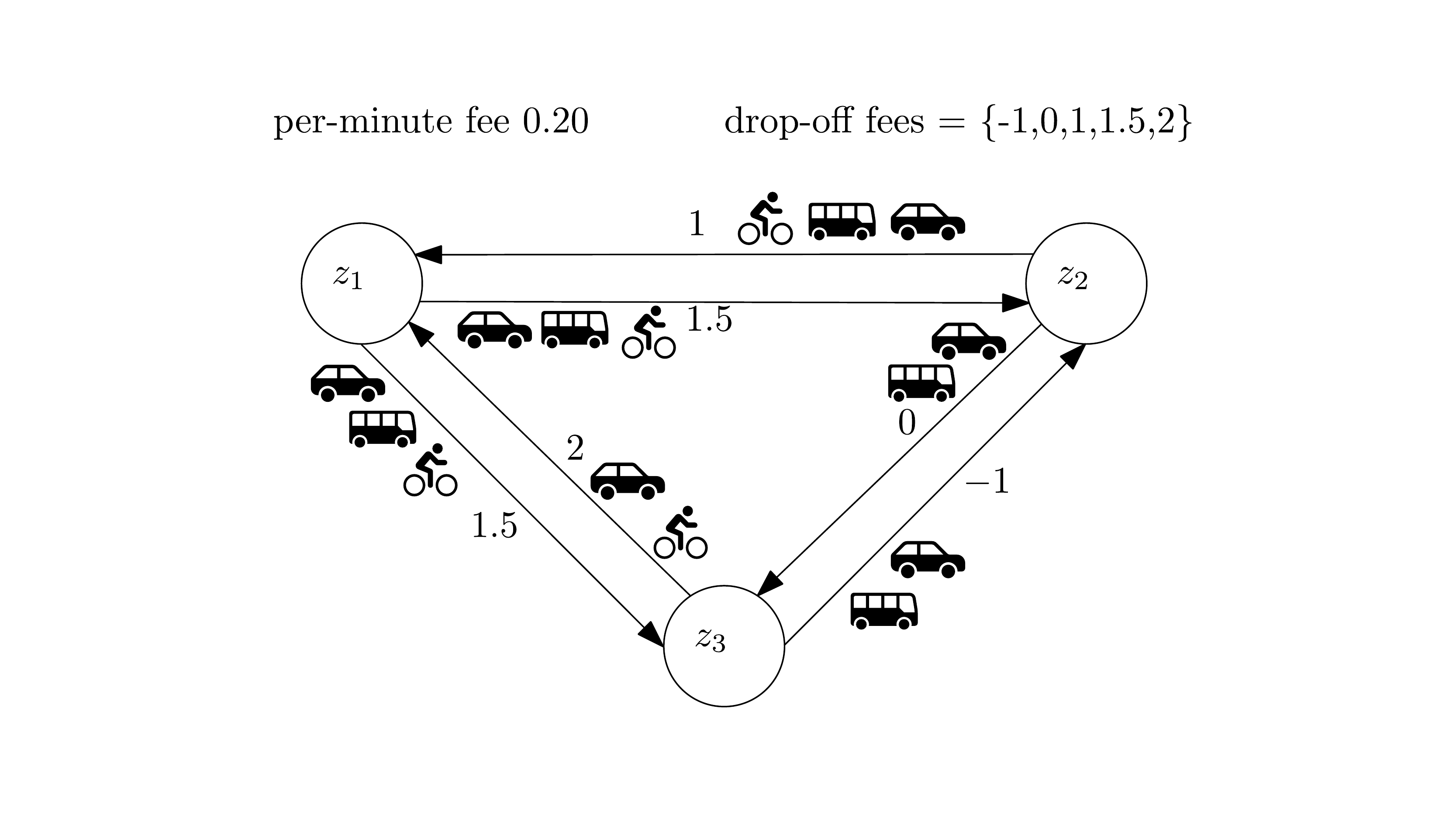}
  \caption{Zones, fees and alternative transport services. This example counts three zones, and three alternative transport services. Not all transport services are available between each pair of zones. Prices are expressed in Euro.}
  \label{fig:zones}
\end{figure}

\begin{description}
\item[A0-Target periods] The operating hours are partitioned into a number of distinct \textit{target periods},
  that is, portions of the operating hours in which the CSO may, in general, apply different prices and distributions of the fleet, see \Cref{fig:target_period}.
  Before each target period, and the CSO must decide i) the prices to apply during the target period and ii) the relocations to perform in sight of the uncertain rental demand during the target period.
  The CSO plans for each target period independently based on updated system information (e.g., fleet distribution and individual vehicles' status) and demand outlook within the target period.
  In the example in \Cref{fig:target_period}, before 12.00 the CSO must decide the prices to apply in target period $1$ (12.00 -- 13.00) and the relocations to perform before that.
  Thus, in this case the target period lasts for one hour. Similarly, one may consider longer target periods,
  e.g., morning hours and afternoon hours, as well as shorter ones, e.g., $30$ minutes target periods, depending on how often it is sensible to adjust prices in the specific context.
\item[A1-Business area] The operating area is made of a finite set of locations, henceforth \textit{zones}, see \Cref{fig:zones}.
  If the carsharing service is station-based a zone naturally represents a station. If the service is free-floating we assume that the business area is suitably partitioned into a number of zones, and each zone is represented by a suitable geographical location.    
\item[A2-Pricing scheme] The price is made of a \textit{per-minute fee} and a \textit{drop-off fee}. The per-minute fee is valid throughout the day (i.e., in all target periods) and is independent of the origin and destination of the trip.
  Instead, the drop-off fee can be different in each target period and for each origin and destination. \Cref{fig:zones} provides an example where the per-minute fee is Euro $0.2$, independently of the origin and destination, while for each pair of zones a different drop-off fee is set.
  In the example, a drop-off fee of Euro $1.5$ is charged if the car is picked up in zone $z_1$ and returned in zone $z_3$, while a drop-off fee of Euro $-1$ is charged if the car is picked up in zone $z_3$ and returned in zone $z_2$.
  Thus, we assume that the drop-off fee may also be negative to encourage desired movements of cars and increase demand.
  This setup generalizes the pricing schemes adopted in a number of carsharing services which typically charge a positive drop-off fee only if the customer returns the car in specific, unfavorable, zones of the city, or provide incentives, such as free driving minutes, to pick up a car in specific, unfavorable, zones.
  In order to keep the pricing scheme easy to communicate to customers, we assume the CSO must choose among only a finite set of possible drop-off fees. In the example of \Cref{fig:zones} this set is Euro $\{-1,0,1,1.5,2\}$.
\item[A3-Alternative transport services] The business area offers a number of alternative transport services (e.g., public transport and bicycles) outside the control of the CSO.
  The alternative services may be different for each pair of zones. Each alternative service has unlimited capacity (i.e., each customer can choose any alternative service without decreasing their availability). In the example in \Cref{fig:zones}, carsharing is offered on all origin-destination pairs,
  while busses are not an available alternative for moving from zone $z_3$ to zone $z_1$, and riding a bicycle is not an option between $z_2$ and $z_3$ due to e.g., the absence of suitable bicycle lanes.
\item[A4-Customers are informed] The CSO is able to inform customers about the current price from their location to every other zone, prior to rentals.
  In the example of \Cref{fig:zones}, the CSO is able to inform a user in $z_1$ (e.g., on the mobile application used to locate the car) that, if the car is returned in $z_2$, there will be a drop-off fee of Euro $1.5$, in addition to the per-minute fee. 
  Customers are also aware about the availability of alternative transport services. Considering the example in \Cref{fig:zones}, a customer moving between $z_1$ and $z_2$ knows that they may use bus, bicycle, and carsharing.
  For all possible transport modes (including carsharing) the user knows the respective prices and characteristics (e.g., waiting time and travel time).
\item[A5-Closed market] A customer chooses exactly one transport service among the available ones. This corresponds to saying that a customer does not give up their trip. In the example of \Cref{fig:zones}, a customer moving from $z_1$ to $z_2$ will eventually choose to travel either by bicycle, carsharing or bus, and complete its journey.
\item[A6-Customers preferences] The CSO is able to describe a portion of customers travel preferences as a function of different observable characteristics of the available transport services (e.g., travel time, price and waiting time). Nevertheless, the choice of each customer depends also on a number of additional elements not observable by the CSO. Therefore, customers preferences are partially unknown to the CSO. The unknown part of customers preferences is fully described by a probability distribution. 
\item[A7-Direct rentals] Customers traveling with shared cars drive directly from their origin to their destination zone. This assumption is made for simplicity and is without loss of generality.
  Different travel patterns can be included simply by modeling customer-specific travel times in \eqref{eq:utility}.
\item[A8-Homogeneous fleet] All shared vehicles are identical. This assumption is made for the sake of simplicity in the exposition of the reformulation of the second-stage problem and is without loss of generality.
  Throughout the text we will comment on the necessary modifications in case of a heterogeneous fleet.
\item[A9-Profit maximization] The CSO maximizes profits. While other objectives may be considered, such as maximizing demand served, or minimizing zonal deficit of cars, profits are the central objective of private carsharing operators.
\item[A10-One-way trips] For the sake of simplicity, we assume one-way trips. That is, customers move from their origin zone to a different zone.
  The model presented in \Cref{sec:model} can however accommodate also round trips, provided a suitable specification of the parameters of the trip (e.g., duration).
\end{description}

Based on these assumptions, the problem can be briefly stated as follows. 
Given (a) a target period, (b) the cumulative mobility demand between each pair of zones in the target period, (c) usage and relocation costs, (d) the current distribution of cars, (e) a model of customers preferences including a probability distribution describing customer preferences unknown to the CSO, the CSO is to decide
i) the drop-off fees to apply during the target period and ii) the relocations to perform in sight of the uncertain rental demand during the target period in order to maximize expected profits.

\section{Mathematical model}\label{sec:model}
Consider a urban area represented by a finite set $\mtc{I}$ of zones (e.g., charging stations or a suitable partition of the business area)
and a CSO offering a finite set of shared vehicles $\mtc{V}$.
Before the beginning of the target period, the CSO is to decide the drop-off fee between each pair of zones and the relocations to perform to better serve demand in the target period.
At the time of planning, the fleet is geographically dispersed in the urban area as the result of previous rentals.
Let decision variable $z_{vi}$ be equal to $1$ if vehicle $v$ is made available for rental in (possibly relocated to) zone $i$ in the target period, $0$ otherwise.
Let $C^R_{vi}$ be the relocation cost born by the CS company to make vehicle $v$ available in zone $i$. This cost is zero if the vehicle is initially in zone $i$ and positive otherwise.
Let $\mtc{L}$ be a finite set of drop-off fees the CSO may apply. Let decision variable $\lambda_{ijl}$ be equal to $1$ if fee $l$ is applied between zone $i$ and and zone $j$, $0$ otherwise.
Finally, let $z:=(z_{vi})_{i\in\mtc{I},v\in\mtc{V}}$ and $\lambda:=(\lambda_{ijl})_{i,j\in\mtc{I},l\in\mtc{L}}$.
The carsharing pricing and relocation problem is thus 
\begin{subequations}
  \label{eq:1S}
  \begin{align}
    \label{eq:1S:obj}&\max-\sum_{v\in\mtc{V}}\sum_{i\in\mtc{I}}C^R_{vi}z_{vi}+Q(z,\lambda)\\
    \label{eq:1S:c1}&\sum_{i\in\mtc{I}}z_{vi} = 1   & v\in\mtc{V}\\
    \label{eq:1S:c2}&\sum_{l\in\mtc{L}}\lambda_{ijl}=1& i,j\in\mtc{I}\\
                      &z_{vi}\in\{0,1\}              & i\in\mtc{I},v\in\mtc{V}\\
                      &\lambda_{ijl}\in\{0,1\} & i,j\in\mtc{I},l\in\mtc{L}.
  \end{align}
\end{subequations}
Constraints \eqref{eq:1S:c1} ensure that each vehicle is made available in exactly one zone.
Constraints \eqref{eq:1S:c2} state that exactly one drop-off fee can be selected between each origin $i$ and destination $j$.
The objective function \eqref{eq:1S:obj} represents the expected profit obtained in the target period.
The first term consists of the total relocation cost while the second term $Q(z,\lambda)$ represents the expected revenue from rentals as a result of pricing and relocation activities.
The meaning of $Q(x,\lambda)$ will be made explicit by the end of this section as a result of the definition of the second-stage problem which we are now introducing. 

Once relocation ($z$) and pricing ($\lambda$) decisions have been made, the CSO observes the consequent customers rentals.
The business area offers a set $\mtc{A}$ of alternative transport services, outside the control of the CSO, such as metro, busses and private bicycles.
Each service has, in general, a different price and different characteristics. 
Let decision variable $p_{vij}$ be the price of service $v\in\mtc{V}\cup\mtc{A}$ between zones $i$ and $j$.
The price of a carsharing ride between zones $i$ and $j$ is
\begin{equation}\label{eq:pV}
  p_{vij} = P^VT^{CS}_{ij} + \sum_{l\in\mtc{L}}L_{l}\lambda_{ijl}\qquad \forall v\in\mtc{V}, i,j\in\mtc{I}
\end{equation}
where parameter $P^V$ is the carsharing per-minute fee,
$T^{CS}_{ij}$ the driving time between zones $i$ and $j$ and $L_{l}$ the value of drop-off fee at level $l\in\mtc{L}$ in some currency.
Note that, in case of a heterogeneous fleet, it is simply necessary to make the per-minute fee and the driving time vehicle-dependent.
Instead, the price of alternative services is entirely exogenous, that is
\begin{equation}\label{eq:pA}
  p_{vij} = P_{vij}\qquad \forall v\in\mtc{A}, i,j\in\mtc{I}
\end{equation}
where parameter $P_{vij}$ is the price of alternative service $v\in\mtc{A}$ between $i$ and $j\in\mtc{I}$.

Let $\mtc{K}$ be the set of customers, with $\mtc{K}_{i}\subseteq\mtc{K}$ being the set of customers traveling from zone $i\in\mtc{I}$ and $\mtc{K}_{ij}\subseteq\mtc{K}_i$ the set of customers traveling from $i\in\mtc{I}$ to $j\in\mtc{I}$ in the target period.

Consider an individual customer $k$. The customer is faced with a choice among a finite number of alternative transport services that can bring them to their destination.
Using the fairly standard assumption that customers maximize their utility, we can state that each service will provide the customer a different utility, and that the customer will choose the transport
service that provides them the highest utility. This utility is known to the customer but not to the CSO.

Consider now the CSO. As we said, the CSO is not aware of the utility provided by the different services to each customer.
Rather, the CSO is aware of a number of characteristics of the different services, primarily the price, $p_{vij}$ and a some additional characteristics, say $\pi^1_{vij},\ldots,\pi^N_{vij}$ for service $v$ between $i$ and $j$ (e.g., travel time and waiting time),
as well as possibly some characteristics of the decision maker.
Based on this, the CSO can specify a function that relates these, known, characteristics to the utility obtained by the customer.
We denote this function as
$$F_k(p_{vij},\pi_{vij}^1,\ldots,\pi_{vij}^N)$$
However, there are additional elements that influence the utility that the CSO does not or cannot observe.
For this reason, the utility is better represented by
$$F_k(p_{vij},\pi_{vij}^1,\ldots,\pi_{vij}^N)+\tilde{\xi}_{kv}$$
where $\tilde{\xi}_{kv}$ is a random variable that captures the difference between the utility that the CSO is able to model and the true utility observed by the customer.
Different distributions for $\tilde{\xi}_{kv}$ will lead to different choice models. As an example, the popular Logit model is obtained when each $\tilde{\xi}_{kv}$ follows, independently, an identical extreme value distribution (Gumbel type I), and the Probit model is obtained when it follows a multivariate normal distribution.
These, and additional choice models, as well as a discussion of their fundamental assumptions and limitations are discussed in, e.g., \cite{Tra09}. See also \cite{BenB99} for an exposition related to transport choices. 

Let now $u_{ijkv}$ be a decision variable which captures the utility obtained by customer $k\in\mtc{K}$ when moving from $i$ to $j\in \mtc{I}$ using service $v\in \mtc{V}\cup\mtc{A}$.
Given a realization $\xi_{kv}$ of the random term $\tilde{\xi}_{kv}$ the utility is determined by
\begin{equation}\label{eq:utility}
  u_{ijkv} = F_k(p_{vij},\pi^1_{vij},\ldots,\pi^N_{vij})+\xi_{kv}\qquad \forall i,j\in \mtc{I},k\in\mtc{K}_{ij},v\in \mtc{V}\cup\mtc{A}
\end{equation}
Note that, since the observed characteristics of the different transport services, $\pi^1_{vij},\ldots,\pi^N_{vij}$, are given,
constraints \eqref{eq:utility} are linear if $F_k(\cdot)$ is linear in $p_{vij}$, which is instead a decision variable.

Based on the utility provided by the different transport services, customers will make their choices.
Let decision variable $w_{ijkv}$ be equal to $1$ if customer $k\in\mtc{K}_{ij}$ chooses service $v\in\mtc{V}\cup\mtc{A}$, $0$ otherwise.
A customer will choose exactly one service (see Assumption A5 in \Cref{sec:prob})
\begin{equation}\label{eq:onlyOne}
  \sum_{v\in \mathcal{V}\cup\mtc{A}}w_{ijkv} = 1 \qquad \forall i,j\in\mathcal{I},k\in \mathcal{K}_{ij}
\end{equation}

In order for a customer to choose a service, the service must be available.
Let thus binary variable $y_{ikv}$ be equal to $1$ if service $v\in\mtc{V}\cup\mtc{A}$ is offered to customer $k\in\mtc{K}_{i}$, $0$ otherwise.
Alternative services $v\in\mtc{A}$ are always offered to customers whenever they are available at all, that is
\begin{equation}\label{eq:availA}
  y_{ikv} = Y_{vi}\qquad\forall i\in \mathcal{I},k\in \mathcal{K}_{i},v\in\mathcal{A}\\
\end{equation}
where parameter $Y_{vi}$ is equal to $1$ if alternative service $v$ is available in zone $i$, $0$ otherwise.
Conversely, a shared car $v\in\mtc{V}$ may be offered to customers in zone $i$ only if it is physically available at $i$, that is
\begin{equation}\label{eq:availV}
  y_{ikv} \leq z_{iv}\qquad\forall i\in \mathcal{I},k\in \mathcal{K}_{i},v\in\mathcal{V}
\end{equation}
In addition, each car $v\in \mathcal{V}$ can be rented by only one customer. If more than one customers wish to use car $v$,
the car is taken by the first customer arriving at the car. We assume that customers are indexed according to their arrival time at the car,
i.e., customer $k$ arrives before $q$ if $k<q$. We impose that a vehicle is offered to a customer only if it is offered also to the customer arriving before them
(who perhaps did not take it), that is:
\begin{equation}\label{eq:availBefore}
  y_{ikv}\leq y_{i(k-1)v} \qquad \forall i\in \mathcal{I},k\in \mathcal{K}_{i},v\in\mathcal{V}
\end{equation}
A vehicle becomes unavailable for a customer if any customer has arrived before them and rented the car, that is:
\begin{equation}\label{eq:taken}
  z_{iv}-y_{ikv}= \sum_{j\in\mtc{I}}\sum_{q\in\mtc{K}_{ij}:q < k}w_{ijqv}\qquad \forall i\in \mtc{I},k\in\mtc{K}_{i},v\in\mtc{V}
\end{equation}
that is, if car $v$ is in zone $i$ ($z_{iv}=1$), but it is not offered to customer $k$ ($y_{ikv}=0$)
we obtain
$$1 = \sum_{j\in\mtc{I}}\sum_{q\in\mtc{K}_{ij}:q < k}w_{ijqv}$$
meaning that one customer has arrived before $k$ and rented the car. On the other hand,
if the car is offered to customer $k$, ($y_{ikv}=1$), then it must be in zone $i$ ($z_{iv}=1$ -- see \eqref{eq:availV}),
and we obtain
$$0 = \sum_{j\in\mtc{I}}\sum_{q\in\mtc{K}_{ij}:q < k}w_{ijqv}$$
meaning that no customer arriving before $k$ has taken the car.
The same equality holds if the vehicle is not available at all ($z_{iv}=0$ and $y_{ijkv}=0$).
Now that we have clarified how the availability of rental cars is regulated, we can state that a service can be chosen only if it is offered to the customer
\begin{equation}
  \label{eq:chooseIfAvailable}
  w_{ijkv}\leq y_{ikv}\qquad\forall i,j\in\mathcal{I},k\in \mathcal{K}_{ij},v\in\mathcal{V}\cup\mtc{A}
\end{equation}
Among the available services, the customer will chose the one yielding the highest utility.
Therefore, for a given zone $i\in\mtc{I}$, let decision variable $\nu_{ivwk}$ be equal to $1$ if both services $v$ and $w$ in $\mtc{V}\cup\mtc{A}$ are available to customer
$k\in\mtc{K}_{i}$, $0$ otherwise, and decision variable $\mu_{ijvwk}$ be equal to one if service $v\in\mtc{V}\cup\mtc{A}$
yields a greater utility than service $w\in\mtc{V}\cup\mtc{A}$ to customer $k\in\mtc{K}_{ij}$ moving from $i$ to $j$, $0$ otherwise.
The following constraints state that $\nu_{ivwk}$ is equal to one when both services $v$ and $w$ are available
\begin{align}
\label{eq:nu1}  &y_{ikv}+y_{ikw}\leq 1 + \nu_{ivwk}&\forall i\in\mathcal{I},k\in \mathcal{K}_{i},v,w\in\mathcal{V}\cup\mtc{A},\\
\label{eq:nu2}  &\nu_{ivwk}\leq y_{ikv}&\forall i\in\mathcal{I},k\in \mathcal{K}_{i},v,w\in\mathcal{V}\cup\mtc{A},\\
\label{eq:nu3}  &\nu_{ivwk}\leq y_{ikw}&\forall i\in\mathcal{I},k\in \mathcal{K}_{i},v,w\in\mathcal{V}\cup\mtc{A}.
\end{align}
A service is chosen only if it yields the highest utility
\begin{equation}\label{eq:chooseHighest}
  w_{ijkv}\leq \mu_{ijvwk}\qquad\forall i,j\in\mathcal{I},k\in \mathcal{K}_{ij},v,w\in\mathcal{V}\cup\mtc{A}
\end{equation}
that is, as soon as $\mu_{ijvwk}$ is set to $0$ for some index $w$, $w_{ijkv}$ is forced to take value $0$ and service $v$ is not chosen by customer $k$ on $i$-$j$.
The following constraints ensure that decision variable $\mu_{ijvwk}$ takes the correct value according to the utility
\begin{align}\label{eq:choice1}
  M_{ijk}\nu_{ivwk}-2M_{ijk}\leq u_{ijkv}-&u_{ijkw} - M_{ijk}\mu_{ijvwk}\\
  &\nonumber\forall i,j\in\mathcal{I},k\in \mathcal{K}_{ij},v,w\in\mathcal{V}\cup\mtc{A}
\end{align}
and
\begin{align}\label{eq:choice2}
  u_{ijkv}-u_{ijkw} - M_{ijk}\mu_{ijvwk}\leq &(1 -\nu_{ivwk})M_{ijk}\\
  &\nonumber\forall i,j\in\mathcal{I},k\in \mathcal{K}_{ij},v,w\in\mathcal{V}\cup\mtc{A}
\end{align}
where constant $M_{ijk}$ represents the greatest difference in utility between two services on $i-j$ for customer $k\in\mtc{K}_{ij}$,
that is $M_{ijk}\geq |u_{ijkv}-u_{ijkw}|,\forall v,w\in\mtc{V}\cup\mtc{A}$.
Constraints \eqref{eq:choice1}-\eqref{eq:choice2} work as follows.
When both services $v$ and $w$ are available ($\nu_{ivwk}=1$) and $u_{ijkv}>u_{ijkw}$, \eqref{eq:choice2} forces $\mu_{ijvwk}$ to take value $1$,
while \eqref{eq:choice1} reduces to $0\leq  u_{ijkv}-u_{ijkw}$. When both service $v$ and $w$ are available and $u_{ijkv}<u_{ijkw}$,
\eqref{eq:choice1} forces $\mu_{ijvwk}$ to take value $0$, while \eqref{eq:choice2} reduces to $0\geq  u_{ijkv}-u_{ijkw}$.
When one of the two services is not available ($\nu_{ivwk}=0$), constraints \eqref{eq:choice1}-\eqref{eq:choice2} are satisfied irrespective of the value of $\mu_{ijvwk}$.
In case of ties ($u_{ijkv}=u_{ijkw}$) we impose
\begin{equation}
\label{eq:onlyOneMu}  \mu_{ijvwk}+\mu_{ijwvk}\leq 1\qquad\forall i,j\in\mathcal{I},k\in \mathcal{K}_{ij},v,w\in\mathcal{V}\cup\mtc{A}
\end{equation}
A service can be preferred only if offered
\begin{equation}
  \label{eq:muIfAvail}\mu_{ijvwk}\leq y_{ikv} \qquad\forall i,j\in\mathcal{I},k\in \mathcal{K}_{ij},v,w\in\mathcal{V}\cup\mtc{A}
\end{equation}
Let decision variable $\alpha_{ijkvl}$ be equal to $1$ if fare $l$ is applied between $i$ and $j$ and customer $k$ chooses shared car $v\in\mtc{V}$, $0$ otherwise.
The following constraints ensure the relationship between $\lambda_{ijl}$ and $w_{ijkv}$ and $\alpha_{ijkvl}$
\begin{align}
  \label{eq:lin1}&\lambda_{ijl}+w_{ijkv}\leq 1 + \alpha_{ijkvl}&\forall v\in\mtc{V}, i,j\in\mtc{I},k\in\mtc{K}_{ij},l\in\mtc{L}\\
  \label{eq:lin2}&\alpha_{ijkvl}\leq\lambda_{ijl}&\forall v\in\mtc{V}, i,j\in\mtc{I},k\in\mtc{K}_{ij},l\in\mtc{L}  \\
  \label{eq:lin3}&\alpha_{ijkvl}\leq w_{ijkv}&\forall v\in\mtc{V}, i,j\in\mtc{I},k\in\mtc{K}_{ij},l\in\mtc{L}
\end{align}
That is, $\alpha_{ijkvl}$ is forced to take value $1$ as soon as both $\lambda_{ijl}$ and $w_{ijkv}$ take value one, and value $0$ as soon as either $\lambda_{ijl}$ or $w_{ijkv}$ take value $0$.

Finally, for a given realization $\xi:=(\xi_{kv})_{k\in\mtc{K},v\in\mtc{V}\cup\mtc{A}}$ of the random utility term $\tilde{\xi}:=(\tilde{\xi}_{kv})_{k\in\mtc{K},v\in\mtc{V}\cup\mtc{A}}$ the second-stage profit can be formally expressed as
\begin{subequations}
  \label{eq:2SF1}
\begin{align}
  \label{eq:obj:1}Q(z,\lambda,\xi) = \max~ & \sum_{v\in\mtc{V}}\sum_{(i,j)\in\mtc{I}\times\mtc{I}}\left(P^VT^{CS}_{ij}-C^U_{ij}\right)\sum_{k\in\mtc{K}_{ij}}w_{ijkv}\\
  \label{eq:obj:2}   & +\sum_{v\in\mtc{V}}\sum_{(i,j)\in\mtc{I}\times\mtc{I}}\sum_{k\in\mtc{K}_{ij}}\sum_{l\in\mtc{L}}L_{ijl}\alpha_{ijkvl}\\
  \text{s.t. }~&\eqref{eq:pV},\eqref{eq:pA},\eqref{eq:utility},\eqref{eq:onlyOne},\eqref{eq:availA},\eqref{eq:availV},\eqref{eq:availBefore},\eqref{eq:taken},\nonumber\\
                                           &\eqref{eq:chooseIfAvailable},\eqref{eq:nu1},\eqref{eq:nu2},\eqref{eq:nu3},\eqref{eq:chooseHighest},\eqref{eq:choice1},\eqref{eq:choice2},\eqref{eq:onlyOneMu},\eqref{eq:muIfAvail},\eqref{eq:lin1},\eqref{eq:lin2},\eqref{eq:lin3} \nonumber\\
\end{align}
\end{subequations}
where $C^U_{ij}$ is the cost born by the CSO when a vehicle is rented between $i$ and $j$, \eqref{eq:obj:1} represents the net revenue generated by the per-minute fee,
and \eqref{eq:obj:2} represents the income generated by the drop-off fee. 

Thus, we can formally express the expected profit (i.e., the recourse function) as
$$Q(z,\lambda):= \mathbb{E}_{\tilde{\xi}}\bigg[Q(z,\lambda,\xi)\bigg]$$
Problem \eqref{eq:1S} is a two-stage mixed-integer stochastic program with integer decision variables at both stages.

\section{L-Shaped Method}\label{sec:ls}
We propose a multi-cut Integer L-Shaped method to find exact solutions to problem \eqref{eq:1S}.
The original, single-cut, version of the method was introduced by \cite{LapL93}.
Assuming a set $\mtc{S}=\{1,\ldots,S\}$ of scenarios (e.g., an iid sample) of $\tilde{\xi}$, each with probability $\pi_s$,
the Master Problem (MP) can be formulated as 
\begin{subequations}
  \label{eq:F2MP}
  \begin{align}
    \label{eq:ref:obj}&\max-\sum_{v\in\mtc{V}}\sum_{i\in\mtc{I}}C^R_{vi}z_{vi}+\sum_{s\in\mtc{S}}\pi_s\phi_s\\
    \label{eq:ref:c3}&\sum_{i\in\mtc{I}}z_{vis} = 1   & v\in\mtc{V}\\
    \label{eq:ref:c6}&\sum_{l\in\mtc{L}}\lambda_{ijl}=1& i\in\mtc{I},j\in\mtc{J}\\
                      &z_{vi}\in\{0,1\}              & i\in\mtc{I},v\in\mtc{V}\\
                      &\lambda_{ijl}\in\{0,1\} & i\in\mtc{I},j\in\mtc{I},l\in\mtc{L}\\
                      &\phi_s \text{~free~}& s\in\mtc{S}.
  \end{align}
\end{subequations}
Let $\phi:=(\phi_s)_{s\in\mtc{S}}$. For each $s\in \mtc{S}$, the second-stage problem $Q(z,\lambda,\xi_s)$ is solved as a subproblem.
The L-Shaped method consists of solving MP in a Branch\& Cut framework where optimality cuts are added at (integer) nodes of the tree.
Observe that problem \eqref{eq:1S} has relatively complete recourse, that is, the second-stage problem $Q(z,\lambda,\xi_s)$ is feasible for every solution that satisfies the first-stage constraints.
Consequently, the method requires only the definition of optimality cuts.

The practical viability of the method is enabled by a compact reformulation of $Q(z,\lambda,\xi_s)$, provided in \Cref{sec:ls:f2}, and an efficient exact algorithm for its solution, introduced in \Cref{sec:ls:subproblems}.
We then provide the expression of optimality cuts and relaxation cuts in \Cref{sec:ls:oc} and \Cref{sec:ls:relaxation}, respectively. The former are necessary to cut off solutions for which $\phi_s<Q(z,\lambda,\xi_s)$, the latter provide non-trivial lower bounds to $Q(z,\lambda,\xi_s)$ and are
crucial for the efficiency of the algorithm. Final efficiency measures are provided in \Cref{sec:ls:efficiency}.

\subsection{Compact formulation of the second-stage problem}\label{sec:ls:f2}
For a given realization $\xi$ of $\tilde{\xi}$, and first-stage decision $(z,\lambda)$, the second-stage problem can be reformulated by preprocessing customer choices.
The preprocessing phase ensures that the reformulation is linear regardless of whether the utility function adopted is linear in the price.
Thus, any choice model can be used, without any restriction to utility function linear in the price.

We introduce the concept of a request. A request represents a customer who wishes to use carsharing for moving from its origin to its destination.
For a given realization $\xi$, let the set $\mtc{R}(\xi)$ be the set of requests. The set $\mtc{R}(\xi)$ contains a request for each customer $k\in\mtc{K}$ for which there exists at least one drop-off level $l\in\mtc{L}$
such that the customer would prefer carsharing to alternative transport services, that is, for which $u_{ijkv}>u_{ijkw}$ with $v\in\mtc{V}$ and $w\in\mtc{A}$ for some choice of $l\in\mtc{L}$ (note that all shared cars yield the same utility).
Let $i(r)$, $j(r)$ and $k(r)$ be the origin, destination and customer of request $r$, respectively, and $l(r)$ the highest drop-off fee at which customer $k(r)$ would prefer carsharing to other services.
Note that customer $k(r)$ would still prefer carsharing at any drop-off fee lower than $l(r)$ (under the reasonable assumption that the customer is sensitive to price).
For each realization $\xi$ of $\tilde{\xi}$ the set of requests can be populated in $\mtc{O}\big(|\mtc{K}|\times|\mtc{L}|\times |\mtc{A}|\big)$ operations as described in \Cref{alg:Rs}.

\begin{algorithm}[h]
  \caption{Computation of $\mathcal{R}(\xi)$.}
  \label{alg:Rs}
  \begin{algorithmic}[1]
    \STATE Input: $\xi$
    \STATE $\mtc{R}(\xi)\gets \emptyset$
    \FOR{customer $k\in\mtc{K}$}
    \STATE $i\gets i(k)$, $j\gets j(k)$ \COMMENT{{\footnotesize $i(k)$ and $j(k)$ are the origin and destination of customer $k$.}}
    \STATE $l^{MAX}\gets -\infty$ \COMMENT{{\footnotesize The highest drop-off fee at which customer $k$ will choose carsharing.}}
    \STATE $L_{l^{MAX}}\gets -\infty$
    \FOR{drop-off level $l\in\mtc{L}$}
    \STATE $p^{CS} \gets P^VT^{CS}_{ij} + L_{l}$ \COMMENT{{\footnotesize Calculate the price of a carsharing ride.}}
    \STATE $U^{CS} \gets F_k(P^{CS},\pi^1_{vij},\ldots,\pi^N_{vij})+\xi_{kv}$ for some $v\in\mtc{V}$\COMMENT{{\footnotesize Calculate the utility of carsharing.}}
    \FOR{service $v\in\mtc{A}:Y_{vi}=1$}
    \STATE $p^A_v \gets P_{vij}$ \COMMENT{{\footnotesize Calculate the price of a ride with alternative $v$.}}
    \STATE $U^A_v \gets F_k(P^{A}_v,\pi^1_{vij},\ldots,\pi^N_{vij})+\xi_{kv}$ \COMMENT{{\footnotesize Calculate the utility of alternative $v$.}}
    \ENDFOR
    \IF{$U^{CS} > \max_{v\in\mtc{A}}\{U^A_v\}$ and $L_l > L_{l^{MAX}}$ }
    \STATE $l^{MAX}\gets l$
    \STATE $L_{l^{MAX}}\gets L_l$
    \ENDIF
    \ENDFOR
    \IF{$l^{MAX}> -\infty$} 
    \STATE $r\gets |\mtc{R}(\xi)|+1$ \COMMENT{{\footnotesize In this case there exists a drop-off fee at which $k$ prefers carsharing. Thus we create a request.}}
    \STATE $\mtc{R}(\xi)\gets \mtc{R}(\xi)\cup \{r\}$
    \STATE $i(r)\gets i$
    \STATE $j(r)\gets j$
    \STATE $k(r)\gets k$
    \STATE $l(r)\gets l^{MAX}$
    \ENDIF
    \ENDFOR
    \RETURN $\mtc{R}(\xi)$
  \end{algorithmic}
\end{algorithm}
Observe that, in case of a heterogeneous fleet, \Cref{alg:Rs} should be edited at lines $8-9$ and $14$ to account for the fact that each vehicle gives, in general, a different utility.

Let then $R_{rl} = P^VT^{CS}_{i(r),j(r)}-C^U_{i(r),j(r)}+L_l$, for $l\leq l(r)$, be the net revenue generated if request $r$ is satisfied at drop-off fee level $l$.
Let $\mtc{R}_r(\xi)=\{\rho \in\mtc{R}(\xi): i(\rho) = i(r), k(\rho) < k(r)\}$ be the set of requests which have a precedence over $r$.
Let $\mtc{R}_{ij}(\xi)=\{r\in\mtc{R}(\xi):i(r) = i, j(r)=j\}$.
Let $\mtc{L}_r(\xi)=\{l\in\mtc{L}:L_l \leq L_{l(r)}\}$.
Finally, let decision variable $y_{vrl}$ be equal to $1$ if request $r$ is satisfied by vehicle $v$ at level $l$, $0$ otherwise.
We can now reformulate the second-stage problem as follows.
\begin{subequations}
  \label{eq:2SF2}
  \begin{align}
    \label{eq:2SF2:obj}&Q(z,\lambda,\xi)=\max\sum_{r\in\mtc{R}(\xi)}\sum_{v\in\mtc{V}}\sum_{l\in\mtc{L}_r(\xi)}R_{vrl}y_{vrl}\\
    \label{eq:2SF2:c1}&\sum_{v\in\mtc{V}}\sum_{l\in\mtc{L}_r(\xi)}y_{vrl}\leq 1   & r\in\mtc{R}(\xi)\\
    \label{eq:2SF2:c2}&\sum_{r\in\mtc{R}(\xi)}\sum_{l\in\mtc{L}_r(\xi)}y_{vrl}\leq 1   & v\in\mtc{V}\\
    \label{eq:2SF2:c4}&\sum_{l\in\mtc{L}_{r_1}(\xi)}y_{v,r_1,l}  + \sum_{r_2\in\mtc{R}_{r_1}(\xi)}\sum_{l\in\mtc{L}_{r_2}(\xi)}y_{v,r_2,l}\leq  z_{v,i(r_1)}   & r_1\in\mtc{R}(\xi),v\in\mtc{V}\\[5pt]
    \nonumber&y_{v,r_1,l_1}+ \sum_{r_2\in\mtc{R}_{r_1}(\xi)}\sum_{l_2\in\mtc{L}_{r_2}(\xi)}y_{v,r_2,l_2}+ \sum_{v_1\in\mtc{V}:v_1\neq v}y_{v_1,r_1,l_1}  &\\
    \label{eq:2SF2:c5} &\geq \lambda_{i(r_1),j(r_j),l_1} + z_{v,i(r_1)}-1   & r_1\in\mtc{R}(\xi),v\in\mtc{V}, l_1\in\mtc{L}_{r_1}(\xi)\\
    \label{eq:2SF2:c7}&\sum_{v\in\mtc{V}}y_{vrl}\leq\lambda_{i(r),j(r),l}& r\in\mtc{R}(\xi),l\in\mtc{L}_r(\xi)\\
                      &y_{vrl}\in\{0,1\}              & r\in\mtc{R}(\xi),v\in\mtc{V},l\in\mtc{L}_r(\xi)
  \end{align}
\end{subequations}

The objective function \eqref{eq:2SF2:obj} represents the net revenue obtained by the satisfaction customer requests.
Constraints \eqref{eq:2SF2:c1} ensure that each request is satisfied at most once.
Constraints \eqref{eq:2SF2:c2} ensure that each vehicle satisfies at most one request.
Constraints \eqref{eq:2SF2:c4} state that a request can be satisfied by vehicle $v$ only if the vehicle is in zone $i(r_1)$
and the vehicle has not been assigned to a customer with a lower index (that is, arriving at the vehicle before $k(r_1)$).
Constraints \eqref{eq:2SF2:c5} state that a request $r_1$ at a certain level $l_1$ must be satisfied by vehicle $v$
if level $l_1$ has been chosen ($\lambda_{i(r_1),j(r_1),l_1}=1)$ and the vehicle is available at $i(r_1)$ ($z_{v,i(r_1)}=1)$, unless the car has been used to satisfy the request of a customer with a higher priority (second term on the left-hand-side), or $r_1$ has been satisfied by another vehicle (third term on the left-hand-side).
Constraints \eqref{eq:2SF2:c7} state that a request can be satisfied at level $l$ only if level $l$ is applied to all customers traveling between $i$ and $j$.
Note that $z$ and $\lambda$ are input data in problem \eqref{eq:2SF2}.

\subsection{Solution of the second-stage problem}\label{sec:ls:subproblems}

Given a solution $(z,\lambda)$ to MP and a scenario $\xi_s$, the optimal second stage profit $Q(z,\lambda,\xi_s)$ and solution can be computed by the greedy procedure sketched in \Cref{alg:Q}.

\begin{algorithm}[h]
  \caption{Greedy algorithm for computing $Q(z,\lambda,\xi_s)$ and its optimal solution.}
  \label{alg:Q}
  \begin{algorithmic}[1]
    \STATE INPUT: $z$, $\lambda$, $\mtc{R}(\xi_s)$.
    \STATE $\mtc{V}^A_s\gets \mtc{V}$\COMMENT{{\footnotesize $\mtc{V}^A_s$ is the set of available vehicles.}}
    \STATE $Y_{vrl} \gets 0$ $\forall v\in\mtc{V}, r\in\mtc{R}(\xi_s),l\in\mtc{L}_r(\xi_s)$.
    \STATE $Q(z,\lambda,\xi_s)\gets 0$
    \STATE Sort requests $\mtc{R}(\xi_s)$ in non-decreasing order of the customer index $k(r)$\COMMENT{{\footnotesize Remember that customers with a lower index have the precedence over customers with a higher index, see \eqref{eq:taken}}.}
    \FOR{request $r\in \mtc{R}(\xi_s)$}
    \STATE $L_{i(r),j(r)} = \sum_{l\in\mathcal{L}}l\lambda_{i(r),j(r),l}$\COMMENT{{\footnotesize Identify the fee applied between $i(r)$ and $j(r)$.}}
    \IF{$L_{i(r),j(r)}\leq l(r)$}
    \FOR{$v\in\mtc{V}^A_s$}
    \IF{$z_{v,{i(r)}}=1$}
    \STATE $Y_{v,r,L_{i(r),j(r)}}\gets 1$
    \STATE $\mtc{V}^A_s\gets \mtc{V}^A_s\setminus \{v\}$\COMMENT{{\footnotesize Vehicle $v$ becomes unavailable}}
    \STATE $Q(z,\lambda,\xi_s)\gets Q(z,\lambda,\xi_s) + R_{r,L_{i(r),j(r)}}$
    \ENDIF
    \ENDFOR
    \ENDIF
    \ENDFOR
    \RETURN $Q(z,\lambda,\xi_s)$ and $Y_{vrl} \forall v\in\mtc{V}, r\in\mtc{R}(\xi_s),l\in\mtc{L}(\xi_s)$.
  \end{algorithmic}
\end{algorithm}

\Cref{alg:Q} proceeds as follows. Given a solution $(z,\lambda)$ to MP and the requests available in scenario $\xi_s$, the algorithm first initializes the solution $Y_{vrl}$, the objective value $Q(z,\lambda,\xi_s)$ and the set of available vehicles.
Then it sorts the requests in non-decreasing order of the customer index $k(r)$. This is necessary to enforce that customers with a lower index have their request satisfied before customers with a higher index.
The algorithm then iterates over the ordered requests. For each request it first checks whether the fee applied on its origin-destination pair, $L_{i(r),j(r)}$, is lower than then highest drop-off fee acceptable to the customer, $l(r)$.
If this is the case, the algorithm looks for vehicles available at the origin of request $r$, $i(r)$. If one such vehicle $v$ is found, the request is assigned to the vehicle at the current drop-off fee, i.e, $Y_{v,r,L_{i(r),j(r)}}$ is set to $1$,
vehicle $v$ is made unavailable, and the revenue is increased by the revenue of request $r$, that is $R_{r,L_{i(r),j(r)}}$.
The algorithm performs $\mathcal{O}\big(|\mtc{R}(\xi_s)|\times |\mtc{V}|\big)$ operations.

\subsection{Optimality cuts}\label{sec:ls:oc}
We are now concerned with finding a \textit{valid set of optimality cuts}, that is a finite number of optimality cuts which enforce $\phi_s\leq Q(z,\lambda,\xi_s)$ for all $s\in\mtc{S}$.

Assume an upper bound $U_s$ on $\max_{z,\lambda}Q(z,\lambda,\xi_s)$ exists for all $s\in\mtc{S}$.
In the case of a homogeneous fleet a valid upper bound $U_s$ on $\max_{z,\lambda}Q(z,\lambda,\xi_s)$ is
$$U_s=\sum_{r\in\mtc{R}(\xi_s)}\max\{R_{r,l(r)} , 0 \}$$
That is, the upper bound assumes that all, and only, the requests which generate a positive revenue are satisfied, and that these are satisfied at the highest drop-off fee $l(r)$. 
If the fleet is not homogeneous, we would have a vehicle-specific revenue $R_{vrl}$ for satisfying request $r$ at level $l$. In this case,
a valid upper bound can be obtained by assuming all requests are satisfied by the vehicle which yields the highest non-negative revenue, i.e.,
$$U_s=\sum_{r\in\mtc{R}(\xi_s)}\max\bigg\{ \max_{v\in\mtc{V}}\{R_{v,r,l(r)}\} , 0 \bigg\}$$

Given a solution $(z^t,\lambda^t)$ to MP (e.g., at a given node $t$ of the Branch \& Cut procedure),
let $\mtc{Z}^+_t\subseteq \mtc{V}\times\mtc{I}$ and $\mtc{Z}^-_t\subseteq \mtc{V}\times\mtc{I}$ be the set of tuples $(v,i)$ for which $z_{vi}^t=1$ and $z_{vi}^t=0$, respectively.
Similarly, let $\Lambda^+_t\subseteq\mtc{I}\times\mtc{I}\times\mtc{L}$ and $\Lambda^-_t\subseteq\mtc{I}\times\mtc{I}\times\mtc{L}$ be the set of tuples $(i,j,l)$ for which $\lambda_{ijl}=1$ and $\lambda_{ijl}=0$, respectively.
\Cref{prop:ocs} defines a valid set of optimality cuts.
\begin{proposition}
  \label{prop:ocs}
  Let $(z^t,\lambda^t)$ be the $t$-th feasible solution to MP, and $Q(z,\lambda,\xi_s)$ its second-stage value for scenario $s$. The set of cuts
  \begin{align}
    \label{eq:oc}
    \phi_s \leq &\bigg(Q(z,\lambda,\xi_s) - U_s\bigg)\bigg(\sum_{(v,i)\in\mtc{Z}^+_t}z_{vi}-\sum_{(v,i)\in\mtc{Z}^-_t}z_{vi}+\sum_{(i,j,l)\in\Lambda^+_t}\lambda_{ijl}-\sum_{(i,j,l)\in\Lambda^-_t}\lambda_{ijl}\bigg)\\
    \nonumber &+  U_s - \bigg(Q(z,\lambda,\xi_s) - U_s\bigg)\bigg(|\mtc{Z}^+_t|+|\Lambda^+_t|-1 \bigg)
  \end{align}
  defined for all $(z^t,\lambda^t)$ feasible to MP is a valid set of optimality cuts.
  \begin{proof}
    It is sufficient to observe that, for $(z,\lambda)=(z^t,\lambda^t)$, we have
    $$\bigg(\sum_{(v,i)\in\mtc{Z}^+_t}z_{vi}-\sum_{(v,i)\in\mtc{Z}^-_t}z_{vi}+\sum_{(i,j,l)\in\Lambda^+_t}\lambda_{ijl}-\sum_{(i,j,l)\in\Lambda^-_t}\lambda_{ijl}\bigg)=|\mtc{Z}^+_t|+|\Lambda^+_t|$$
    and optimality cut \eqref{eq:oc} reduces to
    $$\phi_s\leq Q(z,\lambda,\xi_s), \qquad \forall s\in\mtc{S}$$
    
    On the other hand, if $(z,\lambda)\neq(z^t,\lambda^t)$ we get
    $$\bigg(\sum_{(v,i)\in\mtc{Z}^+_t}z_{vi}-\sum_{(v,i)\in\mtc{Z}^-_t}z_{vi}+\sum_{(i,j,l)\in\Lambda^+_t}\lambda_{ijl}-\sum_{(i,j,l)\in\Lambda^-_t}\lambda_{ijl}\bigg)\leq|\mtc{Z}^+_t|+|\Lambda^+_t|-1$$
    and, observing that $Q(z,\lambda,\xi_s) - U_s\leq 0$, the right-hand-side of the cut becomes greater than or equal to $U_s$.

    Thus, since MP is a maximization problem, the set of cuts enforces $\phi_s\leq Q(z,\lambda,\xi_s)$ when $(z,\lambda)=(z^t,\lambda^t)$ and yields a valid upper bound on the remaining solutions.
  \end{proof}
\end{proposition}

\subsection{Relaxation cuts}\label{sec:ls:relaxation}
Ordinary Benders decomposition cuts can be derived by solving the LP relaxation of the subproblems $Q(z,\lambda,\xi_s)$. Relaxation cuts are not to be confused with optimality cuts as they are, in general, not tight at the point $(z^t,\lambda^t)$ at which they are generated.
However, relaxation cuts provide a, possibly, non-trivial upper bound on $Q(z,\lambda,\xi_s)$, that is an upper bound which might be lower than $U_s$ for several $(z,\lambda)$ solutions.
A relaxation cut is obtained as follows
\begin{subequations}
  \label{eq:rc}
  \begin{align}
    \phi_s \leq& \sum_{r\in\mtc{R}(\xi_s)}\pi^A_r+\sum_{v\in\mtc{V}}\pi^B_{v}+\sum_{r\in\mtc{R}(\xi_s)} \sum_{v\in\mtc{V}} \pi^C_{rv}z_{v,i(r)}\\
               &+\sum_{r\in\mtc{R}(\xi_s)}\sum_{v\in\mtc{V}}\sum_{l\in\mtc{L}_{r}(\xi_s)}\pi^D_{rvl}\bigg(\lambda_{i(r),j(r),l} + z_{v,i(r)}-1\bigg)+ \sum_{r\in\mtc{R}(\xi_s)}\sum_{l\in\mtc{L}_r(\xi_s)}\pi^E_{rl}\lambda_{i(r),j(r),l}
  \end{align}
\end{subequations}
where $\pi^A_r$, $\pi^B_v$, $\pi^C_{rv}$, $\pi^D_{rvl}$, $\pi^E_{rl}$ are the values of the dual solution to $Q(z,\lambda,\xi_s)$ corresponding to constraints \eqref{eq:2SF2:c1}, \eqref{eq:2SF2:c2}, \eqref{eq:2SF2:c4}, \eqref{eq:2SF2:c5} and \eqref{eq:2SF2:c7}, respectively.

\subsection{Other efficiency measures}\label{sec:ls:efficiency}
MP exhibits symmetric solutions. If no vehicle is made available in a given zone $i$, all configurations of the drop-off fees between zone $i$ and the remaining zones $j$ are equivalent.
In fact, no customer will be served on those $(i,j)$ pairs. This problem can be solved by mean of the following constraints
\begin{align}\label{eq:VI}
  &\sum_{v\in\mtc{V}}z_{vi}+\lambda_{ij1}\geq 1& \forall i,j\in\mtc{I}
\end{align}
Constraints \eqref{eq:VI} enforce that, when no vehicle is available in zone $i$ (i.e., $\sum_{v\in\mtc{V}}z_{vi}=0$) we arbitrarily chose drop-off fee number $1$ (i.e., $\lambda_{ij1}=1$). On the other hand, if $\sum_{v\in\mtc{V}}z_{vi}\geq 1$ constraints \eqref{eq:VI} are satisfied regardless of the choice of a drop-off fee.

\section{Test instances}\label{sec:instances}
In this section we present the test instances we used to run a computational study whose results are presented in \Cref{sec:results}.
The test instances mimic carsharing services in the Italian city of Milan. The city hosts a number of carsharing companies and, according to the municipality of Milan \cite{Mil20},
in $2018$ there were a total of $3 108$ free-floating shared vehicles, with an average of $16 851$ daily rentals, and $149$ station-based shared vehicles with an average of $108$ daily rentals.
We start by describing how the instances were constructed and finally we clarify the specific choices of control parameters for our tests.
For the sake of replicability, an instance generator is made publicly available at the address \url{https://github.com/GioPan/instancesPricingAndRepositioningProblem}.

\subsection{Zones and alternative transport services}
\begin{figure}
  \centering
  \includegraphics[width=0.6\textwidth]{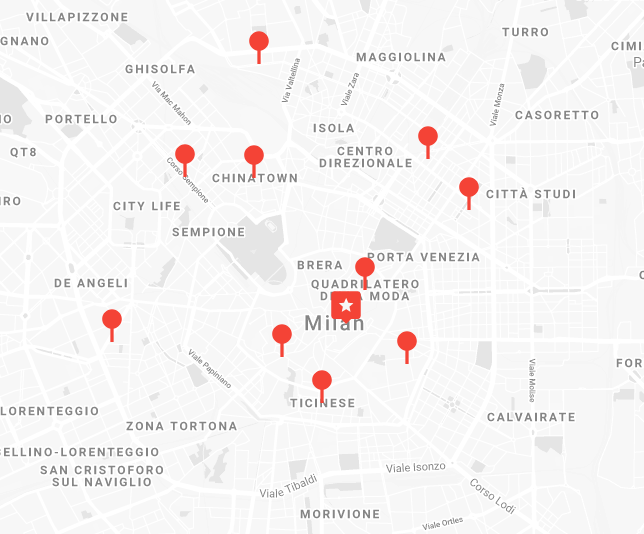}
  \caption{Municipality of Milan. Pins identify the ten locations of the city of Milan in the instances of \cite{HanP18}. The star identifies the Dome of Milan used as a reference point of the center of the city. }
  \label{fig:mi}
\end{figure}

We build upon, and expand, the instances used by \cite{HanP18}.
The authors consider ten key locations in the business area of the city of Milan which we use as representatives of as many zones, thus setting $\mtc{I}=\{1,\ldots,10\}$, see \Cref{fig:mi}.
The authors consider as alternative transport services \textit{public transport} (PT -- consisting of a combination of busses, metro and superficial trains) and \textit{bicycles} (B).
Therefore we set $\mtc{A}=\{PT,B\}$.
For each pair of zones, the authors provide all the information necessary to calculate the utility as further explained in \Cref{sec:instances:utility}.

\subsection{Customers and utility functions}\label{sec:instances:utility}
According to \eqref{eq:utility} each customer $k$ is characterized by a known utility function and a random variable $\tilde{\xi}_{kv}$ which represents the portion of the preferences of the customer with respect to service $v$ that the CSO cannot explain. We start by introducing the portion of the utility estimated by the CSO.

 We adopt the utility function described in \cite{ModS98}.
This function provides an estimate of the utility of each transport service as a function of price, travel time, walking time (e.g., to reach the service and from the service to destination) and waiting time. These represent the characteristics observable by the CSO.

For each customer $k\in\mtc{K}$ traveling between $i$ and $j$ with transportation service $v$ the utility function is
\begin{align}\label{eq:cs:utility}
  \nonumber F_k(p_{vij},T^{CS}_{vij},T^{PT}_{vij},&T^B_{vij},T^{Walk}_{vkij}, T^{Wait}_{vij})= \beta_k^{P} p_{vij} + \beta_k^{CS} T_{vij}^{CS}+ \beta_k^{PT} T_{vij}^{PT}\\
                                                  &+\tau(T_{vij}^{B})\beta_k^{B} T_{vij}^{B}+\tau(T_{vij}^{Walk})\beta_k^{Walk} T_{vij}^{Walk}+
                                                    \beta_k^{Wait} T_{vij}^{Wait}  
\end{align}
The meaning of each parameter and coefficient of function \eqref{eq:cs:utility} is clarified in \Cref{tab:instances:utilityparams} and function  $\tau:\mathbb{R}\rightarrow\mathbb{R}$, defined as $\tau(t)=\lceil \frac{t}{10}\rceil$, allows us to model the utility of cycling and walking as a piece-wise linear function: the utility of walking and cycling decreases faster as the walking and cycling time increases, see \cite{ModS98}.
\begin{table}
  \centering
  \caption{Parameters and coefficients of the utility function.}\label{tab:instances:utilityparams}
  \begin{tabular}{p{0.1\linewidth}|p{0.8\linewidth}}
    \toprule
    Parameter & Meaning\\
    \midrule
    $T^{CS}_{vij}$ & Time spent riding a shared car between $i$ and $j$ when using service $v$. This quantity is strictly positive only when $v$ is a carsharing service, otherwise it is $0$.\\
    $T^{PT}_{vij}$ & Time spent in public transportation between $i$ and $j$ when using service $v$. This quantity is strictly positive only when $v$ is PT, otherwise it is $0$.\\
    $T^{B}_{vij}$ & Time spent riding a bicycle between $i$ and $j$ when using service $v$. This quantity is strictly positive only when $v$ is B, otherwise it is $0$.\\
    $T^{Walk}_{vij}$& Walking time necessary when moving with transportation service $v$ between $i$ and $j$. This includes the walking time to the nearest service (e.g., shared car or bus stop),
  between connecting means (e.g., when switching between bus and metro to reach the final destination), and from to the final destination.\\
    $T^{Wait}_{vij}$& Waiting time when using service $v$ between $i$ and $j$, and includes the waiting time for the service (e.g., bus or metro) as well as for connection.\\
    $beta_k^{P}$ & Price sensitivity of customer $k$.\\
    $\beta_k^{CS}$ & Time sensitivity of customer $k$ when using a shared vehicle.\\
    $\beta_k^{PT}$ &Time sensitivity of customer $k$ when using public transport.\\
    $\beta_k^{B}$ & Time sensitivity of customer $k$ when riding a bicycle.\\
    $\beta_k^{Walk}$ & Time sensitivity of customer $k$ when walking.\\
    $\beta_k^{Wait}$ & Time sensitivity of customer $k$ when waiting.\\
    \bottomrule
  \end{tabular}
\end{table}

We use the $\beta$ coefficients of the original utility function provided by \cite{ModS98} and marginally adapted to the carsharing context by \cite{HanP18} (e.g., price sensitivity has been adapted from Italian Lira to Euro). The values of the coefficients are the following: $\beta^{CS} = -1$, $\beta^{PT} =-2$, $\beta^{B}=-2.5$, $\beta^{Walk}=-3$ and $\beta^{Wait}=-6$.
For the price sensitivity $\beta^P$ the authors create two customer segments. They assign $\beta^P = - 188.33$ if a customer belongs to the \textit{lower-middle class} or $\beta^P = - 70.63$ if a customer belongs to the \textit{upper-middle class}. We randomly assigned customers (with equal probability) to the either upper-middle class ($\beta^P_k=-70.63$) or lower-middle class ($\beta^P_k=-188.33$).
In more general cases, the parameters of utility functions can be estimated, provided the availability of data records on actual customers choices. The estimation procedure itself depends on several elements and underlying assumptions. As an example, a classical procedure to estimate the parameters of a Logit model is to maximize the log-likelihood function. Alternative methods include maximum simulated likelihood, simulated moments as well as Bayesian estimation. This topic is treated in detail in, e.g., \cite{Tra09}.

For the time parameters ($T$-parameters) in \Cref{tab:instances:utilityparams} we use the values estimated by \cite{HanP18} on the actual transport services in the city of Milan in $2017$.
These values can also be found in the files accompanying the instance generator we make available online at \url{https://github.com/GioPan/instancesPricingAndRepositioningProblem}.
It should be noted that, in more general applications, the $T$-parameters for the different services might change significantly during the day as a result of issues such as different traffic patterns, road congestion or time-varying public transport schedules.
Therefore, the $T$-parameters should be understood as specific for the target period under consideration.

The price parameters are set as follows. The price of a bicycle ride is set to $P_{Bij}=0$ for all $(i,j)$ pairs, the price for public transport services is $P_{PT,ij}=2$ (in Euro) for all $(i,j)$ corresponding to the current price of an ordinary ticket valid for $90$ minutes between each origin and destination within the municipality of Milan (price valid on November 2020).
The per-minute fee of carsharing is set to $P^V_v=0.265$ Euros per minute (the average of current per-minute fees offered by the CSOs in the city of Milan).
The drop-off fees considered are $L_1=-2$, $L_2=-1$, $L_3=0$, $L_4=1$, $L_5=2$ Euro in the base case ($\mtc{L}=\{1,2,3,4,5\}$). Further analysis on the drop-off fee will be described in \Cref{sec:results}.

\subsection{Individual customers profiles}\label{sec:instances:individual}
The utility function in \Cref{sec:instances:utility} entails that all customers within a given class (upper- or lower-middle class) are characterized by identical preferences with respect to travel time, price, waiting time and walking time. Customers are told apart by their preferences with respect to unobserved features of the services, captured by $\tilde{\xi}_{kv}$.

However, the availability of large amounts of customer data may allow the CSO to profile customers at the individual level, i.e., to assign each customer an individual utility function.
We are not aware of publicly available utility functions which are able distinguish between individual customers.
Therefore, in order to assess the effect of individual customer profiles, at least on the the performance of the algorithm,
we use an additional configuration in which an individual utility function for each customer is obtained by applying a random perturbation to the coefficients provided by \cite{HanP18}.
Particularly, for each customer $k$, $\beta^P_k$ will be uniformly drawn in $[-188.33,-70.63]$, where $-188.33$ is the $\beta^P_k$ coefficient for lower-middle class customers and $-70.63$ is the $\beta^P_k$ coefficient of upper-middle class customers in the general case, see \Cref{sec:instances:utility}.
This allows us to obtain customers which can be anywhere between the upper- and lower-middle class. 

The remaining $\beta$ coefficients will be uniformly drawn in $[0.8\beta,1.2\beta]$, where $\beta$ is the value provided by \cite{HanP18}.
As an example, for each $k$ we will draw $\beta^{PT}_k$ in $[-1.6,-2.4]$. The lower $\beta^{PT}_k$ the less utility the customer will obtain for each minute spent in public transportation.

\subsection{Uncertainty}\label{sec:uncertainty}
The random term of the utility $\tilde{\xi}_{kv}$ is modeled as a Gumbel (Extreme Value type I) distribution with mean $0$ and standard deviation $\sigma$.
This corresponds to using a Logit choice model (see \cite{Tra09,BenB99}).
The value of $\sigma$ is set as the empirical standard deviation of $U_{ijkv}= F_k(p_{vij},T^{CS}_{vij},T^{PT}_{vij},T^B_{vij},T^{Walk}_{vij}, T^{Wait}_{vij})$ for all $i,j\in\mtc{I},v\in\mtc{V}\cup\mtc{A},k\in\mtc{K}_{ij}$.
This entails that the expectation term in the objective function of \eqref{eq:1S} (i.e., $Q(z,\lambda)$) is a multidimensional integral that makes the solution of the problem prohibitive.
For this reason we approximate $\tilde{\xi}_{kv}$ by iid samples drawn from its the underlying Gumbel distributions. The resulting discrete stochastic program goes under the name of \textit{Sample Average Approximation} (SAA), see \cite{KleSH02}. Its optimal objective value provides an unbiased estimator of the true objective value. The full model of the SAA is provided in \Cref{sec:saa}.

\subsection{Position of customers and vehicles}
We partition customers into sets $\mtc{K}_i$ and then further into sets $\mtc{K}_{ij}$ is such a way to test different configurations of demand, e.g., center to outskirt and vice-versa.
Each one of the ten zones in our instances is characterized by a degree of ``centrality''. We use the walking distance from the \textit{Dome of Milan} as a proxy of centrality, see \Cref{fig:mi}.
Let $d_i$ be the walking distance from zone $i\in\mtc{I}$ to the Dome. Customers $\mtc{K}$ are first randomly partitioned
into disjoint subsets $\mtc{K}_{i}$ with a probability $\pi_i$ which depends on the centrality of the zone as follows
\begin{equation}
  \label{eq:probability}\pi_i = \frac{\gamma_id_i}{\sum_{i \in\mtc{I}}\gamma_id_i}
\end{equation}
where $\gamma_i=e^{-\alpha^{FROM}\Delta_i}$ with $\alpha^{FROM}\in[0,1]$ and $\Delta_i=d_i-\sum_{i \in\mtc{I}}d_i/|\mtc{I}|$ is the deviation from the mean distance.
In words, as $\alpha^{FROM}$ increases, the zones closer to the center (negative $\Delta_i$) will receive a higher probability and the zones far from the center a lower probability, resulting in a higher concentration of customers in the central zones.
Further, all customers assigned to a given zone $i$ will be randomly assigned a destination zone $j$, and thus inserted into subset $\mtc{K}_{ij}$, with a probability \eqref{eq:probability}.
This time $\gamma_j=e^{-\alpha^{TO}\Delta_i}$ with $\alpha^{TO}\in[0,1]$. Again, as $\alpha^{TO}$ increases more customers will be directed to central zones.
As an example, setting a low value of $\alpha^{FROM}$ and a high value of $\alpha^{TO}$ will create instances with higher demand from the outskirt to the center.
The partitioning of customers is sketched in \Cref{alg:partitionK}.

\begin{algorithm}[h]
  \caption{Algorithm for the partition of customers into subsets $\mtc{K}_i$ and $\mtc{K}_{ij}$, $(i,j)\in\mtc{I}\times\mtc{I}$.}
  \label{alg:partitionK}
  \begin{algorithmic}[1]
    \STATE Input: $\mtc{K}$, $\mtc{I}$, $d_i$ for $i\in \mtc{I}$, $\alpha^{FROM}\in[0,1]$, $\alpha^{TO}\in[0,1]$
    \STATE $\mtc{K}_i=\mtc{K}_{ij}\gets\emptyset$ for $(i,j)\in\mtc{I}\times\mtc{I}$
    \FOR{ zone $i\in \mtc{I}$}
    \STATE Calculate $\Delta_i=d_i-\sum_{i \in\mtc{I}}d_i/|\mtc{I}|$
    \STATE Calculate $\gamma_i^{FROM}=e^{-\alpha^{FROM}\Delta_i}$ and $\gamma_i^{TO}=e^{-\alpha^{TO}\Delta_i}$
    \ENDFOR
    \FOR{ zone $i\in \mtc{I}$}
    \STATE Calculate $\pi_i^{FROM} = \frac{\gamma_i^{FROM}d_i}{\sum_{i \in\mtc{I}}\gamma_i^{FROM}d_i}$
    \STATE Calculate $\pi_i^{TO} = \frac{\gamma_i^{TO}d_i}{\sum_{i \in\mtc{I}}\gamma_i^{TO}d_i}$
    \ENDFOR      
    \FOR{Customer $k\in \mtc{K}$}
    \STATE Draw an origin zone $i$ from $\mtc{I}$ according to the probability distribution $(\pi_i)^{FROM}_{i\in\mtc{I}}$
    \STATE $\mtc{K}_{i}\gets\mtc{K}_{i}\cup\{k\}$
    \STATE Draw a destination zone $j$ from $\mtc{I}$ according to the probability distribution $(\pi_i)^{TO}_{i\in\mtc{I}}$
    \STATE $\mtc{K}_{ij}\gets\mtc{K}_{ij}\cup\{k\}$
    \ENDFOR
    \RETURN $\mtc{K}_i$ and $\mtc{K}_{ij}$ for $(i,j)\in\mtc{I}\times\mtc{I}$
  \end{algorithmic}
\end{algorithm}

We assume the decision maker is a CSO with a fleet of $|\mtc{V}|$ homogeneous vehicles. Each vehicle $v$ is randomly assigned to an initial zone $i$ according to probability \eqref{eq:probability},
where $\gamma_i=e^{-\alpha^V\Delta_i}$ with $\alpha^V\in[0,1]$. Also in this case, as $\alpha^V$ increases more cars will be initially located in central zones.

\subsection{Costs}
We assume a fleet of Fiat 500 cars with classical combustion engine. The relocation cost $C^R_{ij}$, equal for all vehicles, is set as the cost of the fuel necessary for a ride between $i$ and $j$, plus the per-minute salary of the driver multiplied by the driving time.
The per-minute salary of the driver is set to $0.20$ Euro/minute. It is calculated from the Italian national collective contract for logistics services valid at October 1st 2019 (available at \url{https://www.lavoro-economia.it/ccnl/ccnl.aspx?c=328}) as follows:
the average per minute salary of the five lowest salary levels is increased by $30\%$ to account for e.g., night shifts and holidays, yielding approximately $0.20$ Euro/minute.
Finally, the cost $C^U_{ij}$ is set equal to the fuel necessary for a ride between $i$ and $j$. The fuel consumption is calculated based on the specifics of a Fiat 500 petrol car and assuming an average speed of $50$km/h and a fuel price of $1.60$ Euro/liter.
All the data data necessary to generate the instances described, as well as an instance generator implemented in Java, are made available at \url{https://github.com/GioPan/instancesPricingAndRepositioningProblem}.

\subsection{Control parameters}
The control parameters for the instances used in the tests are summarized in \Cref{tab:instances:controlParams}.
For each control parameter we report the different values used in the tests. The control parameters were chosen in order to test different
vehicles-to-customers ratios, ranging from $1/2$ to $1/12$, and different absolute values for the number of customers and vehicles. 
In addition, the different configurations of the parameters $\alpha^{FROM}$, $\alpha^{TO}$ and $\alpha^{V}$ yield different configurations of demand (e.g., center to outskirt and outskirt to center) and of the carsharing system (e.g., vehicles located in the center and in the outskirt).
The number of scenarios $|\mtc{S}|$ is arbitrarily set equal to $10$ on all instances.
In \Cref{sec:results:ls} we provide insights on how the L-Shaped method scales with the number of scenarios.
As explained in \Cref{sec:uncertainty} scenarios represent i.i.d. samples of the underlying Gumbel distribution.

\begin{table}
  \caption{Control parameters used to generate the test instances.}
  \label{tab:instances:controlParams}
  \begin{tabular}{p{0.25\linewidth}|p{0.55\linewidth}|p{0.15\linewidth}}
    \toprule
    Parameter& Meaning& Values\\
    \midrule
    $|\mtc{V}|$& the number of vehicles & $50$, $100$, $200$\\
    $|\mtc{K}|$& the number of customers& $200$,$400$,$600$\\
    $\alpha^{FROM}$& initial location of the customers&$0.2$, $0.8$\\
    $\alpha^{TO}$& destination of the customers&$0.2$, $0.8$\\
    $\alpha^{V}$& initial location of the vehicles&$0.2$, $0.8$\\
    Individual profiles &Whether each customer is profiled individually & Yes, No \\
    \bottomrule
  \end{tabular}
\end{table}

\section{Results}\label{sec:results}
This section is divided into three parts. First, we provide implementation details and setup of the experiments in \Cref{sec:results:setup}.
Following, in \Cref{sec:results:ls} we report on the performance of L-Shaped method especially in comparison with a commercial solver.
Finally, in \Cref{sec:results:managerial} we provide an analysis of the solutions and comment on their managerial implications. 

\subsection{Experiments setup}\label{sec:results:setup}
The L-Shaped method and the extensive form the SAAs (see \Cref{sec:saa}) were implemented in Java using the CPLEX $12.10$ libraries.
Particularly, the L-Shaped method was implemented by solving the master problem in a Branch \& Cut framework and adding optimality cuts as lazy constraints at integer nodes.
Unless otherwise specified, we used CPLEX's default parameters both when solving the extensive SAA and when using the L-Shaped method.
This entails, e.g., a target relative optimality gap of $0.01\%$. The only exception, unless otherwise specified, is a time limit of $1800$ seconds.
In the L-Shaped method, relaxation cuts and optimality cuts were applied only at integer nodes throughout the entire tree.
Particularly, relaxation cuts were crucial to the implementation. The performance of the algorithm without relaxation cuts was extremely poor.
Tests were run on machines with $2\times 2.4$ GHz \texttt{AMD Opteron 2431} 6 core CPU and $24$Gb RAM.
We remind the reader that, unless otherwise specified, the SAAs are solved with ten iid samples (scenarios) of random variable $\tilde{\xi}$, see \Cref{sec:uncertainty}.
We stress that the number of scenarios does not represent the number of instances, as in a scenario-analysis procedure.
Rather, by definition, the stochastic program takes into account all scenarios simultaneously. The impact of the number of scenarios on the computational complexity of the problems solved is assessed in \Cref{sec:results:ls}.

\subsection{Analysis of the L-Shaped Method}\label{sec:results:ls}
In the first part of the tests we compared the performance of the L-Shaped method to that of CPLEX for solving the SAA on all configurations of the control parameters in \Cref{tab:instances:controlParams} for which $|\mtc{K}|>|\mtc{V}|$. The scope of our experiments is thus to obtain empirical evidences as to whether, and to what extent, the L-Shaped method scales better than using CPLEX without any decomposition strategy. 
The tables in this section report the optimality gap (\texttt{gap}) and elapsed time (\texttt{t} in seconds) for both CPLEX and the L-Shaped method.
For the L-Shaped method they also report the optimality gap at the root node (\texttt{gapR}) and after $50\%$ of the time limit -- i.e., $50\%$ of $1800$ seconds -- (\texttt{gap50}).
All gaps are expressed as percentages and are calculated as $100 * |\texttt{bestbound}-\texttt{bestinteger}|/|\texttt{bestinteger}|$.
The size of the SAA problem without decomposition is reported in \Cref{sec:app:size}. 

\begin{longtable}{rrrr|rr|rrrr}
  \caption{Comparison of CPLEX and L-Shaped method on the instances with and $\alpha^V=0.2$.}\label{tab:resultsIC02}\\
  \toprule
  & & & & \multicolumn{2}{c}{CPLEX}&\multicolumn{4}{c}{L-Shaped}\\
  $|\mtc{V}|$ & $|\mtc{K}|$  &  $\alpha^{FROM}$ &  $\alpha^{TO}$  & \texttt{gap} & \texttt{t} & \texttt{gap} &    \texttt{gapR} &   \texttt{gap50} & \texttt{t}  \\
  \midrule
    50 &  200 &       0.2 &     0.2 &   0.0088 &   452.47 &  0.6967 &  18.6752 &  0.7466 & 1801.00 \\
  50 &  400 &       0.2 &     0.2 &   4.7820 &  1801.68 & 14.4808 &  74.9725 & 17.2825 & 1800.66 \\
  50 &  600 &       0.2 &     0.2 &        - &        - & 34.9493 &  89.2833 & 39.1895 & 1802.00 \\
 100 &  200 &       0.2 &     0.2 &   0.0084 &   553.92 &  0.0975 &   9.9177 &  0.0975 & 1801.38 \\
 100 &  400 &       0.2 &     0.2 &        - &        - &  2.5186 &  23.4470 &  5.3168 & 1801.07 \\
 100 &  600 &       0.2 &     0.2 &        - &        - & 14.1692 &  35.2426 & 14.8637 & 1821.75 \\
 200 &  400 &       0.2 &     0.2 &        - &        - &  0.0000 &   6.3509 &       - &  156.78 \\
 200 &  600 &       0.2 &     0.2 &        - &        - &  0.0785 &  10.6491 &  0.7863 & 1806.39 \\
\midrule
  50 &  200 &       0.2 &     0.8 &   0.0000 &   134.00 &  0.3381 &  17.4624 &  0.3381 & 1800.03 \\
  50 &  400 &       0.2 &     0.8 &   1.8164 &  1803.60 & 11.2732 &  64.1181 & 13.2796 & 1800.05 \\
  50 &  600 &       0.2 &     0.8 & 137.5269 &  1801.60 & 16.2842 &  91.0324 & 25.1233 & 1800.16 \\
 100 &  200 &       0.2 &     0.8 &   0.0000 &   346.80 &  0.2815 &   7.2972 &  0.2815 & 1800.55 \\
 100 &  400 &       0.2 &     0.8 &        - &        - &  0.6280 &  19.5696 &  0.8692 & 1800.42 \\
 100 &  600 &       0.2 &     0.8 &        - &        - & 11.1184 &  34.2304 & 12.2836 & 1804.76 \\
 200 &  400 &       0.2 &     0.8 &        - &        - &  0.0000 &  41.7344 &       - &  119.53 \\
 200 &  600 &       0.2 &     0.8 &        - &        - &  0.3384 &  11.2541 &  0.3625 & 1820.99 \\
\midrule
  50 &  200 &       0.8 &     0.2 &   0.1916 &  1800.45 &  1.9646 &  23.9553 &  2.2150 & 1800.28 \\
  50 &  400 &       0.8 &     0.2 &  10.3621 &  1806.04 & 19.1742 &  91.0136 & 25.5109 & 1800.03 \\
  50 &  600 &       0.8 &     0.2 &        - &        - & 50.0477 & 109.7915 & 52.8586 & 1800.02 \\
 100 &  200 &       0.8 &     0.2 &   0.0863 &  1803.15 &  0.5169 &  15.6382 &  0.5336 & 1801.16 \\
 100 &  400 &       0.8 &     0.2 &        - &        - &  6.5731 &  27.8935 &  9.8736 & 1801.15 \\
 100 &  600 &       0.8 &     0.2 &        - &        - & 24.8446 &  62.7364 & 24.9827 & 1812.75 \\
 200 &  400 &       0.8 &     0.2 &        - &        - &  0.2654 &  13.1386 &  0.4498 & 1810.83 \\
 200 &  600 &       0.8 &     0.2 &        - &        - &  5.1871 &  20.1720 &  5.7967 & 1821.70 \\
\midrule
  50 &  200 &       0.8 &     0.8 &   0.1403 &  1800.78 &  0.7991 &  22.5586 &  0.8374 & 1800.33 \\
  50 &  400 &       0.8 &     0.8 &   3.1983 &  1803.68 & 15.2703 &  57.2480 & 20.6710 & 1802.52 \\
  50 &  600 &       0.8 &     0.8 &        - &        - & 46.6289 & 106.1663 & 48.8839 & 1804.99 \\
 100 &  200 &       0.8 &     0.8 &   0.2490 &  1801.28 &  0.7255 &  12.8107 &  0.7930 & 1801.43 \\
 100 &  400 &       0.8 &     0.8 &        - &        - &  3.4139 &  27.6540 &  5.4211 & 1803.38 \\
 100 &  600 &       0.8 &     0.8 &        - &        - & 24.8277 &  52.2813 & 26.1444 & 1800.02 \\
 200 &  400 &       0.8 &     0.8 &        - &        - &  0.0963 &   7.2607 &  0.0963 & 1801.18 \\
  200 &  600 &       0.8 &     0.8 &        - &        - &  6.9456 &  17.7706 &  7.3125 & 1846.31 \\
  \midrule
  &    &       &     &  12.1823 &   885.47 &  9.8292 &  38.2289 & 12.1067 & 1701.43 \\
  \bottomrule
\end{longtable}

We start by reporting the results on the default setup, that is, in which customers are not profiled at the individual level (i.e., customers have identical sensitivities to prices, driving time, waiting and walking time, see \Cref{sec:instances:utility}).
\Cref{tab:resultsIC02} and \Cref{tab:resultsIC08} report the results on the instances with $\alpha^V=0.2$ (more vehicles initially located in zones far from the center) and $\alpha^V=0.8$ (more vehicles initially located in central zones), respectively.  Each table reports on a total of $32$ instances, one for each configuration of the control parameters. 
The results in \Cref{tab:resultsIC02} and \Cref{tab:resultsIC08} are rather similar. We observe that CPLEX is a viable alternative only for the smallest instances.
As the number of vehicles grows CPLEX fails to deliver a feasible solution, and runs into memory problems. On the other hand, the L-Shaped method is able to provide a solution to all instances tested, and in many cases it provides a high quality solution, with a rather small optimality gap. We can also observe that, while the optimality gap at the root node is on average much higher than the final optimality gap, the gap after $50\%$ of the allowed time is only a few percentage points higher. This illustrates, that the L-Shaped method may also deliver good solutions in a relatively short time ($15$ minutes).

\begin{longtable}{rrrr|rr|rrrr}
  \caption{Comparison of CPLEX and L-Shaped method on the instances with $\alpha^V=0.8$.}\label{tab:resultsIC08}\\
  \toprule
  & & & & \multicolumn{2}{c}{CPLEX}&\multicolumn{4}{c}{L-Shaped}\\
  $|\mtc{V}|$ & $|\mtc{K}|$  &  $\alpha^{FROM}$ &  $\alpha^{TO}$  & \texttt{gap} & \texttt{t} & \texttt{gap} &    \texttt{gapR} &   \texttt{gap50} & \texttt{t}  \\
  \midrule
    50 &  200 &       0.2 &     0.2 &  0.0059 &   479.49 &  0.5028 & 24.4373 &  0.5411 & 1800.83 \\
  50 &  400 &       0.2 &     0.2 & 16.7040 &  1801.95 & 14.0036 & 48.9993 & 19.7977 & 1800.03 \\
  50 &  600 &       0.2 &     0.2 &       - &        - & 33.2558 & 81.8976 & 37.6588 & 1800.91 \\
 100 &  200 &       0.2 &     0.2 &  0.0091 &  1250.03 &  0.0904 & 15.8445 &  0.2180 & 1800.29 \\
 100 &  400 &       0.2 &     0.2 &       - &        - &  2.0635 & 22.9479 &  6.5261 & 1802.54 \\
 100 &  600 &       0.2 &     0.2 &       - &        - & 16.2161 & 37.5899 & 17.4882 & 1808.28 \\
 200 &  400 &       0.2 &     0.2 &       - &        - &  0.7891 & 12.0835 &  2.2191 & 1809.83 \\
 200 &  600 &       0.2 &     0.2 &       - &        - &  6.3879 & 20.9339 & 30.6798 & 1853.58 \\
\midrule
  50 &  200 &       0.2 &     0.8 &  0.0099 &   246.36 &  0.6165 & 25.7909 &  0.9025 & 1800.87 \\
  50 &  400 &       0.2 &     0.8 &  2.0974 &  1804.45 & 11.7004 & 45.4173 & 15.7910 & 1802.50 \\
  50 &  600 &       0.2 &     0.8 &       - &        - & 31.6544 & 95.5335 & 36.7324 & 1801.59 \\
 100 &  200 &       0.2 &     0.8 &  0.0000 &   469.42 &  0.0291 & 13.9523 &  0.0291 & 1800.33 \\
 100 &  400 &       0.2 &     0.8 &       - &        - &  1.9427 & 25.7279 &  4.6547 & 1805.75 \\
 100 &  600 &       0.2 &     0.8 &       - &        - & 16.5597 & 34.3635 & 17.5988 & 1801.26 \\
 200 &  400 &       0.2 &     0.8 &       - &        - &  0.3674 & 11.6829 &  1.8109 & 1809.27 \\
 200 &  600 &       0.2 &     0.8 &       - &        - &  4.5425 & 19.4970 &  6.6028 & 1807.60 \\
\midrule
  50 &  200 &       0.8 &     0.2 &  0.0000 &    76.22 &  0.0389 & 12.3461 &  0.0389 & 1800.51 \\
  50 &  400 &       0.8 &     0.2 &  4.1437 &  1804.33 & 15.1026 & 67.4814 & 21.7608 & 1803.37 \\
  50 &  600 &       0.8 &     0.2 &       - &        - & 44.8552 & 87.1632 & 48.1871 & 1803.69 \\
 100 &  200 &       0.8 &     0.2 &  0.0000 &   208.51 &  0.0000 &  9.4161 &       - &   26.28 \\
 100 &  400 &       0.8 &     0.2 &       - &        - &  0.2335 & 12.5821 &  0.2726 & 1800.85 \\
 100 &  600 &       0.8 &     0.2 &       - &        - &  7.3760 & 25.0311 &  8.8577 & 1817.52 \\
 200 &  400 &       0.8 &     0.2 &       - &        - &  0.0020 &  8.6414 &       - &  206.88 \\
 200 &  600 &       0.8 &     0.2 &       - &        - &  0.0021 & 11.9967 &       - &  790.72 \\
\midrule
  50 &  200 &       0.8 &     0.8 &  0.0083 &    78.21 &  0.0083 & 16.2539 &       - &   34.46 \\
  50 &  400 &       0.8 &     0.8 & 34.2159 &  1802.63 &  9.6329 & 47.5271 & 12.1571 & 1802.80 \\
  50 &  600 &       0.8 &     0.8 &       - &        - & 37.0417 & 90.4674 & 40.3969 & 1800.35 \\
 100 &  200 &       0.8 &     0.8 &  0.0000 &   222.09 &  0.0000 &  9.2810 &       - &   31.41 \\
 100 &  400 &       0.8 &     0.8 &       - &        - &  0.2077 & 19.7303 &  0.3026 & 1805.44 \\
 100 &  600 &       0.8 &     0.8 &       - &        - &  6.8777 & 27.1833 &  6.9310 & 1811.58 \\
 200 &  400 &       0.8 &     0.8 &       - &        - &  0.0000 & 42.1315 &       - &  125.03 \\
  200 &  600 &       0.8 &     0.8 &       - &        - &  0.0294 &  8.4615 &  0.0306 & 1801.72 \\
  \midrule
   &    &         &       &  4.7662 &   539.14 &  8.1916 & 32.2623 & 13.0072 & 1505.25 \\
\bottomrule
\end{longtable}

A sensible reduction of the optimality gap can be achieved by applying valid inequality \eqref{eq:VI}.
Such valid inequality has the effect of removing some symmetric solutions, as explained in \Cref{sec:ls:efficiency}.
We observed that the addition of \eqref{eq:VI} decreased the average optimality gap from $9.82\%$ to $9.24\%$ on the instances with $\alpha^V=0.2$ and
from $8.19\%$ to $7.10\%$ on the instances with $\alpha^V=0.8$.
The addition of \eqref{eq:VI} resulted particularly beneficial on the instances which yielded the largest optimality gaps reported in \Cref{tab:resultsIC02} and \Cref{tab:resultsIC08}.
All the results on the performance of the L-Shaped method with the addition of \eqref{eq:VI} are reported in \Cref{sec:app:ICplusVI}.

We turn now our attention on the performance of the algorithm when each customer is profiled individually (see \Cref{sec:instances:individual}).
The results in \Cref{tab:resultsDC02} and \Cref{tab:resultsDC08} are obtained with $\alpha^V=0.2$ and $0.8$, respectively.
Valid inequality \eqref{eq:VI} is always added to the models.

\begin{longtable}{rrrr|rr|rrrr}
  \caption{Comparison of CPLEX and L-Shaped method on the instances with $\alpha^V=0.2$ and individual customer profiles.}\label{tab:resultsDC02}\\
  \toprule
  & & & & \multicolumn{2}{c}{CPLEX}&\multicolumn{4}{c}{L-Shaped}\\
  $|\mtc{V}|$ & $|\mtc{K}|$  &  $\alpha^{FROM}$ &  $\alpha^{TO}$  & \texttt{gap} & \texttt{t} & \texttt{gap} &    \texttt{gapR} &   \texttt{gap50} & \texttt{t}  \\
  \midrule
    50 &  200 &       0.2 &     0.2 &  0.0014 &   297.99 &  0.4747 & 21.5553 &  0.4747 & 1801.00 \\
  50 &  400 &       0.2 &     0.2 &  0.4151 &  1804.09 &  6.5981 & 59.3308 &  8.9622 & 1800.13 \\
  50 &  600 &       0.2 &     0.2 &       - &        - & 13.6751 & 52.5981 & 18.9363 & 1800.03 \\
 100 &  200 &       0.2 &     0.2 &  0.0000 &   221.91 &  0.0000 &  8.9488 &       - &   16.47 \\
 100 &  400 &       0.2 &     0.2 &       - &        - &  0.4829 & 26.1941 &  0.8372 & 1801.67 \\
 100 &  600 &       0.2 &     0.2 &       - &        - &  9.2186 & 40.2920 &  9.2848 & 1805.22 \\
 200 &  400 &       0.2 &     0.2 &       - &        - &  0.0000 &  9.9094 &       - &  229.96 \\
 200 &  600 &       0.2 &     0.2 &       - &        - &  0.0571 & 14.7626 &  0.0579 & 1801.39 \\
\midrule
  50 &  200 &       0.2 &     0.8 &  0.0067 &   381.59 &  0.5159 & 20.9177 &  0.5159 & 1800.11 \\
  50 &  400 &       0.2 &     0.8 &  0.3913 &  1800.30 &  3.5957 & 43.2969 &  4.3436 & 1800.03 \\
  50 &  600 &       0.2 &     0.8 &       - &        - & 14.1098 & 63.4863 & 16.2944 & 1802.22 \\
 100 &  200 &       0.2 &     0.8 &  0.0000 &   199.66 &  0.1800 & 10.6250 &  0.1800 & 1800.13 \\
 100 &  400 &       0.2 &     0.8 &       - &        - &  0.1798 & 22.5452 &  0.3266 & 1800.02 \\
 100 &  600 &       0.2 &     0.8 &       - &        - &  9.1325 & 36.2728 &  9.8433 & 1808.65 \\
 200 &  400 &       0.2 &     0.8 &       - &        - &  0.0000 & 59.2478 &       - &  165.06 \\
 200 &  600 &       0.2 &     0.8 &       - &        - &  0.1327 & 15.6153 &  0.4830 & 1800.02 \\
\midrule
  50 &  200 &       0.8 &     0.2 &  0.0000 &   145.84 &  0.0632 & 22.1339 &  0.0632 & 1800.60 \\
  50 &  400 &       0.8 &     0.2 &  0.4584 &  1801.37 &  6.0223 & 57.8394 &  9.1647 & 1800.69 \\
  50 &  600 &       0.8 &     0.2 &  2.9868 &  1804.64 & 20.5953 & 74.1159 & 23.3852 & 1800.03 \\
 100 &  200 &       0.8 &     0.2 &  0.0000 &   374.57 &  0.1444 & 14.1475 &  0.1444 & 1800.12 \\
 100 &  400 &       0.8 &     0.2 &  0.2382 &  1806.73 &  0.9336 & 12.6012 &  1.4826 & 1800.09 \\
 100 &  600 &       0.8 &     0.2 &       - &        - & 14.7882 & 35.0665 & 14.8822 & 1815.08 \\
 200 &  400 &       0.8 &     0.2 &       - &        - &  0.1239 & 11.8485 &  0.1239 & 1800.50 \\
  200 &  600 &       0.8 &     0.2 &       - &        - &  4.6691 & 23.0727 &  5.0362 & 1827.23 \\
  \midrule
  50 &  200 &       0.8 &     0.8 &  0.0000 &   137.42 &  0.2137 & 21.9851 &  0.2137 & 1800.11 \\
  50 &  400 &       0.8 &     0.8 &  0.2161 &  1801.02 &  4.0149 & 53.8589 &  5.2497 & 1801.68 \\
  50 &  600 &       0.8 &     0.8 &       - &        - & 14.4585 & 70.6343 & 20.4297 & 1802.88 \\
 100 &  200 &       0.8 &     0.8 &  0.0000 &   367.32 &  0.0327 & 14.5185 &  0.0327 & 1800.03 \\
 100 &  400 &       0.8 &     0.8 &       - &        - &  2.0523 & 30.1037 &  3.4881 & 1800.07 \\
 100 &  600 &       0.8 &     0.8 &       - &        - & 10.3710 & 34.8283 & 17.2050 & 1814.54 \\
 200 &  400 &       0.8 &     0.8 &       - &        - &  0.0001 & 11.5481 &       - &  170.66 \\
  200 &  600 &       0.8 &     0.8 &       - &        - &  1.7683 & 15.4424 &  2.2644 & 1807.74 \\
  \midrule
   &    &        &      &  0.3367 &   867.80 &  3.8517 & 32.0441 &  6.0126 & 1498.71 \\
  \bottomrule
\end{longtable}

The results illustrated in \Cref{tab:resultsDC02} and \Cref{tab:resultsDC08} are similar to those observed earlier in \Cref{tab:resultsIC02} and \Cref{tab:resultsIC08}.
CPLEX remains a viable option only for the smallest instances, while the L-Shaped method is able to find a solution, in some cases a high quality one, to all instances and to solve a number of them.
We observe also that the average optimality gap appears sensibly lower compared to the results in \Cref{tab:resultsIC02,tab:resultsIC08}.
Also in this case high quality solutions can be obtained already after $15$ minutes.

A pattern in the optimality gaps reported in \Cref{tab:resultsIC02,tab:resultsIC08,tab:resultsDC02,tab:resultsDC08} can be observed.
The optimality gap appears inversely correlated with the vehicles-to-customers ratio.
That is, the instances which yield the largest optimality gaps ($50$ vehicles and $400$ customers, $50$ vehicles and $600$ customers, $100$ vehicles and $600$ customers) are those with the smallest vehicles-to-customers ratios among the instances tested ($1/8$, $1/12$, $1/6$, respectively).
Supposedly, when vehicles are scarce compared to the number of customers, it becomes more challenging for the algorithm to identify, within the $30$ minutes provided,
a relocation and pricing plan which is able to satisfy demand in such a way to yield the highest profit.
On the contrary, as the ratio increases, the model has more freedom to satisfy customers demand, and especially those requests generating the highest revenue.

\begin{longtable}{rrrr|rr|rrrr}
  \caption{Comparison of CPLEX and L-Shaped method on the instances with $\alpha^V=0.8$ and individual customer profiles.}\label{tab:resultsDC08}\\
  \toprule
  & & & & \multicolumn{2}{c}{CPLEX}&\multicolumn{4}{c}{L-Shaped}\\
  $|\mtc{V}|$ & $|\mtc{K}|$  &  $\alpha^{FROM}$ &  $\alpha^{TO}$  & \texttt{gap} & \texttt{t} & \texttt{gap} &    \texttt{gapR} &   \texttt{gap50} & \texttt{t}  \\
  \midrule
    50 &  200 &       0.2 &     0.2 &  0.0044 &   311.01 &  3.8564 &  17.7774 &  4.2676 & 1800.18 \\
  50 &  400 &       0.2 &     0.2 &  1.5001 &  1802.42 & 13.5540 &  60.3835 & 14.4840 & 1800.75 \\
  50 &  600 &       0.2 &     0.2 &       - &        - & 25.1954 & 105.0348 & 27.2104 & 1801.63 \\
 100 &  200 &       0.2 &     0.2 &  0.0099 &  1246.59 &  1.1923 &   7.5390 &  1.3117 & 1801.87 \\
 100 &  400 &       0.2 &     0.2 &       - &        - &  7.0935 &  33.3830 &  9.9477 & 1802.14 \\
 100 &  600 &       0.2 &     0.2 &       - &        - & 14.8632 &  29.7073 & 18.9278 & 1803.56 \\
 200 &  400 &       0.2 &     0.2 &       - &        - &  1.7218 &  22.1296 &  4.9610 & 1805.73 \\
 200 &  600 &       0.2 &     0.2 &       - &        - &  8.6047 &  17.3350 &  8.6047 & 1833.83 \\
\midrule
  50 &  200 &       0.2 &     0.8 &  0.0097 &   375.72 &  1.6477 &  20.1198 &  2.0232 & 1800.04 \\
  50 &  400 &       0.2 &     0.8 &  0.9894 &  1802.12 & 11.4544 &  53.4891 & 13.7308 & 1801.20 \\
  50 &  600 &       0.2 &     0.8 &       - &        - & 20.1773 &  89.9520 & 24.8459 & 1803.77 \\
 100 &  200 &       0.2 &     0.8 &  0.0091 &   671.58 &  0.1778 &   6.9976 &  0.2153 & 1801.05 \\
 100 &  400 &       0.2 &     0.8 &       - &        - &  2.9674 &  32.4136 &  5.1617 & 1800.03 \\
 100 &  600 &       0.2 &     0.8 &       - &        - & 15.1294 &  32.5442 & 16.0408 & 1805.42 \\
 200 &  400 &       0.2 &     0.8 &       - &        - &  4.7056 &  28.3808 &  4.7289 & 1801.26 \\
 200 &  600 &       0.2 &     0.8 &       - &        - & 10.8675 &  26.4194 & 10.8675 & 1800.02 \\
\midrule
  50 &  200 &       0.8 &     0.2 &  0.0003 &    73.80 &  0.0175 &  17.0926 &  0.0175 & 1800.08 \\
  50 &  400 &       0.8 &     0.2 &  0.2689 &  1802.31 &  2.5979 &  32.0731 &  3.4004 & 1800.03 \\
  50 &  600 &       0.8 &     0.2 &       - &        - & 14.2154 &  75.7595 & 17.8579 & 1800.19 \\
 100 &  200 &       0.8 &     0.2 &  0.0000 &   167.21 &  0.0000 &  14.6898 &       - &   22.02 \\
 100 &  400 &       0.8 &     0.2 &       - &        - &  0.2919 &  18.9859 &  0.5424 & 1803.09 \\
 100 &  600 &       0.8 &     0.2 &       - &        - & 10.3980 &  40.7039 & 10.4926 & 1810.93 \\
 200 &  400 &       0.8 &     0.2 &       - &        - &  0.0000 &  61.0074 &       - &  150.43 \\
 200 &  600 &       0.8 &     0.2 &       - &        - &  0.0296 &  12.2654 &  0.0328 & 1804.32 \\
\midrule
  50 &  200 &       0.8 &     0.8 &  0.0068 &    64.66 &  0.2138 &   4.6541 &  0.2138 & 1800.14 \\
  50 &  400 &       0.8 &     0.8 &  0.0302 &  1804.08 &  1.7691 &  42.8840 &  2.3329 & 1800.18 \\
  50 &  600 &       0.8 &     0.8 &       - &        - & 10.0485 &  78.3238 & 13.5323 & 1800.02 \\
 100 &  200 &       0.8 &     0.8 &  0.0000 &   165.74 &  0.0000 &   1.4263 &       - &   23.85 \\
 100 &  400 &       0.8 &     0.8 &       - &        - &  0.1842 &   6.4020 &  0.1842 & 1800.06 \\
 100 &  600 &       0.8 &     0.8 &       - &        - &  7.8771 &  36.8295 &  9.4600 & 1800.01 \\
 200 &  400 &       0.8 &     0.8 &       - &        - &  0.0000 &  61.6141 &       - &  129.12 \\
  200 &  600 &       0.8 &     0.8 &       - &        - &  0.0469 &  13.3682 &  0.0471 & 1802.95 \\
  \midrule
   &    &         &       &  0.2357 &   791.33 &  5.7899 &  33.6217 &  8.1178 & 1594.61 \\
  \bottomrule
\end{longtable}

\Cref{tab:resultsLongICVI02} and \Cref{tab:resultsLongICVI08} report, for the instances with $\alpha^V=0.2$ and $0.8$, respectively, the results obtained by letting the L-Shaped method run for up to $5$ hours ($18 000$ seconds) with a target $1\%$ optimality gap. We observe that the optimality gap goes down from an average of $9.82\%$ to an average of $6.48\%$ for the case with $\alpha^V= 0.2$ and from an average of $8.19\%$ to an average of $5.14\%$ for the case with $\alpha^V= 0.8$. For the case with individual customer profiles we obtain an average optimality gap of $2.48\%$ and $3.31\%$ for the case with $\alpha^V= 0.2$ and $0.8$, respectively. All results on the instances with individual customer profiles are reported in \Cref{app:resultsLongDC}.
These results are possibly of little practical use since, in a business context, solutions are most likely required in much shorter time.
Nevertheless, they show that the model can provide useful bounds that may serve a reference point for example in the development of faster heuristic methods.

\begin{longtable}{rrrr|rrrr}
  \caption{Results of the L-Shaped method with the addition of \cref{eq:VI} on the instances with $\alpha^V=0.2$ with a time limit of $18 000$ seconds and $1\%$ target optimality gap.}\label{tab:resultsLongICVI02}\\
  \toprule
  $|\mtc{V}|$ & $|\mtc{K}|$  &  $\alpha^{FROM}$ &  $\alpha^{TO}$  & \texttt{gap} &    \texttt{gapR} &   \texttt{gap50} & \texttt{t}  \\
  \midrule
    50 &  200 &       0.2 &     0.2 &  0.9890 &  18.6752 &       - &   159.36 \\
  50 &  400 &       0.2 &     0.2 & 10.2092 &  66.5573 & 11.1192 & 18000.12 \\
  50 &  600 &       0.2 &     0.2 & 22.9573 &  87.2440 & 25.9743 & 18003.41 \\
 100 &  200 &       0.2 &     0.2 &  0.6386 &   9.5524 &       - &    29.63 \\
 100 &  400 &       0.2 &     0.2 &  0.9219 &  21.2273 &       - &  3653.12 \\
 100 &  600 &       0.2 &     0.2 &  7.4259 &  35.8381 &  8.8381 & 18000.29 \\
 200 &  400 &       0.2 &     0.2 &  0.0000 &   6.3509 &       - &   157.42 \\
 200 &  600 &       0.2 &     0.2 &  0.8479 &  10.4516 &       - &   717.90 \\
\midrule
  50 &  200 &       0.2 &     0.8 &  0.8965 &  17.6529 &       - &    21.23 \\
  50 &  400 &       0.2 &     0.8 &  6.6519 &  62.4842 &  7.6630 & 18000.62 \\
  50 &  600 &       0.2 &     0.8 & 13.2604 &  90.2582 & 13.8926 & 18000.10 \\
 100 &  200 &       0.2 &     0.8 &  0.3083 &   7.2972 &       - &    23.84 \\
 100 &  400 &       0.2 &     0.8 &  0.9374 &  19.6630 &       - &  1164.27 \\
 100 &  600 &       0.2 &     0.8 &  6.3696 &  34.0425 &  6.9128 & 18008.18 \\
 200 &  400 &       0.2 &     0.8 &  0.0000 &  41.7344 &       - &   128.15 \\
 200 &  600 &       0.2 &     0.8 &  0.7697 &  12.6665 &       - &   705.42 \\
\midrule
  50 &  200 &       0.8 &     0.2 &  1.3911 &  23.9553 &  1.5307 & 18000.20 \\
  50 &  400 &       0.8 &     0.2 & 15.1496 &  85.1665 & 15.6460 & 18000.37 \\
  50 &  600 &       0.8 &     0.2 & 36.5566 & 108.3757 & 37.5433 & 18000.09 \\
 100 &  200 &       0.8 &     0.2 &  0.9165 &  15.4557 &       - &   151.94 \\
 100 &  400 &       0.8 &     0.2 &  2.7552 &  22.7406 &  3.2210 & 18008.25 \\
 100 &  600 &       0.8 &     0.2 & 13.9898 &  42.5987 & 14.5641 & 18010.22 \\
 200 &  400 &       0.8 &     0.2 &  0.4081 &  12.4280 &       - &   332.82 \\
 200 &  600 &       0.8 &     0.2 &  1.2695 &  15.3980 &  2.6620 & 18000.07 \\
\midrule
  50 &  200 &       0.8 &     0.8 &  0.9616 &  22.5871 &       - &   839.48 \\
  50 &  400 &       0.8 &     0.8 & 11.8177 &  55.7576 & 13.0508 & 18001.72 \\
  50 &  600 &       0.8 &     0.8 & 33.2157 & 106.9156 & 34.4212 & 18000.09 \\
 100 &  200 &       0.8 &     0.8 &  0.9309 &  13.3519 &       - &   364.48 \\
 100 &  400 &       0.8 &     0.8 &  1.8436 &  17.1371 &  2.0311 & 18001.11 \\
 100 &  600 &       0.8 &     0.8 & 10.3464 &  46.7222 & 15.3117 & 18004.26 \\
 200 &  400 &       0.8 &     0.8 &  0.5470 &   7.2607 &       - &   296.25 \\
  200 &  600 &       0.8 &     0.8 &  2.3070 &  16.7484 &  2.7073 & 18009.94 \\
  \midrule
     &      &           &         &  6.4872 &  36.0717 & 12.7700 &  9837.32 \\
  \bottomrule
\end{longtable}

\begin{longtable}{rrrr|rrrr}
  \caption{Results of the L-Shaped method with the addition of \cref{eq:VI} on the instances with $\alpha^V=0.8$ with a time limit of $18 000$ seconds and $1\%$ target optimality gap.}\label{tab:resultsLongICVI08}\\
  \toprule
  $|\mtc{V}|$ & $|\mtc{K}|$  &  $\alpha^{FROM}$ &  $\alpha^{TO}$  & \texttt{gap} &    \texttt{gapR} &   \texttt{gap50} & \texttt{t}  \\
  \midrule
    50 &  200 &       0.2 &     0.2 &  0.9052 & 37.8092 &       - &   186.24 \\
  50 &  400 &       0.2 &     0.2 & 10.1669 & 42.1383 & 10.4452 & 18000.24 \\
  50 &  600 &       0.2 &     0.2 & 21.5968 & 77.8433 & 22.0765 & 18000.11 \\
 100 &  200 &       0.2 &     0.2 &  0.9726 &  6.1394 &       - &   106.15 \\
 100 &  400 &       0.2 &     0.2 &  1.2153 & 22.6182 &  1.4967 & 18000.48 \\
 100 &  600 &       0.2 &     0.2 &  5.8809 & 35.6069 &  7.0281 & 18000.17 \\
 200 &  400 &       0.2 &     0.2 &  0.9991 &  5.3451 &       - &  1793.44 \\
 200 &  600 &       0.2 &     0.2 &  0.9964 & 10.6059 &  1.3911 & 14025.40 \\
\midrule
  50 &  200 &       0.2 &     0.8 &  0.9928 & 33.4059 &       - &   904.69 \\
  50 &  400 &       0.2 &     0.8 &  8.4548 & 42.0692 &  8.9256 & 18000.29 \\
  50 &  600 &       0.2 &     0.8 & 20.4095 & 96.0829 & 21.1479 & 18000.08 \\
 100 &  200 &       0.2 &     0.8 &  0.7544 &  4.7786 &       - &    37.44 \\
 100 &  400 &       0.2 &     0.8 &  0.9978 & 20.3964 &  1.2145 & 16084.62 \\
 100 &  600 &       0.2 &     0.8 &  6.5475 & 35.6683 &  6.8735 & 18000.45 \\
 200 &  400 &       0.2 &     0.8 &  0.5994 &  9.3896 &       - &   396.00 \\
 200 &  600 &       0.2 &     0.8 &  0.8804 & 17.6274 &       - &  8150.06 \\
\midrule
  50 &  200 &       0.8 &     0.2 &  0.5667 & 12.5954 &       - &    12.95 \\
  50 &  400 &       0.8 &     0.2 &  9.1978 & 65.9956 &  9.6357 & 18000.12 \\
  50 &  600 &       0.8 &     0.2 & 28.8099 & 82.1348 & 30.6454 & 18002.77 \\
 100 &  200 &       0.8 &     0.2 &  0.4533 &  9.4214 &       - &    22.96 \\
 100 &  400 &       0.8 &     0.2 &  0.9253 &  4.0455 &       - &   204.01 \\
 100 &  600 &       0.8 &     0.2 &  3.7135 & 27.6410 &  4.3096 & 18000.61 \\
 200 &  400 &       0.8 &     0.2 &  0.0020 &  8.6414 &       - &   220.92 \\
 200 &  600 &       0.8 &     0.2 &  0.5837 & 11.5793 &       - &   613.09 \\
\midrule
  50 &  200 &       0.8 &     0.8 &  0.6258 & 14.9683 &       - &    12.69 \\
  50 &  400 &       0.8 &     0.8 &  6.2098 & 44.5269 &  7.6010 & 18002.40 \\
  50 &  600 &       0.8 &     0.8 & 27.1409 & 99.1525 & 27.7092 & 18000.46 \\
 100 &  200 &       0.8 &     0.8 &  0.9110 &  1.3040 &       - &    22.28 \\
 100 &  400 &       0.8 &     0.8 &  0.9539 & 17.2446 &       - &   276.26 \\
 100 &  600 &       0.8 &     0.8 &  2.2269 & 13.2578 &  2.6462 & 18005.33 \\
 200 &  400 &       0.8 &     0.8 &  0.0000 & 42.1315 &       - &   124.68 \\
  200 &  600 &       0.8 &     0.8 &  0.0754 &  8.4615 &       - &   419.24 \\
  \midrule
     &      &           &         &  5.1489 & 30.0196 & 10.8764 &  8675.83 \\
  \bottomrule
\end{longtable}

\begin{figure}
  \centering
  \begin{subfigure}[b]{0.45\textwidth}
    \centering
    \includegraphics[width=\textwidth]{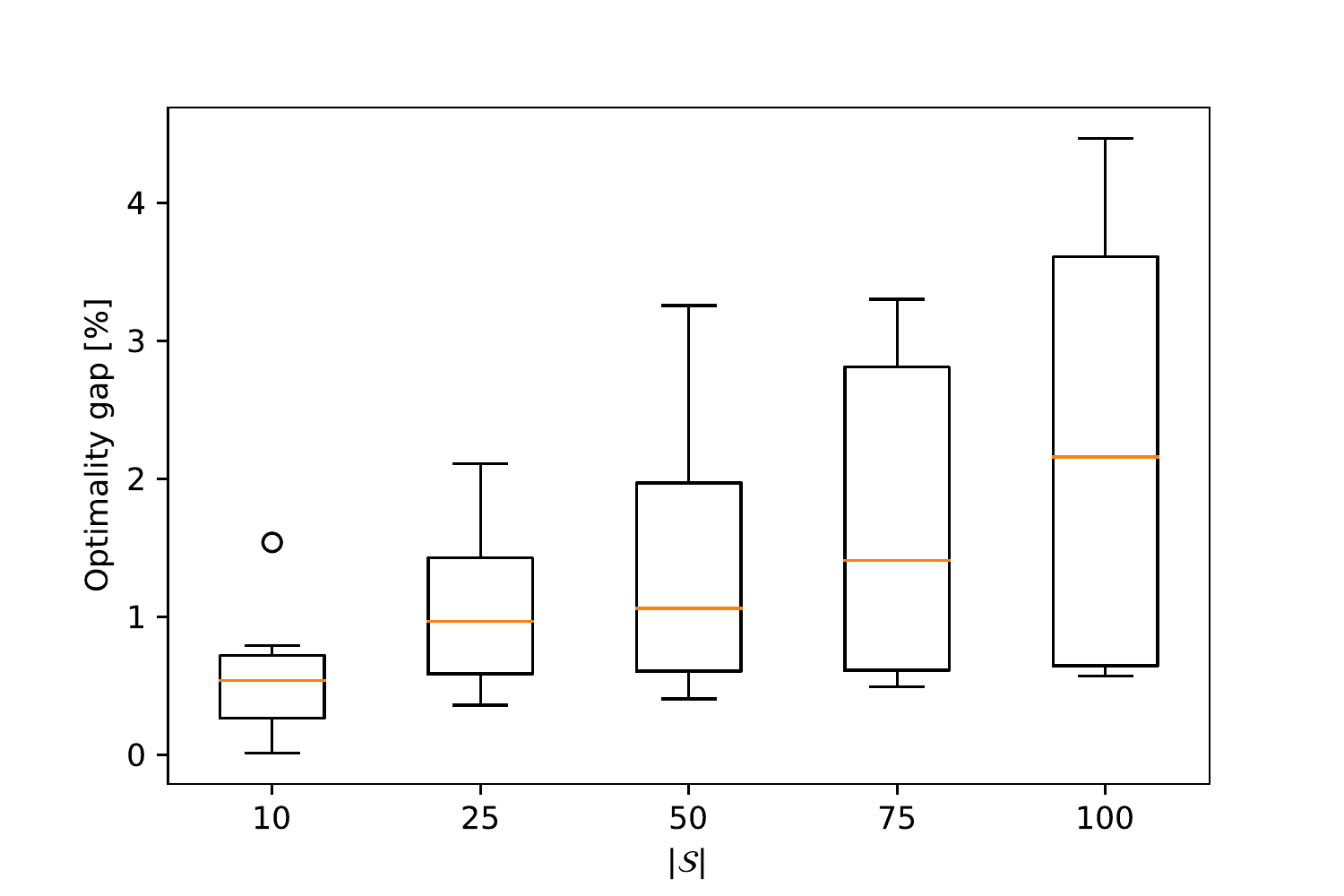}
    \caption{L-Shaped method after $30$ minutes.}
    \label{fig:scenarios_50_200_ls}
  \end{subfigure}
  ~
  \begin{subfigure}[b]{0.45\textwidth}
    \centering
    \includegraphics[width=\textwidth]{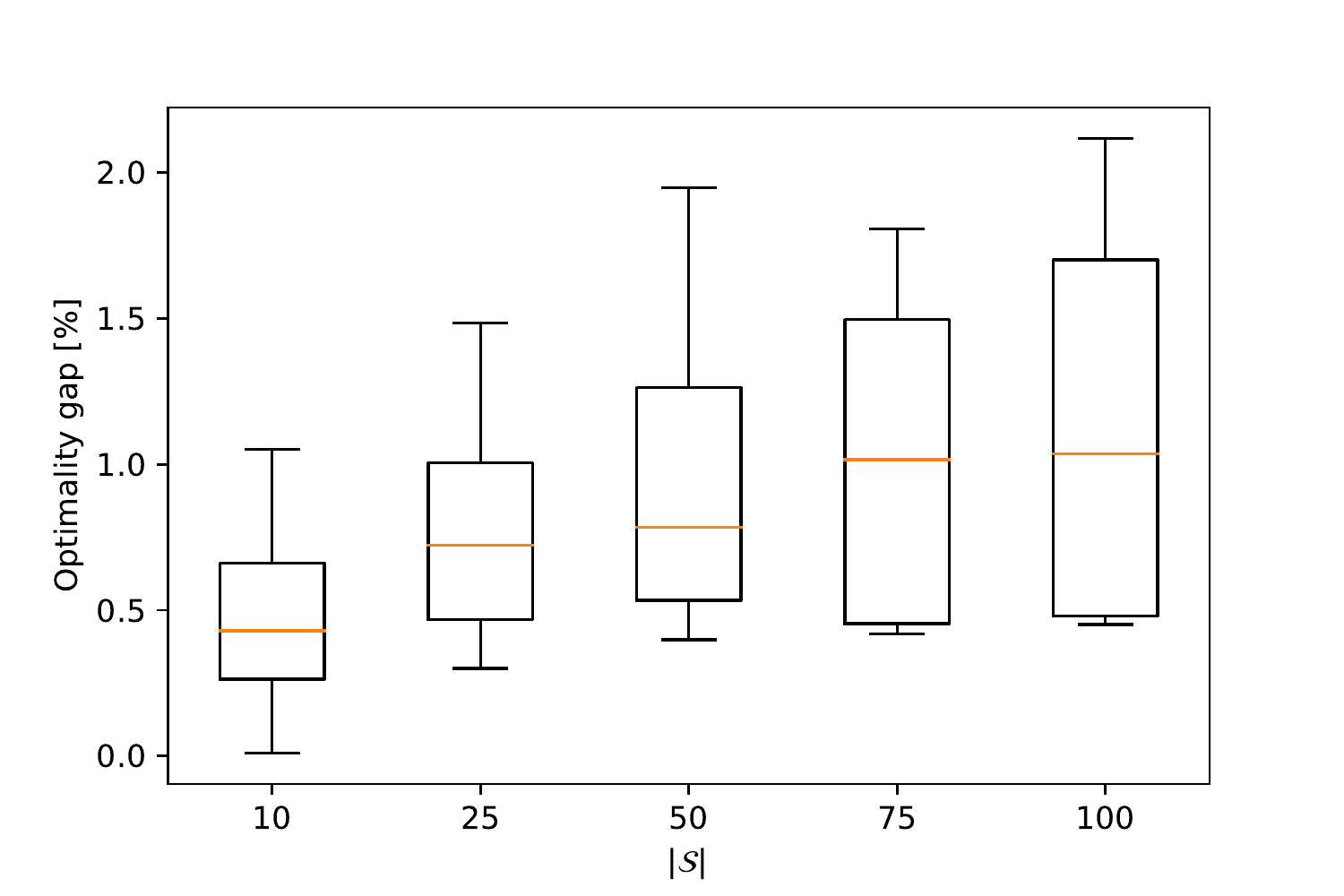}
    \caption{L-Shaped method after $5$ hours. }
    \label{fig:scenarios_50_200_5h}
  \end{subfigure}\\
    \begin{subfigure}[b]{0.45\textwidth}
    \centering
    \includegraphics[width=\textwidth]{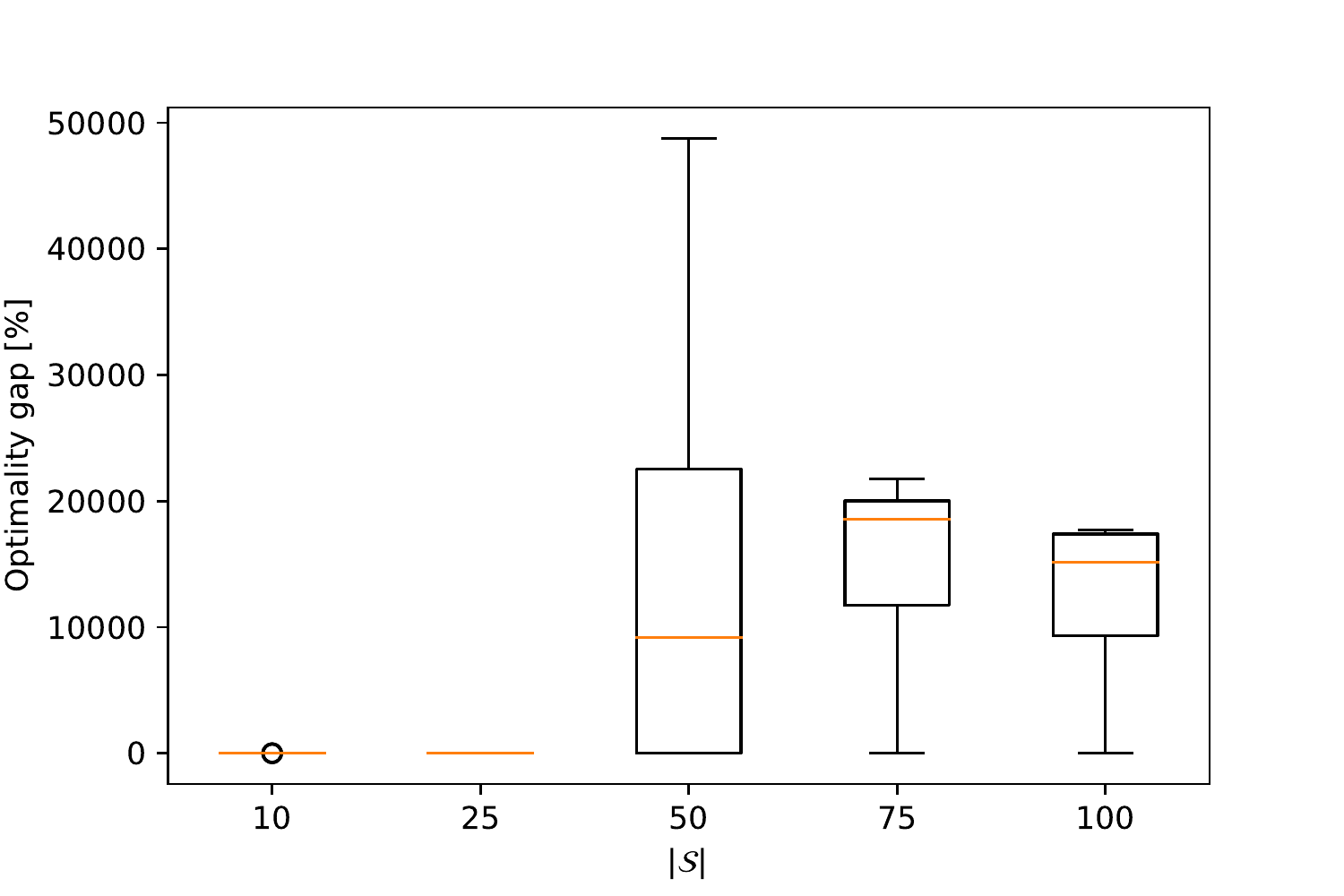}
    \caption{CPLEX after $30$ minutes.}
    \label{fig:scenarios_50_200_cplex}
  \end{subfigure}
  ~
  \begin{subfigure}[b]{0.45\textwidth}
    \centering
    \includegraphics[width=\textwidth]{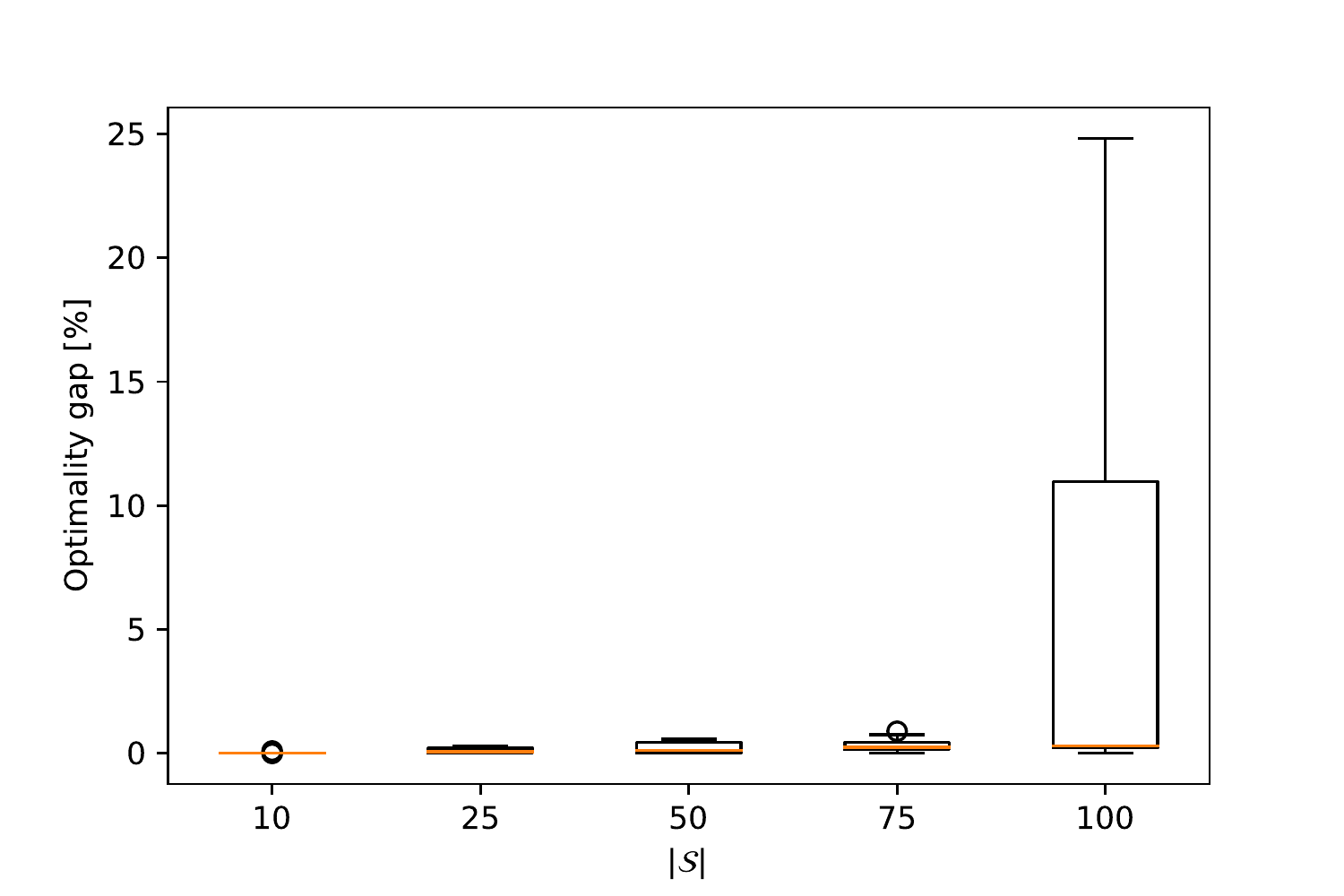}
    \caption{CPLEX after $5$ hours. Plot obtained removing an outlier with an optimality gap greater than $1000$\%.}
    \label{fig:scenarios_50_200_cplex_5h}
  \end{subfigure} 
  \caption{Optimality gap obtained when using the L-Shaped method and CPLEX for different sample sizes $|\mtc{S}|$ on the instances with $|\mtc{V}|=50$ and $|\mtc{K}|=200$. The optimality gap is computed as the average over all combinations of $\alpha^V$, $\alpha^{FROM}$ and $\alpha^{TO}$ assuming identical customers and no additional valid indequality. }
  \label{fig:scenarios_50_200}
\end{figure}

\begin{figure}
  \centering
  \begin{subfigure}[b]{0.45\textwidth}
    \centering
    \includegraphics[width=\textwidth]{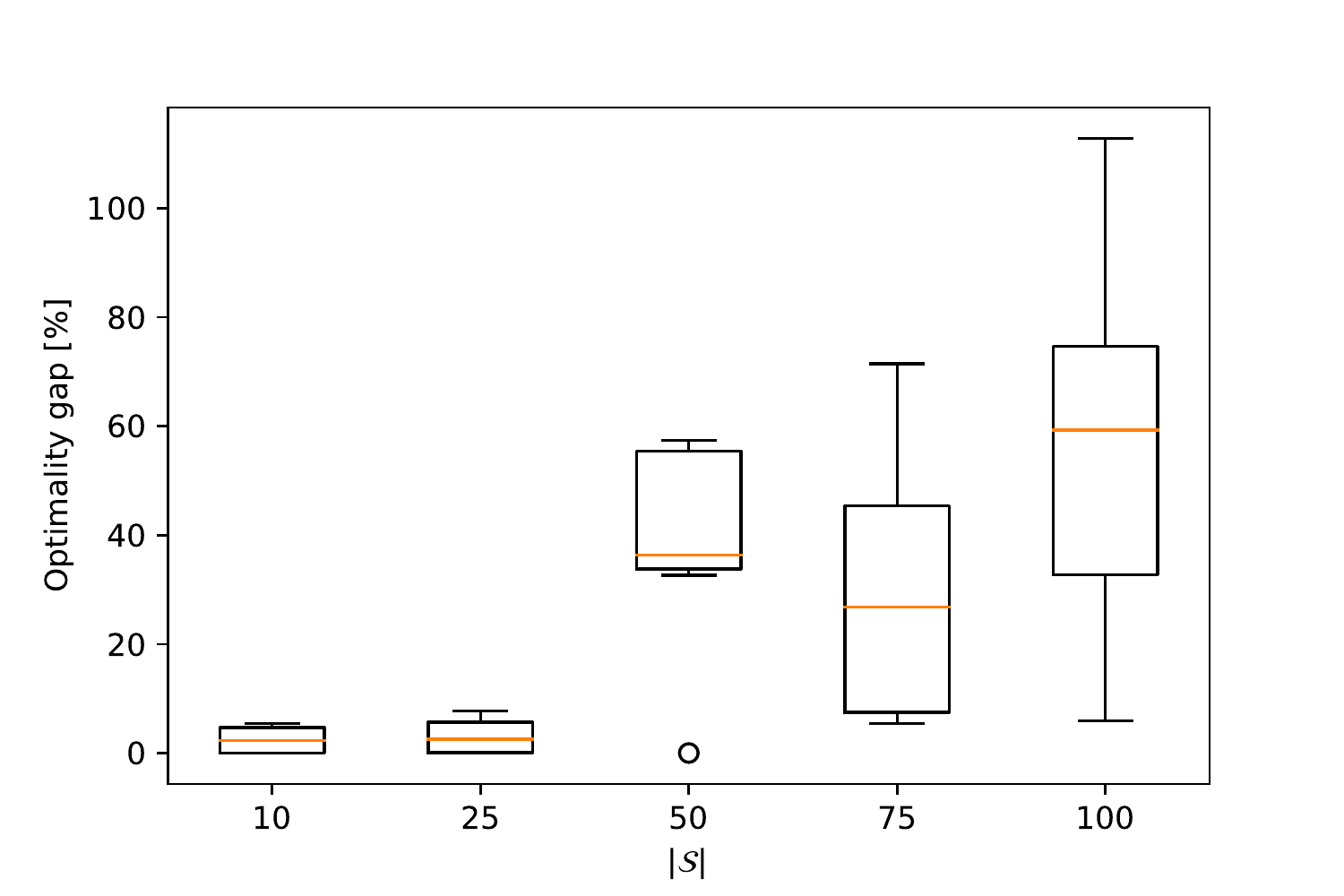}
    \caption{L-Shaped method after $30$ minutes.}
    \label{fig:scenarios_200_60030min}
  \end{subfigure}
  ~
  \begin{subfigure}[b]{0.45\textwidth}
    \centering
    \includegraphics[width=\textwidth]{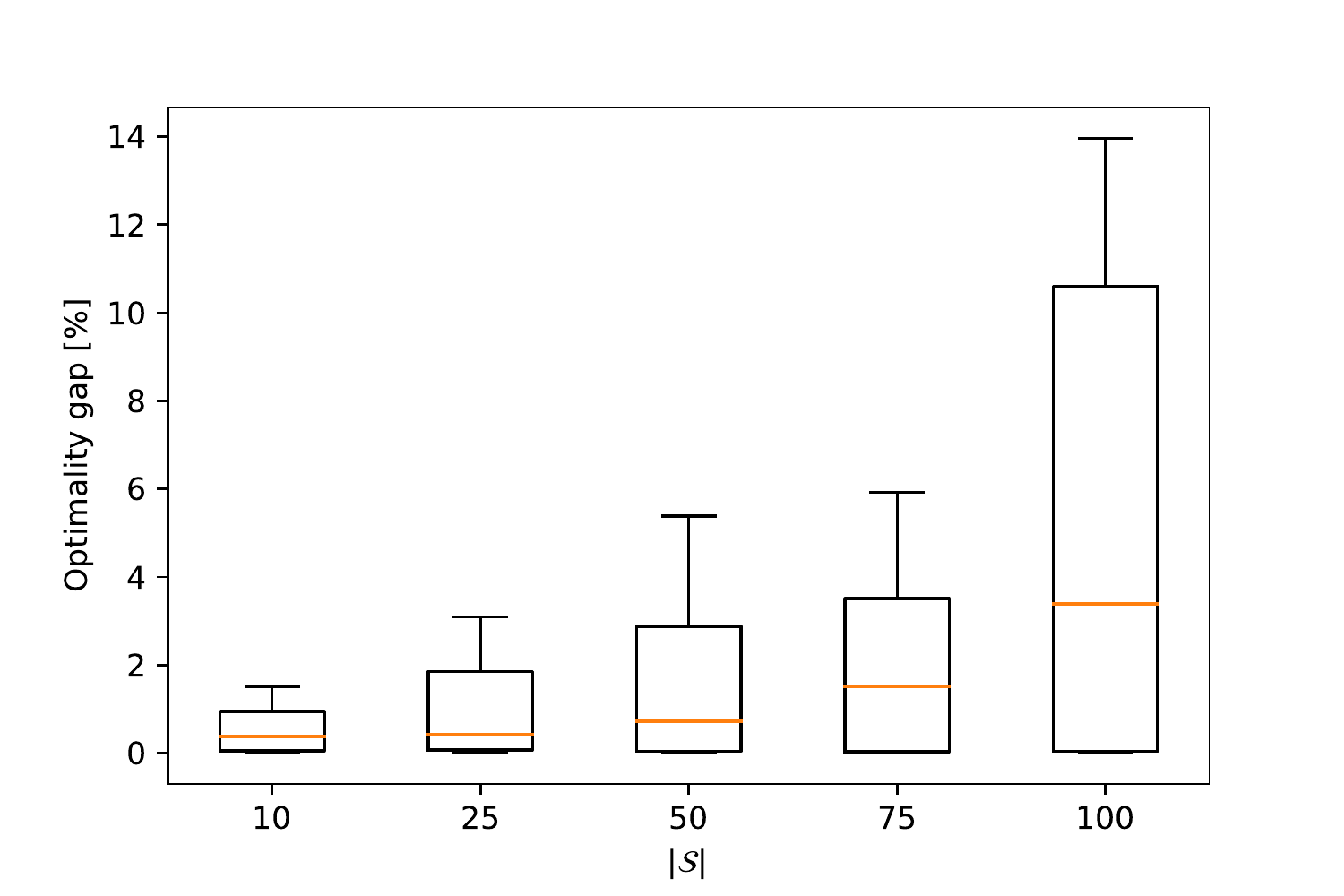}
    \caption{L-Shaped method after $5$ hours.}
    \label{fig:scenarios_200_6005h}
  \end{subfigure}
  \caption{Optimality gap obtained when using the L-Shaped method for different sample sizes $|\mtc{S}|$ on the instances with $|\mtc{V}|=200$ and $|\mtc{K}|=600$. The optimality gap is computed as the average over all combinations of $\alpha^V$, $\alpha^{FROM}$ and $\alpha^{TO}$ assuming identical customers and no additional valid indequality. }
  \label{fig:scenarios_200_600}
\end{figure}

Finally, we report on the performance of the L-Shaped method and CPLEX  as the number of scenarios (sample size) increases.
The results reported above in this section have been obtained by arbitrarily using ten scenarios to approximate the underlying continuous random variable.
\Cref{fig:scenarios_50_200,fig:scenarios_200_600} report the optimality gap obtained with the L-Shaped method and CPLEX as the sample size increases, for the smallest ($|\mtc{V}|=50$, $|\mtc{K}|=200$) and largest ($|\mtc{V}|=200$, $|\mtc{K}|=600$) instances, respectively, with identical customers. Particularly, \Cref{fig:scenarios_50_200_ls} and \Cref{fig:scenarios_50_200_5h} report the gap of the L-Shaped method after $30$ minutes and $5$ hours, respectively, while \Cref{fig:scenarios_50_200_cplex} and \Cref{fig:scenarios_50_200_cplex_5h} report the gap of CPLEX after $30$ minutes and $5$ hours, respectively.
Similarly, \Cref{fig:scenarios_200_600} reports the gaps of the L-Shaped method after after $30$ minutes and $5$ hours on the largest instance. Tests on the largest instances were conducted only for the L-Shaped method as it is already evident from \Cref{tab:resultsIC02,tab:resultsIC08} that, on those instances, CPLEX fails to deliver solutions already with a sample of ten scenarios.

As intuition suggests, when using the L-Shaped method, the optimality gap grows with the number of scenarios, both on the smallest (\Cref{fig:scenarios_50_200}) and the largest (\Cref{fig:scenarios_200_600}) instances. For the smallest instances, with a time limit of $30$ minutes the growth appears mild and, with $100$ scenarios, the optimality gap remains in the neighborhood of $4\%$ in the worst case, with an average optimality gap in the neighborhood of $2\%$ (see \Cref{fig:scenarios_50_200_ls}). After $5$ hours, the L-Shaped method drops the optimality gap even further, with a worst case gap in the neighborhoods of $2$\%. This gives room for a more dense approximation of the uncertainty, i.e., a larger sample size, compared to the $10$ scenarios used in our previous tests.
The performance of CPLEX on the same instances is dramatically worse (see \Cref{fig:scenarios_50_200_cplex}). Given a $30$-minute time limit, on the smallest sample sizes, CPLEX outperforms the L-Shaped method, but as the sample size grows to $50$ or higher the solver's optimality gaps are orders of magnitude higher than those of the L-Shaped method (compare \Cref{fig:scenarios_50_200_cplex} and \Cref{fig:scenarios_50_200_ls}). The performance of CPLEX improves with a $5$ hours time limit (see \Cref{fig:scenarios_50_200_cplex_5h}), though performing much worse than the L-Shaped method at least for the largest sample sizes (compare \Cref{fig:scenarios_50_200_cplex_5h} and \Cref{fig:scenarios_50_200_5h}).
On the largest instances, with a $30$-minute time limit, the gap growth for the L-Shaped method remains limited up to a sample size of $25$, but grows dramatically with a larger sample size, see \Cref{fig:scenarios_200_60030min}.
Also the variance of the optimality gap grows with the sample size, limiting the reliability of the method for large numbers of scenarios. Nevertheless, with a longer time limit the L-Shaped method is able to reduce the optimality gap approximately ten times on the largest sample sizes, see \Cref{fig:scenarios_200_6005h}.

Summarizing, the results in \Cref{fig:scenarios_50_200,fig:scenarios_200_600} illustrate that, for small instances, a solution time of $30$ minutes is sufficient to obtain high quality solutions also with a more extensive approximation of the uncertainty. However, especially on the largest instances, a $30$-minute time limit might result too small to accommodate for a better description of the uncertainty. Nevertheless, allowing the L-Shaped method to run for a longer time (e.g., $5$ hours) can yield substantial reductions of the optimaility gap. In any case, \Cref{fig:scenarios_50_200} illustrates that the L-Shaped method scales significantly better than CPLEX as the sample size increases.

The observed performance of the L-Shaped method allows us a conclusive reflection on the envisaged usage of the method.
In \Cref{tab:resultsIC02,tab:resultsIC08,tab:resultsDC02,tab:resultsDC08} we have let both the algorithm and CPLEX run for $30$ minutes
and observed that our algorithm scales significantly better. Therefore, depending on the practical operating needs of the the CSO, the algorithm might already provide a practice-ready tool.
That is, if the CSO is able to wait for a solution for $30$ minutes, our tests provide empirical evidences that the algorithm delivers a solution and often a high-quality one,
with the vehicles-to-customers ratio being a strong driver of the quality of the solution obtained.
\Cref{tab:resultsIC02,tab:resultsIC08,tab:resultsDC02,tab:resultsDC08}
also illustrate that the algorithm was able to find good solutions after $15$ minutes ($50$\% of the solution time, see column \texttt{gap50}).
Therefore, if the CSO has tighter time requirements, a potential strategy is to terminate
the algorithm earlier, knowing that this implies giving up something in terms of quality of the solution.
However, our tests show that the optimality gap after $15$ minutes is typical not dramatically higher than the
optimality gap obtained after $30$ minutes.
Nevertheless, there might arise situations in which waiting for a solution for $30$ or even $15$ minutes might be impractical, e.g., if the demand landscape changes more
frequently and relocation and pricing plans are required more often.
In these cases, the proposed algorithm might be proven inefficient and one might have to consider developing faster heuristic algorithms.
If this is the case, the performance of the L-Shaped method, and the dual bounds it delivers, provide a reliable benchmark. An example is provided in \Cref{app:ILS} where we test a simple Iterated Local Search and assess its performances using the bounds provided by the L-Shaped method. For the smallest instances the heuristic is able to provide primal solutions of quality comparable or even better than the L-Shaped method.
Nevertheless, the quality of the solutions delivered drops significantly as the size of the instance increase, indicating that further refinement is needed. 
Finally, regardless of the solution time, the algorithm proposed may be used by the CSO to obtain solutions
that allow them to support managerial choices or analyze policy implications, e.g., subsidies, plans to expand the fleet or to hire additional staff for relocations activities.
Some insights on the impact of a pricing scheme on profits and demand are provided in \Cref{sec:results:managerial}.

\subsection{Analysis of the solutions}\label{sec:results:managerial}
In this section we present some evidences based on the analysis of the solutions obtained by the proposed model.
The analysis was performed using the default configuration, i.e., without individual customer profiles, as we believe it is a more realistic configuration to achieve by CSOs.

The analysis in this section is based on the results obtained on the instances with the largest number of customers ($600$)
and with different distributions of vehicles and customers, namely
\begin{description}
\item[D1] Vehicles in the outskirt and demand from center to outskirt ($\alpha^V=0.2$, $\alpha^{FROM}=0.8$, $\alpha^{TO}=0.2$)
\item[D2] Vehicles in the center and demand from outskirt to center ($\alpha^V=0.8$, $\alpha^{FROM}=0.2$, $\alpha^{TO}=0.8$)
\item[D3] Vehicles in the center and demand from center to outskirt ($\alpha^V=0.8$, $\alpha^{FROM}=0.8$, $\alpha^{TO}=0.2$)
\item[D4] Vehicles in the outskirt and demand from outskirt to center ($\alpha^V=0.2$, $\alpha^{FROM}=0.2$, $\alpha^{TO}=0.8$)
\end{description}
The number of vehicles is set either to $50$ or to $200$, corresponding to vehicles-to-customers ratios of $1/12$ and $1/3$, respectively.
 It should be noted that, as pointed out by e.g., \cite{HuaHK18,JiaRD19}, the number of available cars in a zone influences customer demand.
In our approach, in which customers are considered at the individual level, the connection between the number of available cars and demand is handled jointly by the utility function and, especially, by the optimization model. 
That is, there is a potential demand, made of the users which, according to their utility function, would choose carsharing, if available, at the drop-off fee level set in the first stage,
and there is a realized demand (i.e., actual rentals), which takes into account that not all potential customers may find an available car.
This is done through constraints \eqref{eq:taken} that state that an available car is taken by the first customer arriving at the car.
It should be further clarified that, in our experiments, we assume that customers do not wait for more cars to become available, i.e., the waiting time is always set to zero in our instances, corresponding to saying that, if no car is available, the customer will immediately choose another transport service.
Indeed, in some real-life carsharing services, some waiting time could be taken into account. That is, when there are no cars available in the zone, some user might decide to wait until a car is returned. However, we believe the assumption that customers do not wait is the most appropriate especially in a free-floating one-way service, where both customers and the CSO have limited information on whether and when a new car will be returned close to the user.

We start by presenting the effect of pricing strategies on profits and relocations.
For each distribution of customers and vehicles, we solved two configurations of the model.
In the first configuration the prices were optimally set by model \eqref{eq:1S}.
In the second configuration the drop-off fee was set to $0$ (the average of the drop-off fees considered) everywhere, corresponding to a situation in which the CSO applies only a per-minute fee and does not adjust prices with respect to the origin and destination and to the time of the day.

\begin{longtable}{ll|cc}
  \caption{Comparison of the solutions with and without dynamic pricing on the instances with $50$ vehicles and $600$ customers.}\label{tab:results:solutions50600}\\
  \toprule
  Distribution&Metric & With dynamic pricing & Without dynamic pricing \\
  \midrule
  D1          &Expected Profit [\%] &100 &81.78\\
              &\% of vehicles Relocated & 26.0 & 10.0\\
              &Min $|\mtc{R}(\xi)|$&167 &80 \\
              &Max $|\mtc{R}(\xi)|$&195 &107\\
              &Expected \% Requests satisfied & 24&42\\
  \midrule
  D2          &Expected Profit [\%] &100 & 66.06\\
              &\% of vehicles Relocated & 22.0 & 2.0\\
              &Min $|\mtc{R}(\xi)|$&168 &81 \\
              &Max $|\mtc{R}(\xi)|$&187 &105\\
              &Expected \% Requests satisfied & 26&49\\
  \midrule
  D3          &Expected Profit [\%] &100 & 65.05\\
              &\% of vehicles Relocated & 18.0 & 6.0\\
              &Min $|\mtc{R}(\xi)|$&167 &80 \\
              &Max $|\mtc{R}(\xi)|$&195 &107\\
              &Expected \% Requests satisfied & 26&49\\
  \midrule
  D4          &Expected Profit [\%] &100 & 66.36\\
              &\% of vehicles Relocated & 10.0 & 0.0\\
              &Min $|\mtc{R}(\xi)|$&168 &81 \\
              &Max $|\mtc{R}(\xi)|$&187 &105\\
              &Expected \% Requests satisfied & 26&48\\
  \bottomrule
\end{longtable}

\Cref{tab:results:solutions50600,tab:results:solutions200600} report a number of solution statistics for the case with $50$ and $200$ customers, respectively.
Expected profits for the case without dynamic pricing are reported as a percentage of the expected profits with dynamic pricing.
In both the case with $50$ and $200$ vehicles we observe that the expected profit without pricing is approximately $65$ to $80\%$ of the expected profit obtained by adjusting prices.
The main driver of the higher expected profit generated by a pricing strategy is the higher number of requests generated, approximately double both in the case with $50$ and in the case with $200$ vehicles. That is, by adjusting prices the CSO is able to attract significantly more demand and increase competition.

\begin{longtable}{ll|cc}
  \caption{Comparison of the solutions with and without dynamic pricing on the instances with $200$ vehicles and $600$ customers.}\label{tab:results:solutions200600}\\
  \toprule
  Distribution&Metric & With dynamic pricing [\%] & Without dynamic pricing [\%]\\
  \midrule
  D1          &Expected Profit [\%] &100 & 70.12\\
              &\% of vehicles Relocated & 0.5 & 1.5\\
              &Min $|\mtc{R}(\xi)|$&167 &80 \\
              &Max $|\mtc{R}(\xi)|$&195 &107\\
              &Expected \% Requests satisfied & 53&90\\
  \midrule
  D2          &Expected Profit [\%] &100 & 70.87\\
              &\% of vehicles Relocated & 1.5 & 0.5\\
              &Min $|\mtc{R}(\xi)|$&168 &82 \\
              &Max $|\mtc{R}(\xi)|$&187 &99\\
              &Expected \% Requests satisfied & 51&91\\
  \midrule
  D3          &Expected Profit [\%] &100 & 73.94\\
              &\% of vehicles Relocated & 0 & 0\\
              &Min $|\mtc{R}(\xi)|$&167 &80 \\
              &Max $|\mtc{R}(\xi)|$&195 &107\\
              &Expected \% Requests satisfied & 56&100\\
  \midrule
  D4          &Expected Profit [\%] &100 & 71.18\\
              &\% of vehicles Relocated & 0 & 0\\
              &Min $|\mtc{R}(\xi)|$&168 &82 \\
              &Max $|\mtc{R}(\xi)|$&187 &99\\
              &Expected \% Requests satisfied & 56&94\\
  \bottomrule
\end{longtable}

\begin{figure}
  \centering
  \begin{subfigure}[b]{0.48\textwidth}
    \centering
    \includegraphics[width=\textwidth]{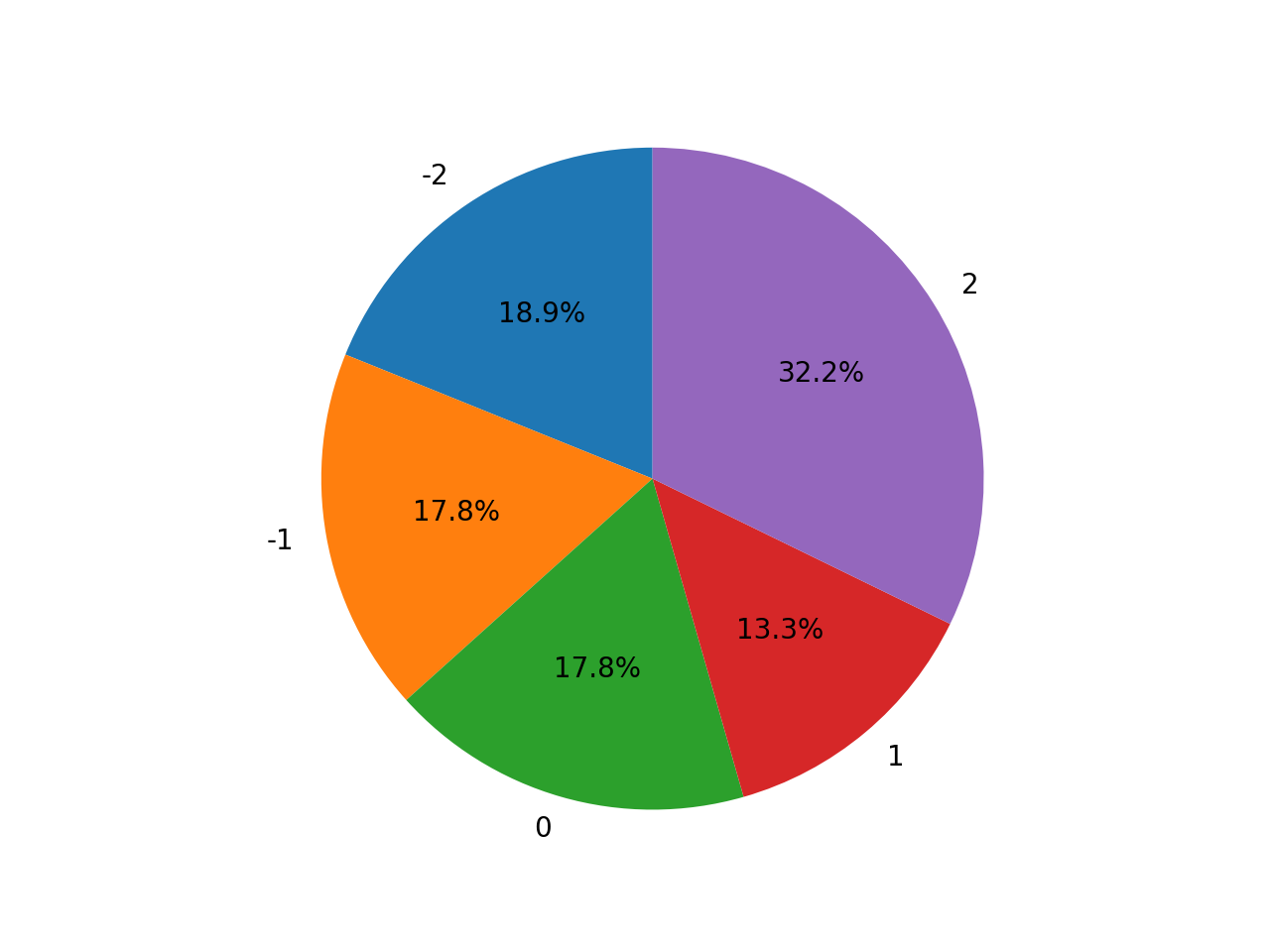}
    \caption{Distribution $D1$}
    \label{fig:200d1}
  \end{subfigure}
  \hfill
  \begin{subfigure}[b]{0.48\textwidth}
    \centering
    \includegraphics[width=\textwidth]{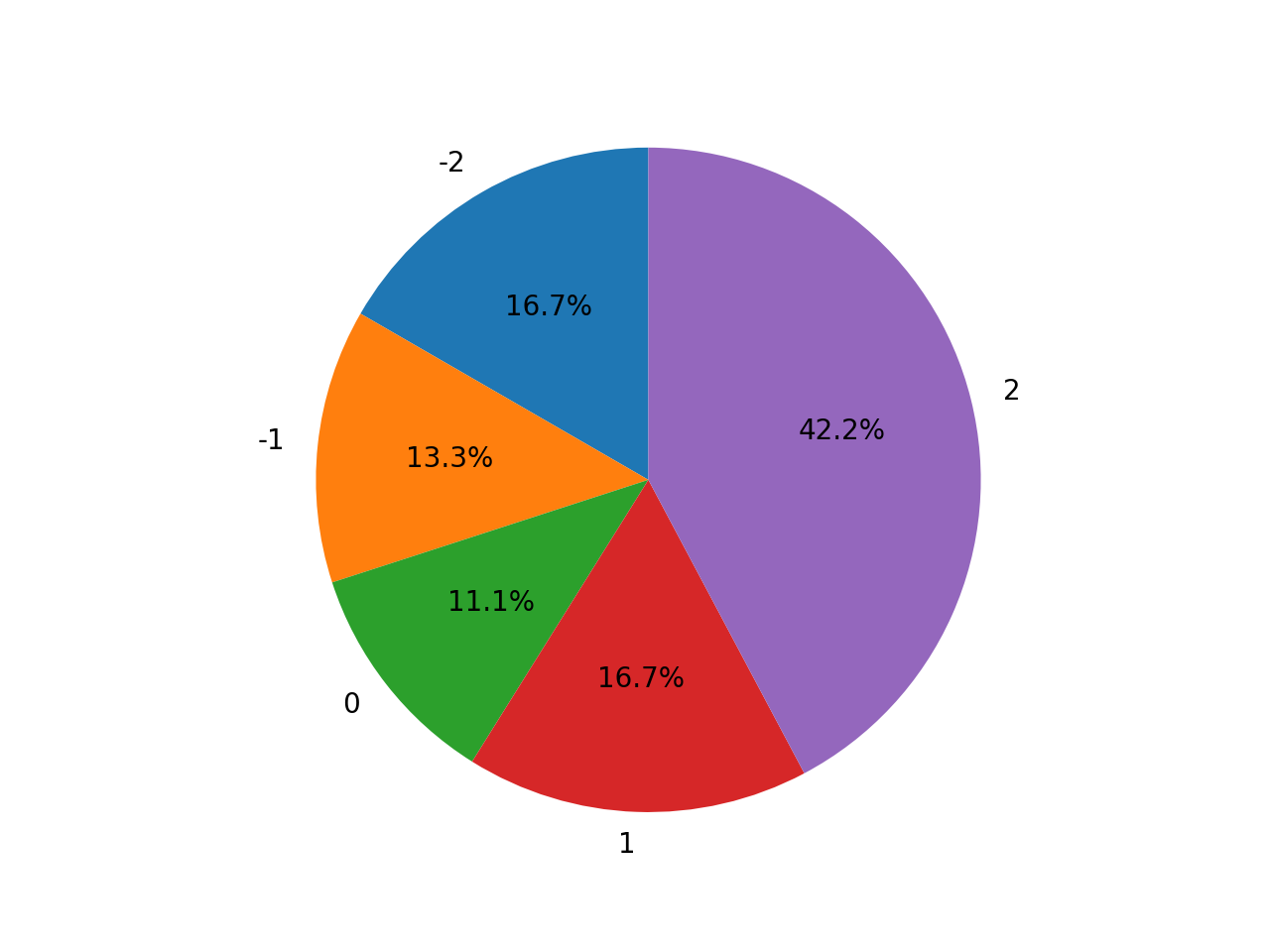}
    \caption{Distribution $D2$}
    \label{fig:200d2}
  \end{subfigure}
  \\
  \begin{subfigure}[b]{0.48\textwidth}
    \centering
    \includegraphics[width=\textwidth]{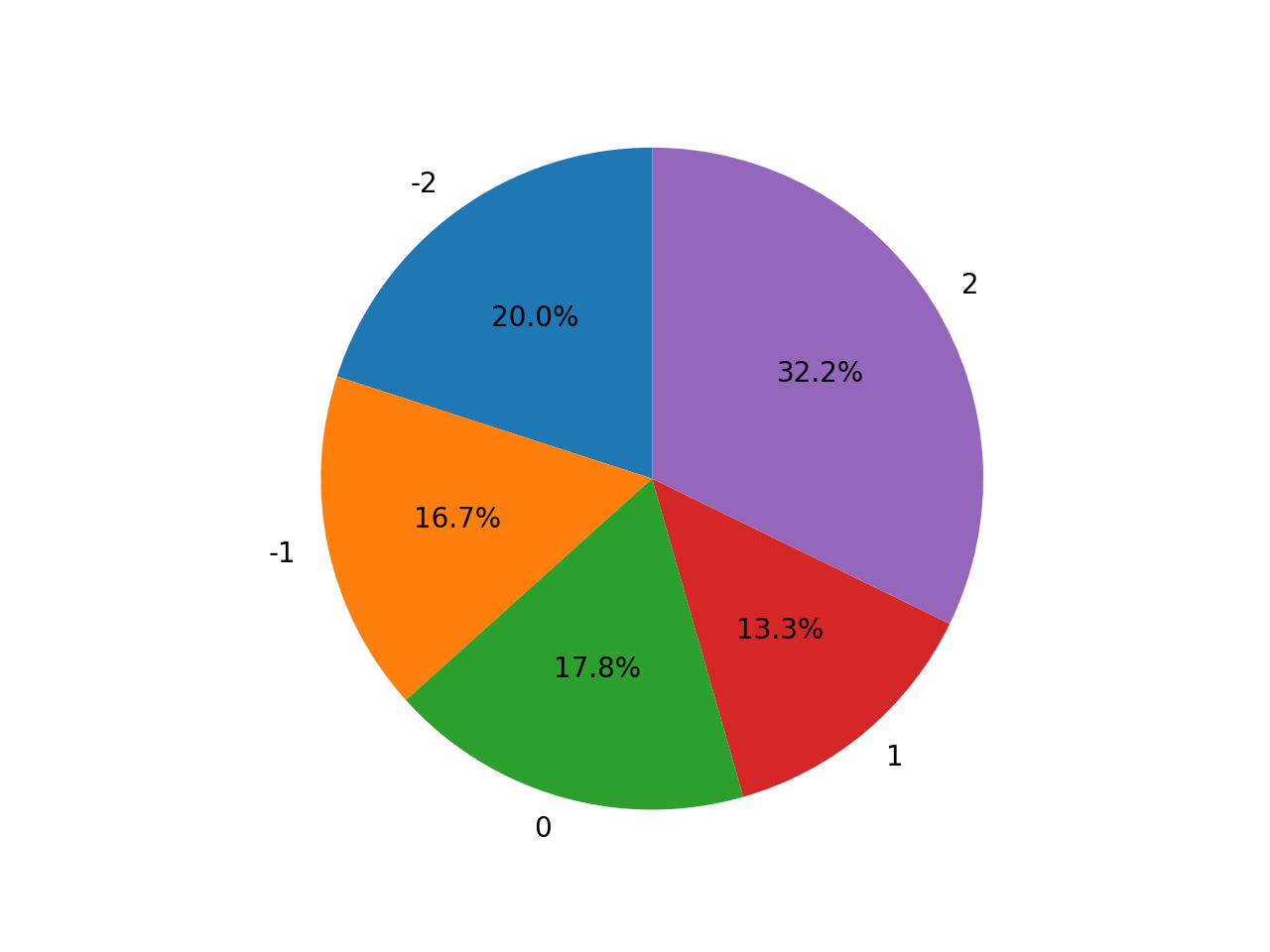}
    \caption{Distribution $D3$}
    \label{fig:200d3}
  \end{subfigure}
  \hfill
  \begin{subfigure}[b]{0.48\textwidth}
    \centering
    \includegraphics[width=\textwidth]{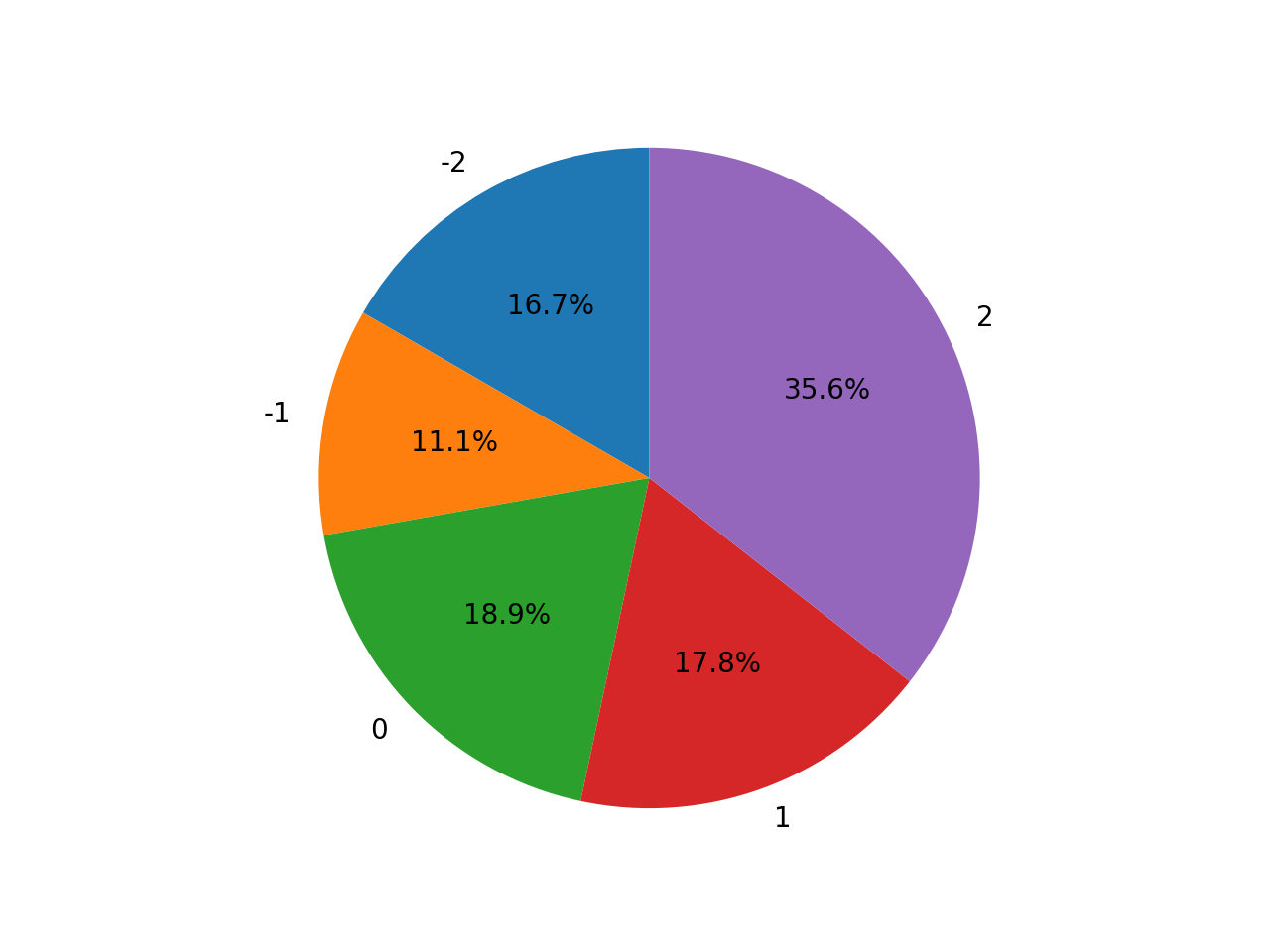}
    \caption{Distribution $D4$}
    \label{fig:200d4}
  \end{subfigure}
  \caption{Drop-off fees adopted in the instances with $200$ vehicles and $600$ customers. Drop-off fees are expressed in Euro.}
  \label{fig:fares_distribuition200}
\end{figure}

In the case with $200$ vehicles and without dynamic pricing, the CSO is able to satisfy the great majority of the requests (more than $90\%$ -- see \Cref{tab:results:solutions200600}) performing very few, if any, relocations -- we will return to this point later.
However, in this case the CSO is able to attract less than (approximately) $1/6$ of the customers. With dynamic pricing, the CSO is able to attract close to $1/3$ of the customers and serve slightly more than $50\%$ of the requests. In total, the number of rentals are approximately the same, with and without dynamic pricing. However, by adjusting prices the CSO is able to increase the revenue. In fact, as shown in \Cref{fig:fares_distribuition200}, the most used drop-off fee is the highest ($2$ Euro), illustrating that the CSO is able to exploit the higher willingness to pay of some customers.

\begin{figure}
  \centering
  \begin{subfigure}[b]{0.48\textwidth}
    \centering
    \includegraphics[width=\textwidth]{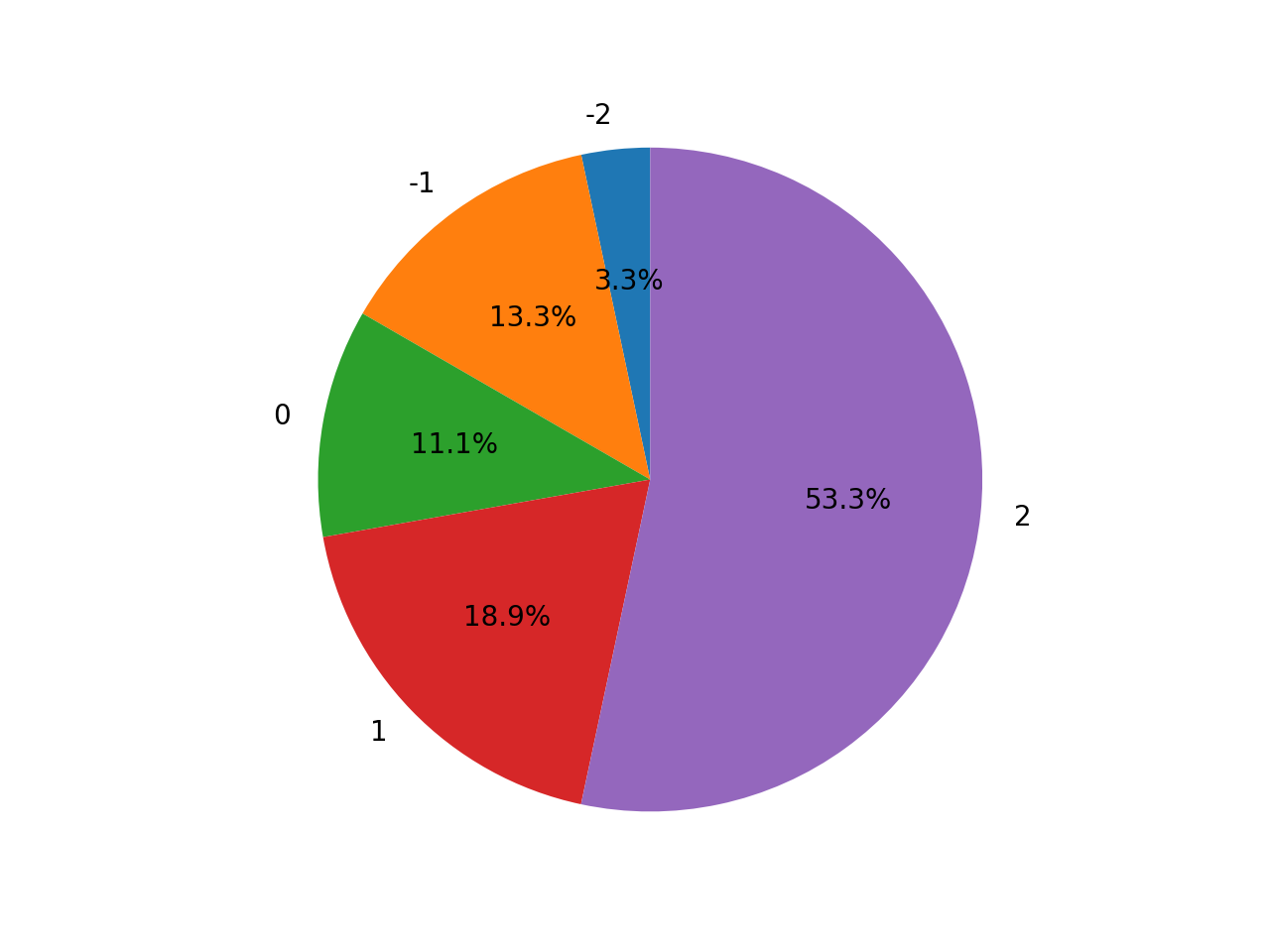}
    \caption{Distribution $D1$}
    \label{fig:d1}
  \end{subfigure}
  \hfill
  \begin{subfigure}[b]{0.48\textwidth}
    \centering
    \includegraphics[width=\textwidth]{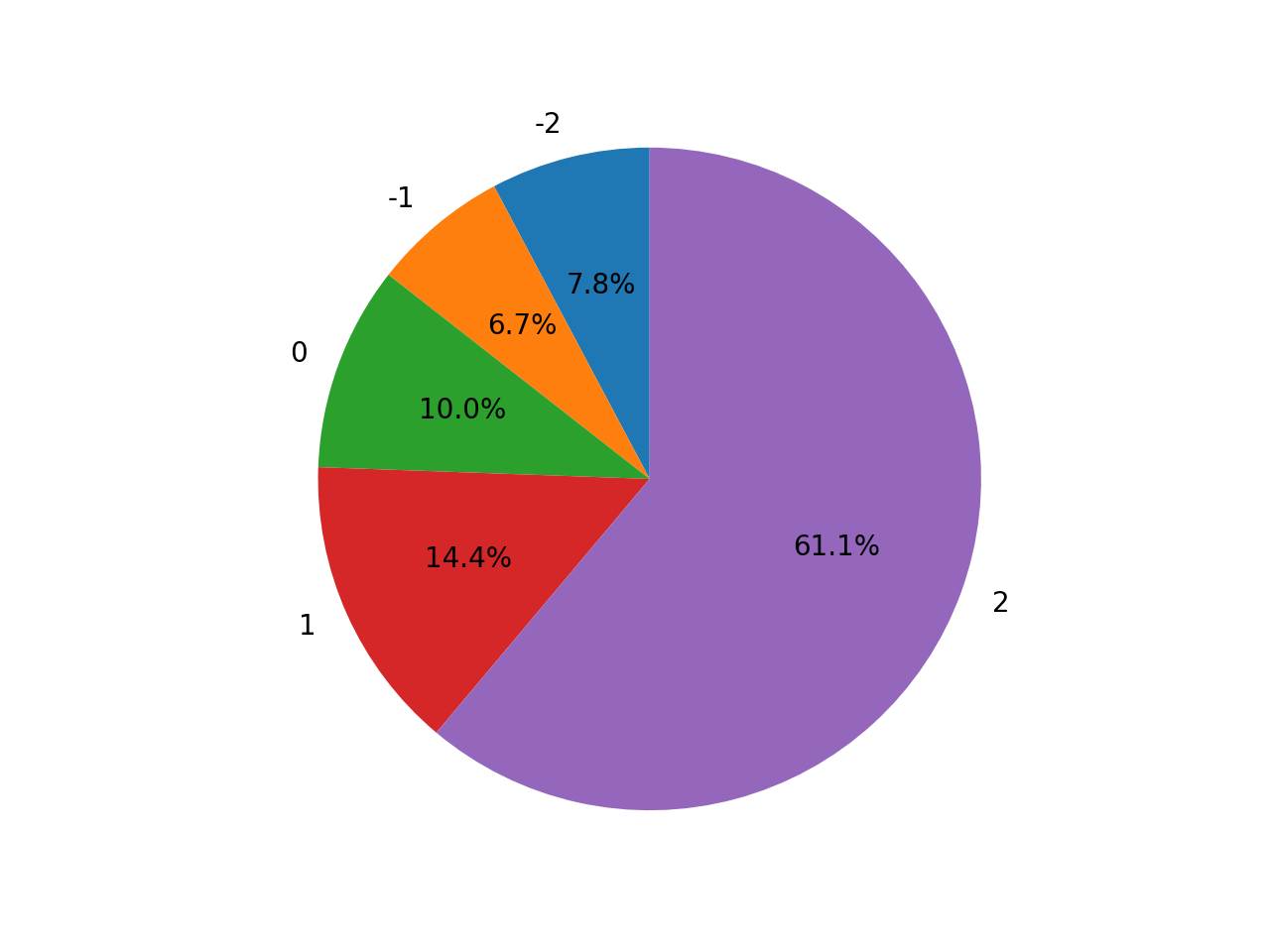}
    \caption{Distribution $D2$}
    \label{fig:d2}
  \end{subfigure}
  \\
  \begin{subfigure}[b]{0.48\textwidth}
    \centering
    \includegraphics[width=\textwidth]{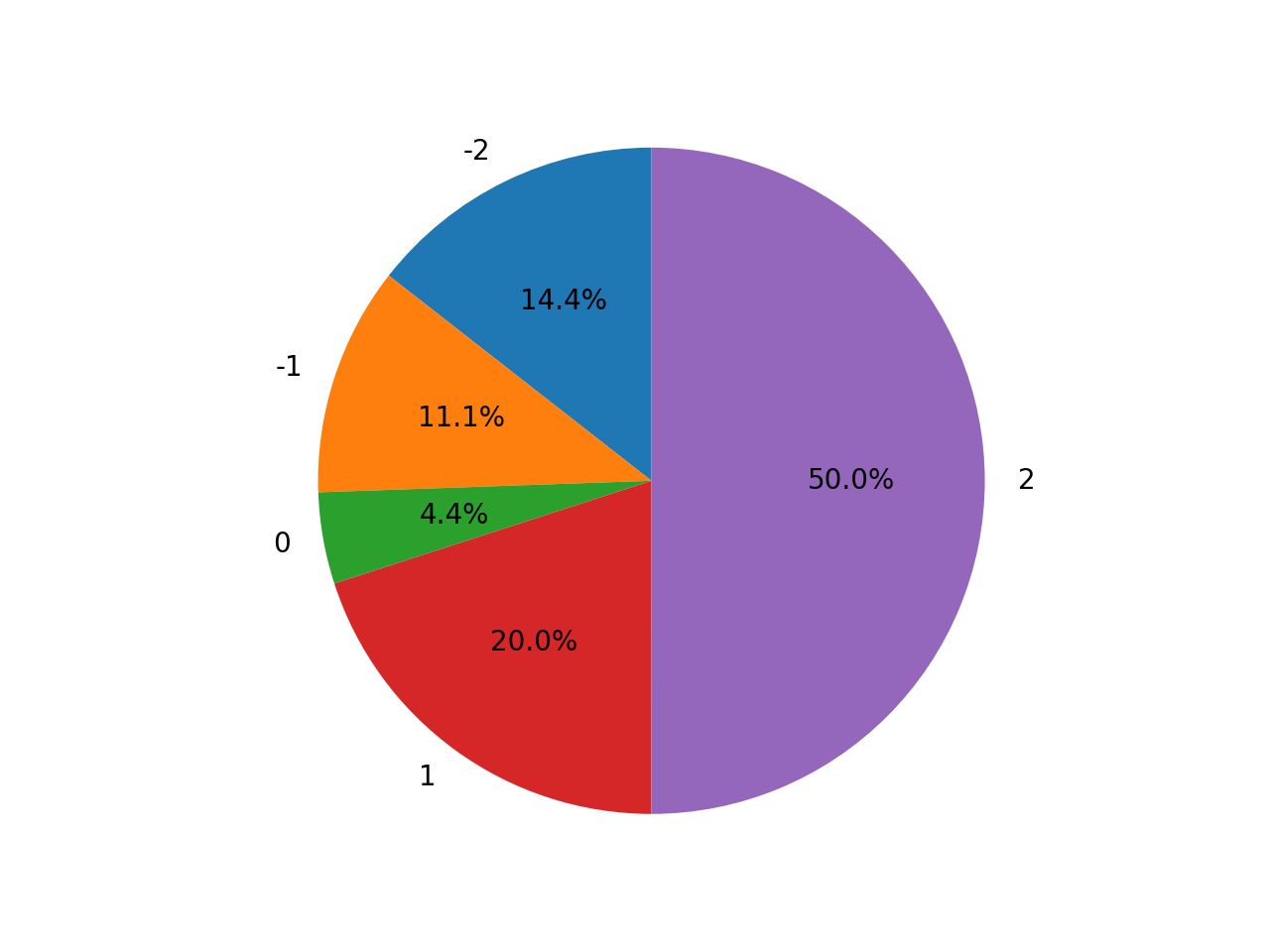}
    \caption{Distribution $D3$}
    \label{fig:d3}
  \end{subfigure}
  \hfill
  \begin{subfigure}[b]{0.48\textwidth}
    \centering
    \includegraphics[width=\textwidth]{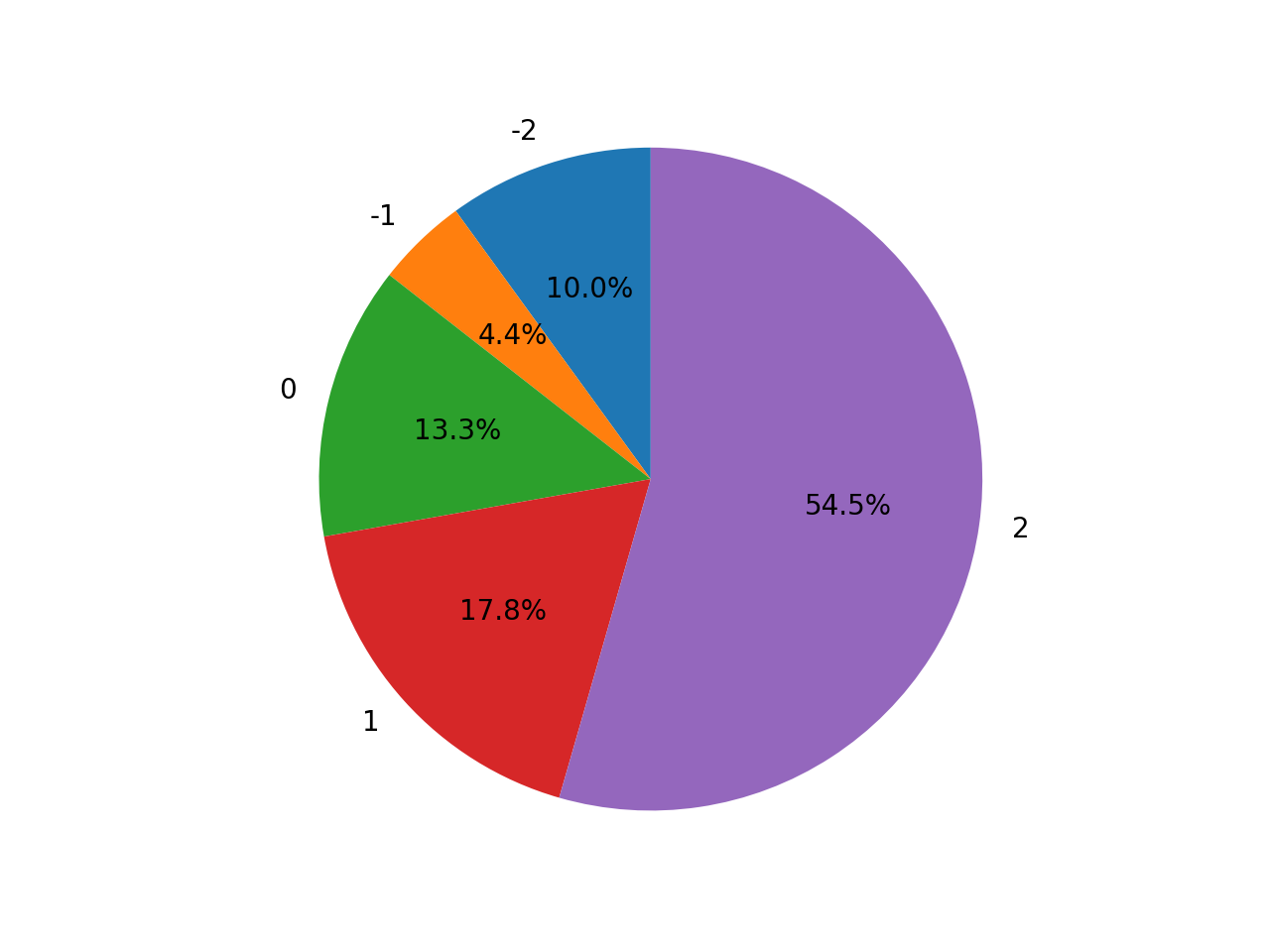}
    \caption{Distribution $D4$}
    \label{fig:d4}
  \end{subfigure}
  \caption{Drop-off fees adopted in the instances with $50$ vehicles and $600$ customers. Drop-off fees are expressed in Euro.}
  \label{fig:fares_distribuition50}
\end{figure}

In the case with $50$ vehicles (thus one vehicle every twelve potential customers), without dynamic pricing the CSO is able to satisfy approximately half of the total requests due to the reduced number of vehicles, see \Cref{tab:results:solutions50600}. Also in this case the number of requests is lower (approximately half) then the number of requests obtained by adjusting prices.
With dynamic pricing, the CSO is able to attract close to $1/3$ of the customers and to satisfy only approximately $25\%$ of them. Also in this case, the CSO benefits from the higher competition.
In more than $50\%$ of the origin-destination pairs the CSO is able to apply the highest drop-off fee ($2$ Euro) see \Cref{fig:fares_distribuition50} and to reposition the fleet in such a way to satisfy the requests with the highest revenue. Thus, a trend we observe in \Cref{tab:results:solutions50600} is that a pricing strategy allows the CSO to attract more demand but to satisfy only part of it. While this allows the CSO to exploit competition, many customers do not see their wish to use carsharing satisfied. This negative user experience might have an impact in the long run. This is however beyond the scope of this study.

Interestingly, in the case with $200$ vehicles (one every three customers -- see \Cref{tab:results:solutions200600}) the need for relocations is almost null, regardless of how prices are set.
The fleet is large enough to cover sufficiently well the entire business area and serve almost all requests.
On the other hand, with a fleet $50$ vehicles (one every twelve customers -- see \Cref{tab:results:solutions50600}) the need for relocations is more evident.
The fleet is now insufficient to cover the entire demand. In the case without dynamic pricing, fewer relocations are needed compared to the case with dynamic pricing.
This is due to the lower demand attracted (approximately $80$ to $107$ requests with a fleet of $50$ vehicles).
Many more relocations are performed when dynamic pricing is applied as a consequence of the higher demand generated (at least $167$ requests for $50$ vehicles). Thus the CSO finds it beneficial to move vehicles where they can generate more revenue.

A natural follow up question is the impact of relocations on profits. Therefore, we focused on those instances where more relocations were suggested, i.e., with $50$ vehicles and dynamic pricing, see \Cref{tab:results:solutions50600}. We solved the same instances, but this time preventing the model from making any relocations.
That is, vehicles were forced to remain in their initial positions.
The results indicate that the expected profit without relocations is only marginally lower.
Particularly, $98.5\%$ for distribution $D1$, $98.3\%$ for $D2$, $97.4\%$ for $D3$ and $98.8\%$ for $D4$.
This is due to the fact that relocations are expensive and can yield only a very minor increase in revenue.
That is, with $50$ vehicles and always more than $167$ requests, vehicles would always be rented, even if not relocated.
By relocating a vehicle the CSO is able to charge a higher drop-off fee, but bears the relocation cost. This results in a marginal profit increase. 

Relocations are however likely to generate a higher impact on profits when the vehicles-to-customers ratio is even smaller. 
We performed the same test with a vehicles-to-customers ratio equal to $1/100$ (i.e., $10$ vehicles $1000$ customers).
The results show that the profit without relocations was $70.07\%$ with distribution $D1$, $75.26\%$ with $D2$,  $76.22\%$ with distribution $D3$ and $78.07\%$ with distribution $D4$ (these percentages are calculated on the best upper bounds since a near optimal solution was available for all instances).
In all cases the percentage of vehicles relocated ranged between $10$ and $20\%$, similar to the case with $50$ vehicles and $600$ vehicles.
This means that the percentage of relocations remained approximately the same with a smaller vehicles-to-customers ratio, but had a much higher impact on profits.
Thus, it appears that a dynamic pricing strategy, coupled with a sufficiently large fleet (say more than one vehicle every twelve customers in our case),
decreases significantly the need for staff-based relocations. Otherwise, relocations remain an important tool even with a dynamic pricing strategy.

\section{Conclusions, limitations and future work}\label{sec:conclusions}
We presented a novel optimization model for jointly deciding carsharing prices and relocations.
The problem is modeled as a two-stage integer stochastic program in order to account for uncertainty in customers preferences.
An exact solution algorithm based on the integer L-Shaped method has been proposed.
Extensive tests have been performed on instances based on the municipality of Milan, in order to assess both the performance of the
solution algorithm and the type of solutions obtained.
The instances have been made available online for the sake of future research.

Results illustrate that, within times compatible with business practice, the method solves or finds a high quality solution to most instances. In addition, it finds a feasible solution to all instances considered. In contrast, CPLEX delivers a solution to only a few, small, instances.

The analysis of the solutions illustrates that a pricing strategy helps the CSO to significantly increase expected profits.
This is due to the increased demand generated and the resulting competition. In our instances the demand was approximately doubled compared to a situation without dynamic pricing.
This in turn generates higher expected profits by exploiting customer's higher willingness to pay.
The results also show that, by adopting a zone-based pricing strategy and employing a large enough fleet, the impact of staff-based relocations on profits becomes marginal.
On the other hand, the impact of relocations becomes more evident as the size of the fleet decreases.

A number of limitations remain to be addressed in future research, as we comment in what follows.
A pricing strategy which varies with each origin and destination, or frequently throughout the day, may not be applicable in all contexts, or raise concerns related to the potential complexity for users who would rather prefer a simpler pricing strategy.
The scope of this article was that of introducing a general model framework which could then be adapted to specific contexts and improved.
For example, the model proposed can be easily adapted to different time and space resolutions, i.e., it is possible to define the length of the target period and the discretization of the business area based on the specific needs.
In addition, the model may be easily modified to enforce that e.g., drop-off fees vary only according to the pick-up place or only according to the drop-off place. Future research may provide furter modifications and improvements.

The size of the instances used in this study is comparable with the size of the station-based carsharing in Milan which,
according to \cite{Mil20}, in 2018 counted $149$ shared vehicles and, on average $108$ daily rentals (see \Cref{tab:results:solutions200600} for a comparison).
Other examples are the station-based carsharing offered by \textit{Letsgo} (\url{https://letsgo.dk/}), which currently operates a fleet of around $200$ vehicles in Copenhagen, and 
\textit{Vy} (\url{https://www.vy.no/en/travelling-with-us/other-modes-of-transport/city-car}) that operates a fleet of 250 vehicles in Oslo.
Nevertheless, bigger fleets and a higher number of customers are likely to limit the practical efficiency of our exact method and call for faster, e.g., heuristic, methods. 

The performance of the algorithm with respect to a higher number of zones remains to be assessed.
Our instances, generated on the basis of \cite{HanP18}, contained ten zones. Supposedly, a more granular discretization of the business area is likely to have a negative impact on the practical applicability of the method. However, the benefits of a finer partition of the business area into pricing zones is to be addressed by further research,
particularly in the case of free-floating services. Effective discretization strategies and methods are, to our knowledge, still an open research question.

 As reported by \cite{ZoeK16}, the choice of carsharing users is also influenced by elements such as the type of vehicle and its proximity to the user.
In addition, comfort, weather conditions, and purpose of the trip are all factors which might influence customers decisions.
While proximity is considered in the form of walking time in the utility function we used in our experiments, the remaining elements are not captured explicitly, but are rather included in the portion of customer preferences that the CSO cannot explain. 
Future research might be set up to extend the model and utility function used in the tests in order to better capture customers behavior.

Several other sources of uncertainty affect the problem, that have not been considered in this study.
These include, e.g., the total number of customers appearing in each zone, and their destination.
Our model might account for this uncertainty by setting a sufficiently large number of customers.
Travel times, both with carsharing and with alternative transport services are also, to a certain degree, uncertain in practice. 
The impact of this uncertainty on solutions and profits remains to be understood.

The analysis of the solutions indicates that, by dynamically adjusting prices, the CSO is able to attract significantly more demand.
However, with a vehicles-to-customers ratio of $1$ to $12$ (see \Cref{tab:results:solutions50600}), the portion of the demand satisfied was, approximately, only $25\%$.
That is, the majority of the customers who would have used carsharing did not have the chance to do so.
As a consequence, users may perceive a low availability of the service. The effect of this in the long term remains to be clarified. 

Finally, our model is currently unable to use pricing as a preventive measure to encourage a profitable distribution of the fleet.
Consider two subsequent target periods, say $t_1$ and $t_2$, and two zones, say $A$ and $B$. Assume that the CSO expects high demand in zone $A$ in period $t_2$.
They may, consequently, set a lower price or an incentive in $t_1$, for renting cars in zone $B$ and delivering them in zone $A$, and/or disincentivize movements in the opposite direction.
In order to be able to optimize these decisions the proposed model should be extended to account for a multistage decision process.

\appendix
\crefalias{section}{appendix}  

\section{Sample Average Approximation}\label{sec:saa}
Let $\xi_{1},\ldots,\xi_{S}$ be an $S$-dimensional iid sample of $\tilde{\xi}$ and let $\mtc{S}=\{1,\ldots,S\}$.
Let decision variable $y_{vrls}$ be equal to $1$ if request $r\in\mtc{R}(\xi_s)$ is satisfied by vehicle $v$ at level $l$ under realization $s$, $0$ otherwise.
The SAA of problem \eqref{eq:1S} can be stated as follows.
\begin{subequations}
  \label{eq:SAA}
  \begin{align}
    \label{eq:SAA:obj}&\max-\sum_{v\in\mtc{V}}\sum_{i\in\mtc{I}}C^R_{vi}z_{vi}+\frac{1}{|\mtc{S}|}\sum_{s\in\mtc{S}}\sum_{r\in\mtc{R}(\xi_s)}\sum_{v\in\mtc{V}}\sum_{l\in\mtc{L}_{r}(\xi_s)}R_{rl}y_{vrls}\\
    \label{eq:SAA:c3}&\sum_{i\in\mtc{I}}z_{vi} = 1   & v\in\mtc{V}\\
    \label{eq:SAA:c6}&\sum_{l\in\mtc{L}}\lambda_{ijl}=1& i\in\mtc{I},j\in\mtc{J}\\
    \label{eq:SAA:c1}&\sum_{v\in\mtc{V}}\sum_{l\in\mtc{L}_r(\xi_s)}y_{vrls}\leq 1   & r\in\mtc{R}(\xi_s), s\in\mtc{S}\\
    \label{eq:SAA:c2}&\sum_{r\in\mtc{R}(\xi_s)}\sum_{l\in\mtc{L}_{r}(\xi_s)}y_{vrls}\leq 1   & v\in\mtc{V},s\in\mtc{S}\\
    \label{eq:SAA:c4}&\sum_{l\in\mtc{L}_{r_1}(\xi_s)}y_{v,r_1,l,s} - z_{v,i(r_1)} + \sum_{r_2\in\mtc{R}_{r_1}(\xi_s)}\sum_{l\in\mtc{L}_{r_2}(\xi_s)}y_{v,r_2,l,s}\leq 0   & r_1\in\mtc{R}(\xi_s),v\in\mtc{V}.s\in\mtc{S}\\[5pt]
    \nonumber&y_{v,r_1,l_1,s}+ \sum_{r_2\in\mtc{R}_{r_1}(\xi_s)}\sum_{l_2\in\mtc{L}_{r_2}(\xi_s)}y_{v,r_2,l_2,s}+ \sum_{v_1\in\mtc{V}:v_1\neq v}y_{v_1,r_1,l_1,s}  &\\
    \label{eq:SAA:c5} &\geq \lambda_{i(r_1),j(r_j),l_1} + z_{v,i(r_1)}-1   & r_1\in\mtc{R}(\xi_s),v\in\mtc{V}, l_1\in\mtc{L}_{r_1}(\xi_s),s\in\mtc{S}\\
    \label{eq:SAA:c7}&\sum_{v\in\mtc{V}}y_{vrls}\leq\lambda_{i(r),j(r),l}& r\in\mtc{R}(\xi_s),l\in\mtc{L}_{r}(\xi_s),s\in\mtc{S}\\
                     &y_{vrls}\in\{0,1\}              & r\in\mtc{R}(\xi_s),v\in\mtc{V},l\in\mtc{L}_{r}(\xi_s),s\in\mtc{S}\\
                      &z_{vi}\in\{0,1\}              & i\in\mtc{I},v\in\mtc{V}\\
                      &\lambda_{ijl}\in\{0,1\} & i\in\mtc{I},j\in\mtc{I},l\in\mtc{L}.
  \end{align}
\end{subequations}

\section{Size of the instances}\label{sec:app:size}
This section reports the size of the instances for the base case in \Cref{tab:app:sizeIC} and for the case with customers profiled individually in \Cref{tab:app:sizeDC}.
\begin{longtable}{rrrr|rr}
  \caption{Size of the SAA model without decomposition for all instances tested in the base case.}\label{tab:app:sizeIC}\\
  \toprule
  $|\mtc{V}|$ & $|\mtc{K}|$  &  $\alpha^{FROM}$ &  $\alpha^{TO}$  & \# Variables &    \# Constraints \\
  \midrule
  50 &  200 &       0.2 &     0.2 &   95350 &  127528 \\
  50 &  400 &       0.2 &     0.2 &  199800 &  269257 \\
  50 &  600 &       0.2 &     0.2 &  296000 &  399868 \\
 100 &  200 &       0.2 &     0.2 &  190200 &  252478 \\
 100 &  400 &       0.2 &     0.2 &  399500 &  533561 \\
 100 &  600 &       0.2 &     0.2 &  591600 &  791919 \\
 200 &  400 &       0.2 &     0.2 &  798300 & 1061560 \\
 200 &  600 &       0.2 &     0.2 & 1182900 & 1576120 \\
\midrule
  50 &  200 &       0.2 &     0.8 &   89450 &  119419 \\
  50 &  400 &       0.2 &     0.8 &  189750 &  254569 \\
  50 &  600 &       0.2 &     0.8 &  279350 &  375643 \\
 100 &  200 &       0.2 &     0.8 &  178500 &  236520 \\
 100 &  400 &       0.2 &     0.8 &  379200 &  504170 \\
 100 &  600 &       0.2 &     0.8 &  558300 &  743843 \\
 200 &  400 &       0.2 &     0.8 &  758100 & 1003471 \\
 200 &  600 &       0.2 &     0.8 & 1116500 & 1480645 \\
\midrule
  50 &  200 &       0.8 &     0.2 &   94300 &  125998 \\
  50 &  400 &       0.8 &     0.2 &  195900 &  263086 \\
  50 &  600 &       0.8 &     0.2 &  287050 &  387577 \\
 100 &  200 &       0.8 &     0.2 &  188100 &  249448 \\
 100 &  400 &       0.8 &     0.2 &  391800 &  521542 \\
 100 &  600 &       0.8 &     0.2 &  573900 &  767780 \\
 200 &  400 &       0.8 &     0.2 &  783300 & 1038043 \\
 200 &  600 &       0.8 &     0.2 & 1147300 & 1527880 \\
\midrule
  50 &  200 &       0.8 &     0.8 &   90300 &  120337 \\
  50 &  400 &       0.8 &     0.8 &  185750 &  249469 \\
  50 &  600 &       0.8 &     0.8 &  279500 &  375949 \\
 100 &  200 &       0.8 &     0.8 &  180100 &  238237 \\
 100 &  400 &       0.8 &     0.8 &  371300 &  494272 \\
 100 &  600 &       0.8 &     0.8 &  558800 &  744752 \\
 200 &  400 &       0.8 &     0.8 &  742300 &  983773 \\
  200 &  600 &       0.8 &     0.8 & 1117700 & 1482856 \\
  \midrule
    &      &           &         &  452819 &  603174 \\
  \bottomrule
\end{longtable}

\begin{longtable}{rrrr|rr}
  \caption{Size of the SAA model without decomposition for all instances tested with individual customers profiles.}\label{tab:app:sizeDC}\\
  \toprule
  $|\mtc{V}|$ & $|\mtc{K}|$  &  $\alpha^{FROM}$ &  $\alpha^{TO}$  & \# Variables &    \# Constraints \\
  \midrule
  50 &  200 &       0.2 &     0.2 &   84500 &  116461 \\
  50 &  400 &       0.2 &     0.2 &  171000 &  237127 \\
  50 &  600 &       0.2 &     0.2 &  250800 &  347389 \\
 100 &  200 &       0.2 &     0.2 &  168600 &  230662 \\
 100 &  400 &       0.2 &     0.2 &  341600 &  469628 \\
 100 &  600 &       0.2 &     0.2 &  501200 &  687990 \\
 200 &  400 &       0.2 &     0.2 &  682900 &  934729 \\
  200 &  600 &       0.2 &     0.2 & 1002100 & 1369291 \\
  \midrule
  50 &  200 &       0.2 &     0.8 &   77200 &  107281 \\
  50 &  400 &       0.2 &     0.8 &  163300 &  226723 \\
  50 &  600 &       0.2 &     0.8 &  239650 &  333007 \\
 100 &  200 &       0.2 &     0.8 &  153900 &  212381 \\
 100 &  400 &       0.2 &     0.8 &  326100 &  448923 \\
 100 &  600 &       0.2 &     0.8 &  478900 &  659508 \\
 200 &  400 &       0.2 &     0.8 &  651700 &  893323 \\
  200 &  600 &       0.2 &     0.8 &  957500 & 1312609 \\
  \midrule
  50 &  200 &       0.8 &     0.2 &   83100 &  114676 \\
  50 &  400 &       0.8 &     0.2 &  168100 &  232792 \\
  50 &  600 &       0.8 &     0.2 &  243900 &  338821 \\
 100 &  200 &       0.8 &     0.2 &  165900 &  227329 \\
 100 &  400 &       0.8 &     0.2 &  335700 &  460942 \\
 100 &  600 &       0.8 &     0.2 &  487300 &  670921 \\
 200 &  400 &       0.8 &     0.2 &  670900 &  917242 \\
  200 &  600 &       0.8 &     0.2 &  974100 & 1335121 \\
  \midrule
  50 &  200 &       0.8 &     0.8 &   79300 &  110035 \\
  50 &  400 &       0.8 &     0.8 &  161150 &  224428 \\
  50 &  600 &       0.8 &     0.8 &  243800 &  340351 \\
 100 &  200 &       0.8 &     0.8 &  158100 &  217835 \\
 100 &  400 &       0.8 &     0.8 &  321800 &  444378 \\
 100 &  600 &       0.8 &     0.8 &  487300 &  674153 \\
 200 &  400 &       0.8 &     0.8 &  643300 &  884479 \\
  200 &  600 &       0.8 &     0.8 &  974100 & 1341553 \\
  \midrule
 &      &           &         &  389025 &  535065 \\  
\bottomrule
\end{longtable}

\section{Effect of valid inequality}\label{sec:app:ICplusVI}

\begin{longtable}{rrrr|rrrr}
  \caption{Results of the L-Shaped method with the addition of \cref{eq:VI} on the instances with $\alpha^V=0.2$.}\label{tab:resultsICVI02}\\
  \toprule
  $|\mtc{V}|$ & $|\mtc{K}|$  &  $\alpha^{FROM}$ &  $\alpha^{TO}$  & \texttt{gap} &    \texttt{gapR} &   \texttt{gap50} & \texttt{t}  \\
  \midrule
    50 &  200 &       0.2 &     0.2 &  0.6967 &  18.6752 &  0.6967 & 1800.20 \\
  50 &  400 &       0.2 &     0.2 & 17.2850 &  66.5573 & 20.6265 & 1800.52 \\
  50 &  600 &       0.2 &     0.2 & 30.5015 &  87.2440 & 39.2311 & 1809.94 \\
 100 &  200 &       0.2 &     0.2 &  0.0975 &   9.5524 &  0.0975 & 1800.98 \\
 100 &  400 &       0.2 &     0.2 &  1.9480 &  21.2273 &  3.9377 & 1805.67 \\
 100 &  600 &       0.2 &     0.2 & 11.6208 &  35.8381 & 14.7910 & 1801.95 \\
 200 &  400 &       0.2 &     0.2 &  0.0000 &   6.3509 &       - &  139.41 \\
 200 &  600 &       0.2 &     0.2 &  0.0777 &  10.4516 &  0.3662 & 1800.78 \\
\midrule
  50 &  200 &       0.2 &     0.8 &  0.3381 &  17.6529 &  0.3381 & 1800.63 \\
  50 &  400 &       0.2 &     0.8 &  9.6397 &  62.4842 & 11.6808 & 1801.21 \\
  50 &  600 &       0.2 &     0.8 & 19.5455 &  90.2582 & 20.6217 & 1800.03 \\
 100 &  200 &       0.2 &     0.8 &  0.2815 &   7.2972 &  0.2815 & 1800.13 \\
 100 &  400 &       0.2 &     0.8 &  0.7148 &  19.6630 &  1.0807 & 1801.47 \\
 100 &  600 &       0.2 &     0.8 &  9.8896 &  34.0425 & 11.4341 & 1812.93 \\
 200 &  400 &       0.2 &     0.8 &  0.0000 &  41.7344 &       - &  115.74 \\
 200 &  600 &       0.2 &     0.8 &  0.4710 &  12.6665 &  0.6677 & 1814.99 \\
\midrule
  50 &  200 &       0.8 &     0.2 &  1.9850 &  23.9553 &  2.2390 & 1800.13 \\
  50 &  400 &       0.8 &     0.2 & 21.8829 &  85.1665 & 25.2703 & 1800.02 \\
  50 &  600 &       0.8 &     0.2 & 44.9080 & 108.3757 & 62.3967 & 1803.83 \\
 100 &  200 &       0.8 &     0.2 &  0.5238 &  15.4557 &  0.5320 & 1801.09 \\
 100 &  400 &       0.8 &     0.2  &  6.1743 &  22.7406 &  9.0017 & 1806.48 \\
 100 &  600 &       0.8 &     0.2 & 23.6427 &  42.5987 & 23.6588 & 1816.20 \\
 200 &  400 &       0.8 &     0.2 &  0.3072 &  12.4280 &  0.3209 & 1804.22 \\
  200 &  600 &       0.8 &     0.2&  6.7474 &  15.3980 & 17.5717 & 1865.22 \\
\midrule 
  50 &  200 &       0.8 &     0.8 &  0.8509 &  22.5871 &  0.9463 & 1800.15 \\
  50 &  400 &       0.8 &     0.8 & 15.4867 &  55.7576 & 19.9941 & 1801.46 \\
  50 &  600 &       0.8 &     0.8 & 40.1722 & 106.9156 & 41.2931 & 1807.37 \\
 100 &  200 &       0.8 &     0.8 &  0.7932 &  13.3519 &  0.8354 & 1802.32 \\
 100 &  400 &       0.8 &     0.8 &  4.4619 &  17.1371 &  7.3475 & 1803.98 \\
 100 &  600 &       0.8 &     0.8 & 19.1824 &  46.7222 & 23.6637 & 1823.21 \\
 200 &  400 &       0.8 &     0.8 &  0.0963 &   7.2607 &  0.0963 & 1800.78 \\
  200 &  600 &       0.8 &     0.8  &  5.6194 &  16.7484 &  6.8357 & 1823.23 \\
  \midrule
  &    &         &        &  9.2482 &  36.0717 & 12.2618 & 1702.07 \\
  \bottomrule
\end{longtable}

\begin{longtable}{rrrr|rrrr}
  \caption{Results of the L-Shaped method with the addition of \cref{eq:VI} on the instances with $\alpha^V=0.8$.}\label{tab:resultsICVI08}\\
  \toprule
  $|\mtc{V}|$ & $|\mtc{K}|$  &  $\alpha^{FROM}$ &  $\alpha^{TO}$  & \texttt{gap} &    \texttt{gapR} &   \texttt{gap50} & \texttt{t}  \\
  \midrule
  50 &  200 &       0.2 &     0.2  &  0.4613 & 37.8092 &  0.5067 & 1800.37 \\
  50 &  400 &       0.2 &     0.2  & 11.6098 & 42.1383 & 13.2894 & 1802.87 \\
  50 &  600 &       0.2 &     0.2  & 30.3469 & 77.8433 & 39.0602 & 1807.97 \\
 100 &  200 &       0.2 &     0.2  &  0.1113 &  6.1394 &  0.2584 & 1800.27 \\
 100 &  400 &       0.2 &     0.2  &  4.0199 & 22.6182 &  6.3574 & 1800.15 \\
 100 &  600 &       0.2 &     0.2  & 14.6664 & 35.6069 & 15.9543 & 1805.19 \\
 200 &  400 &       0.2 &     0.2  &  0.9991 &  5.3451 &  1.0499 & 1810.37 \\
 200 &  600 &       0.2 &     0.2  &  4.5424 & 10.6059 & 14.3498 & 1889.13 \\
\midrule
  50 &  200 &       0.2 &     0.8  &  0.8298 & 33.4059 &  1.0115 & 1800.11 \\
  50 &  400 &       0.2 &     0.8  & 11.8754 & 42.0692 & 16.4411 & 1801.42 \\
  50 &  600 &       0.2 &     0.8  & 24.6840 & 96.0829 & 32.7535 & 1800.96 \\
 100 &  200 &       0.2 &     0.8  &  0.0291 &  4.7786 &  0.0291 & 1800.47 \\
 100 &  400 &       0.2 &     0.8  &  2.1136 & 20.3964 &  3.4560 & 1805.67 \\
 100 &  600 &       0.2 &     0.8  & 12.3031 & 35.6683 & 14.0160 & 1817.99 \\
 200 &  400 &       0.2 &     0.8  &  0.2123 &  9.3896 &  0.4837 & 1800.23 \\
 200 &  600 &       0.2 &     0.8 &  3.7770 & 17.6274 &  4.4652 & 1833.46 \\
\midrule
  50 &  200 &       0.8 &     0.2 &  0.0389 & 12.5954 &  0.0389 & 1800.08 \\
  50 &  400 &       0.8 &     0.2  & 13.3947 & 65.9956 & 15.1977 & 1800.04 \\
  50 &  600 &       0.8 &     0.2  & 34.8704 & 82.1348 & 38.2156 & 1803.88 \\
 100 &  200 &       0.8 &     0.2 &  0.0000 &  9.4214 &       - &   26.88 \\
 100 &  400 &       0.8 &     0.2  &  0.2335 &  4.0455 &  0.2348 & 1801.63 \\
 100 &  600 &       0.8 &     0.2  &  9.1181 & 27.6410 & 11.0863 & 1801.40 \\
 200 &  400 &       0.8 &     0.2 &  0.0020 &  8.6414 &       - &  205.30 \\
 200 &  600 &       0.8 &     0.2  &  0.0034 & 11.5793 &       - &  786.63 \\
\midrule
  50 &  200 &       0.8 &     0.8  &  0.0083 & 14.9683 &       - &   24.39 \\
  50 &  400 &       0.8 &     0.8  &  9.8202 & 44.5269 & 12.5750 & 1801.95 \\
  50 &  600 &       0.8 &     0.8  & 30.9029 & 99.1525 & 36.1072 & 1808.71 \\
 100 &  200 &       0.8 &     0.8 &  0.0000 &  0.4394 &       - &   21.84 \\
 100 &  400 &       0.8 &     0.8  &  0.1667 & 17.2446 &  0.2030 & 1801.88 \\
 100 &  600 &       0.8 &     0.8 &  6.1475 & 13.2578 &  6.5530 & 1805.95 \\
 200 &  400 &       0.8 &     0.8 &  0.0000 & 42.1315 &       - &  122.34 \\
  200 &  600 &       0.8 &     0.8 &  0.0294 &  8.4615 &  0.0306 & 1802.96 \\
  \midrule
   &    &        &      &  7.1037 & 29.9925 & 10.9125 & 1506.01 \\
\bottomrule  
\end{longtable}

\section{Results on the instances with individual customer profiles after 5 hours}\label{app:resultsLongDC}
\begin{longtable}{rrrr|rrrr}
  \caption{Results of the L-Shaped method with the addition of \cref{eq:VI} on the instances with $\alpha^V=0.2$ and individual customer profiles with a time limit of $18 000$ seconds and a $1\%$ target optimality gap.}\label{tab:resultsLongDCVI02}\\
  \toprule
  $|\mtc{V}|$ & $|\mtc{K}|$  &  $\alpha^{FROM}$ &  $\alpha^{TO}$  & \texttt{gap} &    \texttt{gapR} &   \texttt{gap50} & \texttt{t}  \\
  \midrule
    50 &  200 &       0.2 &     0.2 &  0.9337 & 21.5553 &       - &   120.82 \\
  50 &  400 &       0.2 &     0.2 &  3.0295 & 59.3308 &  4.1428 & 18002.76 \\
  50 &  600 &       0.2 &     0.2 & 10.1187 & 52.5981 & 11.3538 & 18000.13 \\
 100 &  200 &       0.2 &     0.2 &  0.1841 &  8.9488 &       - &    16.69 \\
 100 &  400 &       0.2 &     0.2 &  0.8770 & 26.1941 &       - &   841.63 \\
 100 &  600 &       0.2 &     0.2 &  2.6675 & 40.2920 &  3.4833 & 18000.07 \\
 200 &  400 &       0.2 &     0.2 &  0.1093 &  9.9094 &       - &   224.15 \\
 200 &  600 &       0.2 &     0.2 &  0.3517 & 14.7626 &       - &   479.14 \\
\midrule
  50 &  200 &       0.2 &     0.8 &  0.8707 & 20.9177 &       - &    11.71 \\
  50 &  400 &       0.2 &     0.8 &  1.8429 & 43.2969 &  2.5193 & 18000.92 \\
  50 &  600 &       0.2 &     0.8 &  9.9981 & 63.4863 & 10.4941 & 18000.11 \\
 100 &  200 &       0.2 &     0.8 &  0.2865 & 10.6250 &       - &    19.35 \\
 100 &  400 &       0.2 &     0.8 &  0.9395 & 22.5452 &       - &   535.48 \\
 100 &  600 &       0.2 &     0.8 &  2.8913 & 36.2728 &  4.4123 & 18008.31 \\
 200 &  400 &       0.2 &     0.8 &  0.0000 & 59.2478 &       - &   181.93 \\
 200 &  600 &       0.2 &     0.8 &  0.5768 & 15.6153 &       - &   726.62 \\
\midrule
  50 &  200 &       0.8 &     0.2 &  0.7559 & 22.1339 &       - &    44.38 \\
  50 &  400 &       0.8 &     0.2 &  1.5553 & 58.5075 &  1.7269 & 18000.08 \\
  50 &  600 &       0.8 &     0.2 & 14.8996 & 72.8427 & 17.8340 & 18000.08 \\
 100 &  200 &       0.8 &     0.2 &  0.8890 & 13.9807 &       - &   121.73 \\
 100 &  400 &       0.8 &     0.2 &  0.9990 & 11.8809 &       - &  1379.79 \\
 100 &  600 &       0.8 &     0.2 &  4.8890 & 35.7724 &  5.7060 & 18006.26 \\
 200 &  400 &       0.8 &     0.2 &  0.1760 & 11.8485 &       - &   149.53 \\
 200 &  600 &       0.8 &     0.2 &  0.9180 & 23.4393 &       - &  5103.61 \\
\midrule
  50 &  200 &       0.8 &     0.8 &  0.9921 & 21.9851 &       - &    82.89 \\
  50 &  400 &       0.8 &     0.8 &  1.7121 & 53.8589 &  2.1396 & 18002.06 \\
  50 &  600 &       0.8 &     0.8 &  9.9533 & 70.6343 & 10.2331 & 18000.25 \\
 100 &  200 &       0.8 &     0.8 &  0.9995 & 14.5185 &       - &   101.45 \\
 100 &  400 &       0.8 &     0.8 &  0.9953 & 30.1037 &       - &  4394.49 \\
 100 &  600 &       0.8 &     0.8 &  3.0919 & 34.8283 &  5.3525 & 18000.65 \\
 200 &  400 &       0.8 &     0.8 &  0.0001 & 11.5481 &       - &   164.89 \\
  200 &  600 &       0.8 &     0.8 &  0.9972 & 15.4424 &       - &  4699.97 \\
  \midrule
     &      &           &         &  2.4844 & 31.5289 &  6.6165 &  7356.94 \\
  \bottomrule
\end{longtable}

\begin{longtable}{rrrr|rrrr}
  \caption{Results of the L-Shaped method with the addition of \cref{eq:VI} on the instances with $\alpha^V=0.8$ and individual customer profiles with a time limit of $18 000$ seconds and a $1\%$ target optimality gap.}\label{tab:resultsLongDCVI08}\\
  \toprule
  $|\mtc{V}|$ & $|\mtc{K}|$  &  $\alpha^{FROM}$ &  $\alpha^{TO}$  & \texttt{gap} &    \texttt{gapR} &   \texttt{gap50} & \texttt{t}  \\
  \midrule
    50 &  200 &       0.2 &     0.2 &  1.5133 &  17.7774 &  1.7378 & 18000.07 \\
  50 &  400 &       0.2 &     0.2 &  8.8698 &  60.3835 & 11.4205 & 18000.10 \\
  50 &  600 &       0.2 &     0.2 & 16.0592 & 105.0348 & 19.0821 & 18000.14 \\
 100 &  200 &       0.2 &     0.2 &  0.9991 &   7.5390 &       - &  2345.67 \\
 100 &  400 &       0.2 &     0.2 &  1.9678 &  33.3830 &  2.2880 & 18000.46 \\
 100 &  600 &       0.2 &     0.2 &  6.7921 &  29.7073 &  7.3807 & 18000.10 \\
 200 &  400 &       0.2 &     0.2 &  0.9852 &  22.1296 &       - &  3051.54 \\
 200 &  600 &       0.2 &     0.2 &  1.6326 &  17.3350 &  3.1395 & 18000.54 \\
\midrule
  50 &  200 &       0.2 &     0.8 &  1.1847 &  20.1198 &  1.3004 & 18000.58 \\
  50 &  400 &       0.2 &     0.8 &  8.7233 &  53.4891 &  9.5341 & 18000.09 \\
  50 &  600 &       0.2 &     0.8 & 14.6785 &  89.9520 & 16.2307 & 18000.09 \\
 100 &  200 &       0.2 &     0.8 &  0.5693 &   6.9976 &       - &    43.73 \\
 100 &  400 &       0.2 &     0.8 &  1.7097 &  32.4136 &  1.9677 & 18000.31 \\
 100 &  600 &       0.2 &     0.8 &  6.3776 &  32.5442 &  7.1209 & 18000.09 \\
 200 &  400 &       0.2 &     0.8 &  0.9601 &  28.3808 &       - &  4120.56 \\
 200 &  600 &       0.2 &     0.8 &  4.6562 &  26.4194 &  5.5152 & 18000.11 \\
\midrule
  50 &  200 &       0.8 &     0.2 &  0.7444 &  17.0926 &       - &    15.28 \\
  50 &  400 &       0.8 &     0.2 &  1.5675 &  32.0731 &  1.7184 & 18000.13 \\
  50 &  600 &       0.8 &     0.2 &  8.4393 &  75.7595 &  9.0644 & 18000.30 \\
 100 &  200 &       0.8 &     0.2 &  0.0000 &  14.6898 &       - &    25.29 \\
 100 &  400 &       0.8 &     0.2 &  0.9903 &  18.9859 &       - &   579.88 \\
 100 &  600 &       0.8 &     0.2 &  3.8100 &  40.7039 &  5.1358 & 18000.36 \\
 200 &  400 &       0.8 &     0.2 &  0.0000 &  61.0074 &       - &   164.47 \\
 200 &  600 &       0.8 &     0.2 &  0.3060 &  12.2654 &       - &   441.39 \\
\midrule
  50 &  200 &       0.8 &     0.8 &  0.9938 &   4.6541 &       - &    14.80 \\
  50 &  400 &       0.8 &     0.8 &  0.9997 &  42.8840 &  1.0556 & 12864.17 \\
  50 &  600 &       0.8 &     0.8 &  6.3284 &  78.3238 &  6.5753 & 18000.09 \\
 100 &  200 &       0.8 &     0.8 &  0.2054 &   1.4263 &       - &    23.31 \\
 100 &  400 &       0.8 &     0.8 &  0.9188 &   6.4020 &       - &   104.09 \\
 100 &  600 &       0.8 &     0.8 &  2.6157 &  36.8295 &  3.3959 & 18000.60 \\
 200 &  400 &       0.8 &     0.8 &  0.0000 &  61.6141 &       - &   127.80 \\
  200 &  600 &       0.8 &     0.8 &  0.3572 &  13.3682 &       - &   445.02 \\
  \midrule
     &      &           &         &  3.3111 &  34.4277 &  6.3146 & 10324.10 \\
  \bottomrule
\end{longtable}

\section{Example Heuristic}\label{app:ILS}

In this section we present a simple Iterated Local Search (ILS) to find primal solutions to the problem. In a nutshell, an ILS works as follows: given an initial (current at a generic iteration) solution it performs a local search. To escape local optima, the solution returned by the local search is randomly perturbed and the local search restarted. This procedure is repeated until a stopping criteria is met. In what follows we explain how this procedure is adapted to our problem. 

We encode solutions $\zeta$ to the original problem using
\begin{itemize}
\item a vector $\Pi(\zeta)\in\mathbb{N}^{|\mtc{V}|}$ defining the position of the vehicles. The $v$-th element of the vector is an integer $i\in\mtc{I}$ identifying the zone where vehicle $v$ becomes available (possibly after relocation).
\item a matrix $\Lambda(\zeta)\in \mathbb{N}^{|\mtc{I}|\times|\mtc{I}|}$ defining the drop-off-fees. The $i-j$-th element of the matrix is an integer $l\in\mtc{L}$ identifying the drop-off-bee between zone $i$ and $j$.
\end{itemize}
Let $H(a,b)$ be the Hamming distance between two vectors or, in the case of matrices, concatenations of rows. We define two types of neighborhoods.
\begin{itemize}
\item $\mtc{N}^\Pi(\zeta)=\big\{\zeta'|H(\Pi(\zeta),\Pi(\zeta'))=1,\Lambda(\zeta)=\Lambda(\zeta')\big\}$. In words, it defines all solutions which can be obtained from $\zeta$ by changing solely one vehicle position.
\item $\mtc{N}^\Lambda(\zeta)=\big\{\zeta'|\Pi(\zeta)=\Pi(\zeta'), H(\Lambda(\zeta)),\Lambda(\zeta'))=1\big\}$. In words, it defines all solutions which can be obtained by changing solely one drop-off fee.
\end{itemize}
We define two \textit{first improvement operators}:
\begin{itemize}
\item $f^\Pi(\zeta):\mathbb{N}^{|\mtc{V}|}\times \mathbb{N}^{|\mtc{I}|\times|\mtc{I}|} \to \mtc{N}^\Pi(\zeta)$ scans the neighborhood $\mtc{N}^\Pi(\zeta)$ and returns the first improving solution (i.e., with a higher fitness value) if it exists, $\zeta$ otherwise.
\item $f^\Lambda(\zeta):\mathbb{N}^{|\mtc{V}|}\times \mathbb{N}^{|\mtc{I}|\times|\mtc{I}|} \to \mtc{N}^\Lambda(\zeta)$ scans the neighborhood $\mtc{N}^\Lambda(\zeta)$ and returns the first improving solution (i.e., with a higher fitness value) if it exists, $\zeta$ otherwise.
\end{itemize}
Let $z(\zeta)$ be the fitness function defining how each solution is evaluated. It is defined as the objective function of the original problem \eqref{eq:1S:obj}.
For each solution $\zeta$ considered, the second-stage revenue is computed as illustrated in \Cref{sec:ls:subproblems}.
Given a solution $\zeta$, we define two types of local search
\begin{itemize}
\item $LS^\Pi(\zeta):\mathbb{N}^{|\mtc{V}|}\times \mathbb{N}^{|\mtc{I}|\times|\mtc{I}|}\to \mtc{N}^\Pi(\zeta)$ which performs a local search on the $\mtc{N}^\Pi(\zeta)$ neighborhood using the $f^\Pi(\zeta)$ first improvement operator
\item $LS^\Lambda(\zeta): \mathbb{N}^{|\mtc{V}|}\times \mathbb{N}^{|\mtc{I}|\times|\mtc{I}|}\to \mtc{N}^\Lambda(\zeta)$ which performs a local search on the $\mtc{N}^\Lambda(\zeta)$ neighborhood using the $f^\Lambda(\zeta)$ first improvement operator.
\end{itemize}
\Cref{alg:LS} sketches the local search procedures.

\begin{algorithm}[h]
  \caption{Local Search}
  \label{alg:LS}
  \begin{algorithmic}[1]
    \STATE INPUT: $\zeta$,  Operator $f(\zeta)$ to use (i.e., $f^\Pi(\zeta)$ or $f^\Lambda(\zeta)$).
    \STATE INPUT: \texttt{TIMELIMIT}
    \STATE \texttt{FOUND}=\texttt{TRUE}
    \STATE $\zeta^{CURRENT}=\zeta$
    \WHILE{\textit{FOUND} AND \texttt{ELAPSEDTIME}$\leq$\texttt{TIMELIMIT}}
    \STATE $\zeta^N\gets f(\zeta^{CURRENT})$
    \IF{$z(\zeta^N)>z(\zeta^{CURRENT})$}
    \STATE $\zeta^{CURRENT}\gets \zeta^N$
    \ELSE
    \STATE \texttt{FOUND}=\texttt{FALSE}
    \ENDIF
    \ENDWHILE
    \RETURN $\zeta^{CURRENT}$
  \end{algorithmic}
\end{algorithm}

Finally, we define a function $P(\zeta, R): \mathbb{N}^{|\mtc{V}|}\times \mathbb{N}^{|\mtc{I}|\times|\mtc{I}|}\to\mathbb{N}^{|\mtc{V}|}\times \mathbb{N}^{|\mtc{I}|\times|\mtc{I}|}$ that randomly re-assigns $R$\% of the positions and $R$\% of the fees. The entire Iterated Local Search is sketched in \Cref{alg:ILS}.

\begin{algorithm}[h]
  \caption{Iterated Local Search}
  \label{alg:ILS}
  \begin{algorithmic}[1]
    \STATE INPUT: \texttt{MAXRESTARTSWITHOUTIMPROVEWMENT}, \texttt{TIMELIMIT}, $R$
    \STATE \texttt{\#RESTARTSWITHOUTIMPROVEWMENT}$\gets 0$
    \STATE $\zeta^{CURRENT}\gets$ random solution
    \STATE $\zeta^{BEST}\gets\zeta^{CURRENT}$
    \WHILE{\texttt{\#RESTARTSWITHOUTIMPROVEWMENT}$\leq$\texttt{MAXRESTARTSWITHOUTIMPROVEWMENT} AND \texttt{ELAPSEDTIME}$\leq$\texttt{TIMELIMIT}}
    \STATE $\zeta^{N}\gets LS^\Pi(LS^\Lambda(\zeta^{CURRENT}))$
    \IF{$z(\zeta^N)<z(\zeta^{CURRENT})$}
    \STATE \texttt{\#RESTARTSWITHOUTIMPROVEWMENT}$\gets$ \texttt{\#RESTARTSWITHOUTIMPROVEWMENT}$+1$
    \ELSE
      \STATE \texttt{\#RESTARTSWITHOUTIMPROVEWMENT}$\gets 0$
      \ENDIF
      \IF{$z(\zeta^N)>z(\zeta^{BEST})$}
      \STATE $\zeta^{BEST}\gets \zeta^N$
      \ENDIF
      \STATE $\zeta^{CURRENT}\gets P(\zeta^N,R)$
      \ENDWHILE
    \RETURN $\zeta^{BEST}$
  \end{algorithmic}
\end{algorithm}
In \Cref{alg:ILS} we set \texttt{MAXRESTARTSWITHOUTIMPROVEWMENT} to $3$, \texttt{TIMELIMIT} to $1800$ seconds, and $R$ to $30$\%.

\Cref{tab:ILS02,tab:ILS08} report the optimality gap and solution time of the ILS compared to that of the L-Shaped method on all instances with identical customer profiles.
The optimality gap of the ILS is calculated using the bound delivered by the L-Shaped method as
$$\texttt{gap}=\frac{|\texttt{ILSOBJECTIVE}-\texttt{LSBESTBOUND}|}{|ILSOBJECTIVE|+10^{-10}}$$
The tables show that, for the smallest instances, the performance of the ILS is comparable to that of the L-Shaped method. On a number of instances (e.g., with $|\mtc{V}|=50$) the ILS even delivers better primal solutions than the L-Shaped method. Nevertheless, as the size of the instances increases the performance of the ILS drops.

\begin{longtable}{rrrr|rrrr}
  \caption{Comparison of optimality gap and solution time obtained with ILS and the LS method on the instances with identical customer profiles and with $\alpha^V=0.2$.}\label{tab:ILS02}\\
  \toprule
  $|\mtc{V}|$ & $|\mtc{K}|$  &  $\alpha^{FROM}$ &  $\alpha^{TO}$  & \texttt{gapLS} &    \texttt{gapILS} &   \texttt{tLS} & \texttt{tILS}  \\
\midrule
 50 & 200 &       0.2 &     0.2 &  0.6967 &   1.2615 & 1800.20 &  1801.00 \\
 50 & 400 &       0.2 &     0.2 & 17.2850 &  14.4768 & 1800.52 &  1800.13 \\
 50 & 600 &       0.2 &     0.2 & 30.5015 &  21.5516 & 1809.94 &  1800.05 \\
100 & 200 &       0.2 &     0.2 &  0.0975 &   0.3102 & 1800.98 &  1804.39 \\
100 & 400 &       0.2 &     0.2 &  1.9480 &   3.9251 & 1805.67 &  1802.78 \\
100 & 600 &       0.2 &     0.2 & 11.6208 &  44.7249 & 1801.95 &  1811.11 \\
200 & 400 &       0.2 &     0.2 &  0.0000 & 276.9414 &  139.41 &  1810.90 \\
200 & 600 &       0.2 &     0.2 &  0.0777 & 456.4824 & 1800.78 &  1812.87 \\
\midrule
  50 & 200 &       0.2 &     0.8 &  0.3381 &   0.7718 & 1800.63 &  1800.57 \\
 50 & 400 &       0.2 &     0.8 &  9.6397 &   8.5835 & 1801.21 &  1801.20 \\
 50 & 600 &       0.2 &     0.8 & 19.5455 &  15.2189 & 1800.03 &  1801.78 \\
100 & 200 &       0.2 &     0.8 &  0.2815 &   0.2815 & 1800.13 &  1801.84 \\
100 & 400 &       0.2 &     0.8 &  0.7148 &   3.6945 & 1801.47 &  1801.23 \\
100 & 600 &       0.2 &     0.8 &  9.8896 &  38.4997 & 1812.93 &  1807.95 \\
200 & 400 &       0.2 &     0.8 &  0.0000 & 258.6824 &  115.74 &  1803.26 \\
200 & 600 &       0.2 &     0.8 &  0.4710 & 421.6527 & 1814.99 &  1810.50 \\
\midrule
  50 & 200 &       0.8 &     0.2 &  1.9850 &   3.9133 & 1800.13 &  1800.40 \\
 50 & 400 &       0.8 &     0.2 & 21.8829 &  17.4800 & 1800.02 &  1801.26 \\
 50 & 600 &       0.8 &     0.2 & 44.9080 &  33.0970 & 1803.83 &  1801.01 \\
100 & 200 &       0.8 &     0.2 &  0.5238 &   1.3703 & 1801.09 &  1801.37 \\
100 & 400 &       0.8 &     0.2 &  6.1743 &   8.1636 & 1806.48 &  1800.58 \\
100 & 600 &       0.8 &     0.2 & 23.6427 &  52.5532 & 1816.20 &  1806.71 \\
200 & 400 &       0.8 &     0.2 &  0.3072 & 260.9089 & 1804.22 &  1817.48 \\
  200 & 600 &       0.8 &     0.2 &  6.7474 & 386.4842 & 1865.22 &  1801.24 \\
  \midrule
 50 & 200 &       0.8 &     0.8 &  0.8509 &   0.9891 & 1800.15 &  1801.00 \\
 50 & 400 &       0.8 &     0.8 & 15.4867 &  12.6303 & 1801.46 &  1801.19 \\
 50 & 600 &       0.8 &     0.8 & 40.1722 &  28.5881 & 1807.37 &  1800.41 \\
100 & 200 &       0.8 &     0.8 &  0.7932 &   1.4929 & 1802.32 &  1800.77 \\
100 & 400 &       0.8 &     0.8 &  4.4619 &   5.6145 & 1803.98 &  1807.19 \\
100 & 600 &       0.8 &     0.8 & 19.1824 &  45.0015 & 1823.21 &  1806.86 \\
  200 & 400 &       0.8 &     0.8 &  0.0963 & 242.2407 & 1800.78 &  1806.16 \\
  200 & 600 &       0.8 &     0.8 &  5.6194 & 307.5339 & 1823.23 &  1811.13 \\
  \bottomrule
\end{longtable}

\begin{longtable}{rrrr|rrrr}
  \caption{Comparison of optimality gap and solution time obtained with ILS and the LS method on the instances with identical customer profiles and with $\alpha^V=0.8$.}\label{tab:ILS08}\\
  \toprule
  $|\mtc{V}|$ & $|\mtc{K}|$  &  $\alpha^{FROM}$ &  $\alpha^{TO}$  & \texttt{gapLS} &    \texttt{gapILS} &   \texttt{tLS} & \texttt{tILS}  \\
  \midrule
   50 & 200 &       0.2 &     0.2 &  0.4613 &   1.2477 & 1800.37 &  1800.88 \\
 50 & 400 &       0.2 &     0.2 & 11.6098 &  10.5114 & 1802.87 &  1800.99 \\
 50 & 600 &       0.2 &     0.2 & 30.3469 &  23.8963 & 1807.97 &  1800.32 \\
100 & 200 &       0.2 &     0.2 &  0.1113 &   0.4159 & 1800.27 &  1801.25 \\
100 & 400 &       0.2 &     0.2 &  4.0199 &   6.9075 & 1800.15 &  1806.85 \\
100 & 600 &       0.2 &     0.2 & 14.6664 &  42.8427 & 1805.19 &  1813.64 \\
200 & 400 &       0.2 &     0.2 &  0.9991 & 235.4502 & 1810.37 &  1814.95 \\
200 & 600 &       0.2 &     0.2 &  4.5424 & 347.5242 & 1889.13 &  1817.40 \\
\midrule
  50 & 200 &       0.2 &     0.8 &  0.8298 &   1.2781 & 1800.11 &  1800.27 \\
 50 & 400 &       0.2 &     0.8 & 11.8754 &   8.8758 & 1801.42 &  1800.32 \\
 50 & 600 &       0.2 &     0.8 & 24.6840 &  21.8355 & 1800.96 &  1801.07 \\
100 & 200 &       0.2 &     0.8 &  0.0291 &   1.0756 & 1800.47 &  1800.33 \\
100 & 400 &       0.2 &     0.8 &  2.1136 &   5.9566 & 1805.67 &  1800.81 \\
100 & 600 &       0.2 &     0.8 & 12.3031 &  36.9526 & 1817.99 &  1803.40 \\
200 & 400 &       0.2 &     0.8 &  0.2123 & 219.8087 & 1800.23 &  1811.94 \\
200 & 600 &       0.2 &     0.8 &  3.7770 & 332.4861 & 1833.46 &  1805.89 \\
\midrule
  50 & 200 &       0.8 &     0.2 &  0.0389 &   0.3443 & 1800.08 &  1800.10 \\
 50 & 400 &       0.8 &     0.2 & 13.3947 &  12.2721 & 1800.04 &  1801.02 \\
 50 & 600 &       0.8 &     0.2 & 34.8704 &  26.9711 & 1803.88 &  1802.90 \\
100 & 200 &       0.8 &     0.2 &  0.0000 &   0.0000 &   26.88 &  1800.97 \\
100 & 400 &       0.8 &     0.2 &  0.2335 &   1.6940 & 1801.63 &  1802.90 \\
100 & 600 &       0.8 &     0.2 &  9.1181 &  55.6383 & 1801.40 &  1810.57 \\
200 & 400 &       0.8 &     0.2 &  0.0020 & 232.7449 &  205.30 &  1806.80 \\
200 & 600 &       0.8 &     0.2 &  0.0034 & 331.4251 &  786.63 &  1816.44 \\
\midrule
  50 & 200 &       0.8 &     0.8 &  0.0083 &   0.2765 &   24.39 &  1800.52 \\
 50 & 400 &       0.8 &     0.8 &  9.8202 &   8.4977 & 1801.95 &  1802.57 \\
 50 & 600 &       0.8 &     0.8 & 30.9029 &  21.7163 & 1808.71 &  1800.03 \\
100 & 200 &       0.8 &     0.8 &  0.0000 &   0.0000 &   21.84 &  1800.08 \\
100 & 400 &       0.8 &     0.8 &  0.1667 &   1.3339 & 1801.88 &  1805.07 \\
100 & 600 &       0.8 &     0.8 &  6.1475 &  48.0809 & 1805.95 &  1809.60 \\
200 & 400 &       0.8 &     0.8 &  0.0000 & 224.0099 &  122.34 &  1800.48 \\
200 & 600 &       0.8 &     0.8 &  0.0294 & 285.5657 & 1802.96 &  1815.21 \\
  \bottomrule
\end{longtable}

\bibliographystyle{abbrv} 
\bibliography{cs} 

\end{document}